\documentclass[11pt,reqno,natbib]{amsproc} 
	\allowdisplaybreaks
	\usepackage{graphicx}
	\usepackage{fullpage}
	\usepackage{mathptmx}
	\usepackage{enumerate}
	\usepackage{array}
	\usepackage{float}
	\floatstyle{plaintop}
	\restylefloat{table}
	\usepackage{verbatim}
	\numberwithin{equation}{section}
	\usepackage{natbib} 
	\usepackage{mathtools}
	\graphicspath{{figures/}}
	\usepackage{subfig}
	\usepackage{color}
	\usepackage{tcolorbox}
	\usepackage{booktabs}
	\usepackage{caption,setspace}
	\usepackage[debug=false, colorlinks=true, pdfstartview=FitV, 
	linkcolor=blue, citecolor=blue, urlcolor=blue]{hyperref}

\usepackage{morefloats}

\newlength{\drop}
\definecolor{amethyst}{rgb}{0.6, 0.4, 0.8}
\definecolor{burgundy}{rgb}{0.5, 0.0, 0.13}

\begin{document}
	
	\title{A global sensitivity analysis and reduced order models 
	for hydraulically-fractured horizontal wells}
	
	\author{A.~Rezaei \and
	K.~B.~Nakshatrala \and
	F.~Siddiqui \and
	B.~Dindoruk \and
	M.~Soliman \\
	{\small University of Houston, Houston, Texas, USA 77204--4003.}
	}

\begin{titlepage}
  \drop=0.1\textheight
  \centering
  \vspace*{\baselineskip}
  \rule{\textwidth}{1.6pt}\vspace*{-\baselineskip}\vspace*{2pt}
  \rule{\textwidth}{0.2pt}\\[\baselineskip]
  {\LARGE \textbf{\color{burgundy}
  A global sensitivity analysis and reduced order models 
    \\[0.3\baselineskip] 
    for hydraulically-fractured horizontal wells}}\\[0.2\baselineskip]
    \rule{\textwidth}{0.4pt}\vspace*{-\baselineskip}\vspace{3.2pt}
    \rule{\textwidth}{1.6pt}\\[0.5\baselineskip]
    \scshape
    An e-print of the paper will be made 
    available on arXiv. \par
    \vspace*{0.3\baselineskip}
    Authored by \\[0.5\baselineskip]
	{\Large Ali Rezaei}, {\itshape Postdoctoral Research Associate\par}
	{\itshape Department of Petroleum Engineering \\
	University of Houston, Houston, Texas 77204--4003 \\
		\textbf{phone:} +1-651-245-7918, 
		\textbf{e-mail:} arezaei2@uh.edu\par}
		\vspace*{0.75\baselineskip}
	 {\Large Kalyana B.~Nakshatrala}, {\itshape Associate Professor\par}
    {\itshape Department of Civil \& Environmental Engineering \\
    University of Houston, Houston, Texas 77204--4003. \\ 
    \textbf{phone:} +1-713-743-4418, \textbf{e-mail:} knakshatrala@uh.edu \\
    \textbf{website:} \url{http://www.cive.uh.edu/faculty/nakshatrala}\par}
     \vspace*{0.75\baselineskip}
     {\Large Fahd Siddiqui}, {\itshape Postdoctoral Research Associate\par}
     {\itshape Department of Petroleum Engineering \\
	University of Houston, Houston, Texas 77204--4003 \\
	\textbf{phone:} +1-713-743-6103, \textbf{e-mail:} fsiddiq6@central.uh.edu \par}
     \vspace*{0.75\baselineskip}
      {\Large Birol Dindoruk}, {\itshape Professor\par}
      {\itshape Department of Petroleum Engineering\\
 University of Houston \& Shell International Exploration and Production Inc., 
		Houston, Texas 77204--4003 \\
		\textbf{phone:} +1-832-875-9092, \textbf{e-mail:} birol.dindoruk@shell.com 
		\par}
     \vspace*{0.75\baselineskip}
      {\Large Mohammed Soliman},  {\itshape Professor\par} 
      {\itshape Department Chair and William C. Miller Endowed Chair 
		Professor \\
		Department of Petroleum Engineering \\
		University of Houston, Houston, Texas 77204--4003 \\
		\textbf{phone:} +1-713-743-6103, \textbf{e-mail:} msoliman@central.uh.edu\par}
    \vfill
        {\scshape 2018}\\
        {\small Computational and Applied Mechanics Laboratory}\par
\end{titlepage}

	\begin{abstract}
	Several factors affect the performance and stimulation design of 
	hydraulically-fractured wells. Moreover, the dominant factors 
	vary for different quantities of interest, and vary based on the spatial 
	location and with the time of interest. Thus, it will be beneficial 
	if there is a systematic procedure to identify the dominant factors 
	affecting the quantities of interest.
	To this end, we present a systematic global sensitivity analysis using the 
		Sobol method which can be utilized to rank the variables that affect two  
		quantity of interests -- pore pressure depletion and stress 
		change -- around a hydraulically-fractured horizontal well 
		based on their degree of importance. These variables 
		include rock properties and stimulation design variables. 
		A fully-coupled poroelastic hydraulic fracture model is used 
		to account for pore pressure and stress changes due to 
		production. 
		To ease the computational cost of a simulator, we also 
		provide reduced order models (ROMs), which can be 
		used to replace the complex numerical model with 
		a rather simple analytical model, for calculating the pore pressure 
		and stresses at different locations around hydraulic fractures.
		The two main reason for choosing the Sobol method are 
		that it can capture the individual and interaction effects of 
		input variables on the variance of outputs (which is 
		not the case with local sensitivity analysis techniques). 
		It also furnishes a systematic procedure with strong 
		mathematical underpinning to generate ROMs 
		for various quantities of interests for a given mathematical 
		model and for a given set of input variables.
		The main findings of this research are: (i) mobility, production 
		pressure, and fracture half-length are the main contributors 
		to the changes in the quantities of interest. The percentage 
		of the contribution of each	parameter depends on 
		the location with respect to pre-existing hydraulic fractures 
		and the quantity of interest.   
		(ii) As the time progresses, the effect of mobility 
		decreases and the effect of production pressure increases. 
		(iii) These two variables are also dominant for horizontal 
		stresses at large distances from hydraulic fractures. 
		(iv) At zones close to hydraulic fracture tips or inside 
		the spacing area, other parameters such as fracture 
		spacing and half-length are the dominant factors that 
		affect the minimum horizontal stress. 
		The results of this study will provide useful guidelines for the 
		stimulation design of legacy wells and secondary operations 
		such as refracturing and infill drilling.  
	\end{abstract}

	
	\maketitle

	\section{INTRODUCTION AND MOTIVATION}
	\label{Sec:Sobol_intro}
	Horizontal wells that are drilled in unconventional reservoirs 
	typically 
	show 
	a decline in their initial production rate (IPR) over time. 
	Multiple 
	hydraulic
	fractures are usually placed in these wells towards
	the completion stage to enable the production of hydrocarbons. 
	These 
	fractures create high permeability
	conduits which allow the flow of hydrocarbons from the
	rock matrix to the well. In some of these wells, production
	rate drops below the non-economical threshold which places
	these wells at the risk of abandonment.
	The drop in the production rate is mainly due to the small drainage 
	area of 
	these low permeability reservoirs that is limited to the inner 
	reservoir 
	between the fractures \citep{ozkan2009comparison}. It can also have 
	other 
	sources such as proppant degradation and un-successful initial 
	stimulation. 
	Possible ways of increasing the production from these reservoirs 
	are to 
	re-fracture the horizontal well after the occurrence of the 
	production 
	decline, 
	or drill an infill well parallel to the current well and stimulate 
	it. 
	These 
	methods enable production from a bypassed, or intact area of the 
	reservoir. 
	A 
	problem that exists while doing these practices is the 
	redistribution of 
	stresses in the depleted area of the reservoir that might be in the 
	vicinity 
	of 
	a newly created hydraulic fracture. This redistribution of stress 
	causes 
	the 
	new fractures to behave differently than the initial fractures and 
	sometimes 
	failing in the re-stimulation attempts. Therefore, several factors 
	such as 
	extent and severity of the pore pressure depletion have to be 
	considered 
	while 
	designing a refracturing, or an infill fracturing process.  
	
	Hydraulic fracturing design plays a critical role in the success of 
	any 
	refracturing and infill well fracturing process. Several factors 
	such as 
	the 
	state of in-situ stresses, rock geomechanical properties, operation 
	variables 
	(e.g., pump rate, proppant concentration) should be considered in 
	such a 
	design. The host medium with hydraulic fractures can be considered 
	as a 
	fully-saturated poroelastic rock. Hence, in order to properly study 
	the 
	process of placing hydraulic fractures into
	the formation that contains pre-existing fractures with a depleted 
	area in 
	their vicinity,
	the strong coupling between pore pressure and
	rock deformation should be taken into account.
	It is demonstrated that the pore pressure change
	(which is caused by production or injection)
	redistributes the stress state of the rock
	in the vicinity of hydraulic fractures
	\citep{berchenko1997deviation, roussel2012role, safari2015fracture, 
	Rezaei_URTeC2017, Rezaei_ATCE2017, Rezaei_Journal2018}. 
	The stress redistribution affects further activities such as 
	refracturing 
	and 
	infill well fracturing by affecting the preferred propagation 
	direction of 
	fractures either in the same well or an off-set well. 
	Thus, it is crucial to understand the main variables that 
	contribute to the 
	stress redistribution.

	Sensitivity 
	analysis (SA) 
	is a method for quantifying the importance of each model input 
	parameter 
	on the value of a model output parameter. 
	This method may be used to identify the key input parameters whose 
	variance 
	affects the output parameters the most. 
	Moreover, it can be used to build a computationally faster model 
	than the 
	original model \citep{welch1992screening}. 
	Depending on the application, many methods have been introduced to 
	perform 
	such an analysis \citep[e.g.,][]{hill2006effective,sobol2001global}. 
	A review on the recent advances on sensitivity analysis techniques 
	may be 
	found in 
	\citep{iooss2015review,pianosi2016sensitivity,borgonovo2016sensitivity}.
	
	Generally, these methods may be categorized into two subsets, 
	namely local 
	and global sensitivity analysis \citep{saltelli2004sensitivity, 
	saltelli2008global}. 
	Global SA is a method to study the effect of the entire input 
	parameters on 
	the output parameters uncertainty,
	whereas in local SA the focus is on the output parameters 
	themselves rather 
	than their uncertainties. 
	This method can be categorized into four sub-categories: 
	regression-based, screening based, variance-based, and meta-model 
	sensitivity 
	analysis \citep{tian2013review}. 
	\cite{sobol1993sensitivity, sobol2001global} developed a global SA 
	method 
	for 
	calculating the individual input variable influences on the output 
	of a 
	complicated mathematical model. This method is used in this study 
	to 
	analyze 
	and rank the influencing parameters that affect the performance of 
	a 
	refracturing or infill well fracturing.

	 Several parameters affect the changes in pore pressure 
	around 
	horizontal wells. These factors include rock geomechanical 
	properties, 
	operational variables such as production rate (or pressure), HF 
	design 
	parameters such as spacing and half lengths of the pre-existing 
	hydraulic 
	fractures, well spacing (in the case of infill well fracturing), 
	and 
	reservoir 
	in-situ properties such as initial pore pressure. Identifying the 
	parameters 
	that have the greatest impact on the pore pressure depletion 
	extension and 
	magnitude (subsequently principal stresses magnitude and direction) 
	helps 
	to 
	make better decisions about time and design of refracturing and 
	infill well 
	fracturing. Because of the large uncertainty in reservoir rocks 
	properties, 
	sensitivity analysis is being used repeatedly in oil and gas 
	industry for 
	purposes such as matching production history (rate or pressure) 
	\citep{oliver2011recent}, optimizing operation parameters 
	\citep{yu2014sensitivity}, and forecast analysis 
	\citep{nashawi2010forecasting}. \cite{verde2015global} used global 
	sensitivity 
	analysis to investigate the effect of shear modulus, 
	Poisson's ratio, normal joint stiffness, and minimum horizontal 
	stress on 
	the 
	fluid pressure in an injector well and a producer well. They 
	concluded that 
	shear modulus, normal joint stiffness, minimum horizontal stress, 
	and 
	Poisson's 
	ratio have the greatest sensitivity indices respectively. 
	\cite{dai2014efficient} used a global sensitivity analysis based on 
	polynomial 
	chaos expansions (PCEs) proposed a global SA based on the uncertain 
	parameters 
	that were used in a reservoir simulator. \cite{yu2014sensitivity} 
	performed 
	a 
	local sensitivity analysis on shale gas to optimize hydraulic 
	fracture 
	half-length and spacing. \cite{westwood2017sensitivity} applied a 
	Monte-Carlo 
	approach to study the effect of flow rate, pumping time, and 
	differential 
	pressure on the distance of the fluid penetration, stimulated rock 
	volume 
	(SRV), and minimum distance to avoid reactivation of the fault.

	In this 
	paper, our 
	aim is to first show the sensitivity of pore pressure and stresses 
	changes 
	to 
	rock properties, 
	operation variables, and design parameters. 
	Then, we use a global sensitivity scheme to study the uncertainty 
	that is 
	involved in the geomechanical and in-situ variables. Using this 
	approach, 
	parameters are indexed by their importance on the variation of pore 
	pressure 
	and stresses at an arbitrary point inside the rock. 
	It also has the advantage of capturing both individuals and 
	interaction 
	effects 
	of the parameters that are involved in the problem. This approach 
	helps 
	operators to select design parameters in a way to avoid the 
	occurrence of 
	problems such as stress reversal that may negatively affect any 
	refracturing 
	or 
	infill well fracturing. Finally, we use Sobol method to present a reduced 
	order model for points around hydraulically-fractured well at different times 
	from production.

	The rest of this paper is organized as follows.
	Section \ref{Sec:Sobol_TB} provides the necessary
	theoretical background, which includes the presentation of a 
	fully-coupled 
	poroelastic hydraulic fracture model and a brief description of the 
	Sobol 
	method. Uncertainty of the pore pressure and stresses concerning rock type
	variables is given in Section \ref{sec:example}.
	Global sensitivity analysis based on Sobol method is represented in Section 
	\ref{Sec:Global_analysis} to 
	index the 
	variables in order of their significance concerning the pore 
	pressure. In Section \ref{Sec:reduced_model}, using the dominant Sobol 
	indices, reduced order models for pore pressure and stress are developed at 
	different locations around hydraulic fractures.	Finally, conclusions are drawn
	in Section \ref{Sec:Sobol_CR}.
	
	\section{THEORETICAL BACKGROUND}
	\label{Sec:Sobol_TB}
	Our work hinges on poroelasticity, the
	displacement discontinuity method and
	the Sobol method. These ingredients
	are briefly described below for the
	benefit of the reader and for future
	referencing. 
	
	\subsection{Poroelasticity} 
	A poroelastic medium can be characterized by
	five independent material properties: the shear modulus $G$,
	the drained Poisson's ratio $ \nu $, the undrained
	Poisson's ratio $ \nu_u $, the Skempton's pore
	pressure coefficient $B$, and the permeability
	coefficient $\kappa$ 
	\citep{cleary1977fundamental,detournay1987poroelastic}, 
	which we refer to mobility in this paper.
	The Skempton's coefficient is defined as the ratio
	between the induced pore pressure and the variation
	of the confining pressure under the undrained condition, and the 
	permeability 
	coefficient is defined as the ratio between the rock permeability 
	$k$ and 
	the
	dynamic fluid viscosity $\mu$ (i.e., $\kappa= k/\mu$).
	Other parameters such as rock diffusion coefficient $c$ and Biot's 
	coefficient $\alpha$ can be derived from these independent 
	variables as 
	follows:
	\begin{equation}
	c = \frac{2\kappa B^2 G (1-\nu)(1+\nu_u)^2}{9(1-\nu_u)(\nu_u-\nu)} 
	\,.
	\end{equation}
	\begin{equation}
	\alpha = \frac{3(\nu_u-\nu)}{B(1-2\nu)(1+\nu_u)} = 1 - 
	\frac{K_m}{K_s} \,,
	\label{eq:alpha}
	\end{equation}
	where $K_s$ and $K_m$ are the solid and porous matrix bulk moduli, 
	respectively. 
	The response of a poroelastic medium is governed by
	four sets of equations, which are referred to as the
	\textit{field equations}. The four sets of equations
	are constitutive relations, force-equilibrium equations,
	Darcy's law, and continuity equation
	\citep{biot1941general,cleary1977fundamental}. 
	Constitutive equations relate stress, strain, and pore pressure. 
	Unlike 
	elastic media (which need one 
	constitutive relation), two constitutive
	relations are needed for poroelasticity.
	Of course, the field equations should 
	be augmented with appropriate boundary 
	and initial conditions. 
	
	A brief formulation of the generalized equations of poroelasticity 
	are in 
	order. We denote a spatial point by $\mathbf{x}$. The
	gradient and divergence operators with respect
	to $\mathbf{x}$ are, respectively, denoted by
	$\mathrm{grad}[\cdot]$ and $\mathrm{div}[\cdot]$. 
	The Laplacian differential operator is denoted 
	by $\Delta$. That is, $\Delta (\cdot) = 
	\mathrm{div}[\mathrm{grad}[\cdot]]$. 
	We denote the displacement field by $\mathbf{u}$ and also time of 
	an 
	arbitrary variable $a$ by $\dot{a}$. 
	We employ linearized strain, which is
	defined as follows:
	\begin{align}
	\mathbf{E}_{l} = \frac{1}{2} 
	\left(\mathrm{grad}[\mathbf{u}]
	+ \mathrm{grad}[\mathbf{u}]^{\mathrm{T}}\right)
	\end{align} 
	Note that the strain is a second-order 
	tensor. The volumetric strain is given 
	by $\mathrm{tr}[\mathbf{E}_{l}]$, where 
	$\mathrm{tr}[\cdot]$ denotes the 
	trace of a second-order tensor \citep{chadwick2012continuum}. 
	These equations for the case of plane strain
	quasi-static poroelasticity can be written
	as follows:
  \begin{align}
		\label{eq:navier}
		G \Delta (\mathbf{\mathbf{u}}) + 
		\left\{\frac{G}{1-2\nu_u}\right\}\mathrm{grad}
		[\mathrm{tr}[\mathbf{E}_{l}]] -\left\{ 
		\frac{2GB(1+\nu_u)}{3(1-2\nu_u)}\right\}\mathrm{grad}[
		\zeta] + \mathbf{f_b} = & \;0\\
		\label{eq:diff}
		\dot{\zeta}- c \; \Delta(\zeta)
		- \left\{\frac{\kappa 
		B(1+\nu_u)}{3(1-\nu_u)}\right\}\mathrm{div}[\mathbf{f_b}]
		+ \kappa \; \mathrm{div}[\mathbf{f_f}] - \mathbf{\gamma} = & \; 
		0
		\end{align}
	
	where $\zeta$ is the variation of fluid content 
	defined as the increment of fluid volume per unit volume of the 
	porous 
	medium 
	\citep{biot1941general}. Moreover, $\mathbf{f_b}$ and 
	$\mathbf{f_f}$ are 
	the 
	bulk and fluid body forces, respectively, $\gamma$ is the volume 
	rate of 
	injection from the fluid source. It should be noted that we employ 
	mechanics convention -- tensions are treated as positive and pore 
	pressure 
	is 
	positive in compression.
	Eq. \eqref{eq:navier} is referred to as the Navier equation of 
	poroelasticity 
	\citep{cleary1977fundamental,cheng2016poroelasticity}.
	The two main assumptions in the deriving Equations
	\eqref{eq:navier}--\eqref{eq:diff}
	are as follows:
	\begin{enumerate}[(i)]
		\item the fracturing fluid and the reservoir the fluid has the 
		same 
		rheology, and
		\item the deformation in the rock occurs in
		a quasi-static plane strain condition.
	\end{enumerate}
	Equations (\ref{eq:navier}) and (\ref{eq:diff})
	form the building blocks of the poroelastic
	displacement discontinuity method (DDM). 
	\cite{cleary1977fundamental} 
	presented fundamental solutions of point force and point fluid 
	source
	for the theory of poroelasticity, and 
	\citet{curran1987displacement} and 
	\cite{detournay1987poroelastic} developed displacement 
	discontinuity 
	solutions of a poroelastic medium using these equations.
	
	\subsection{The poroelastic displacement discontinuity method} 
	\label{sec_poroddm}
\begin{figure}[!htbp]
		\centering
		\includegraphics[width = 0.4\textwidth]{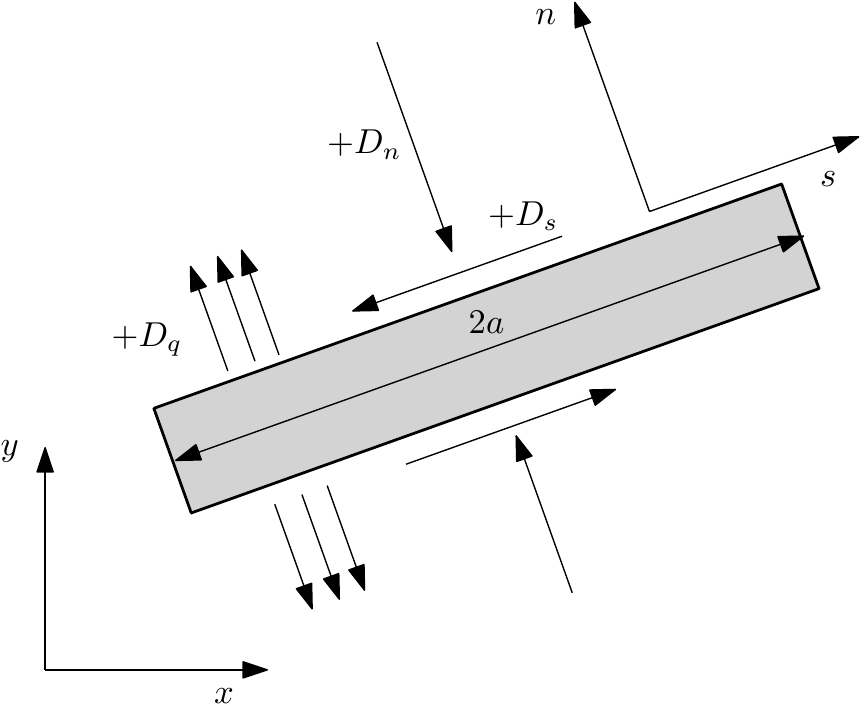}
		\caption{Definition of the shear, normal and flow 
		discontinuities on a 
		fracture element (the convention on the positive verses is 
		shown)} 
		\label{fig:DDM_sample} 
	\end{figure}

	In this section we describe the Poroelastic Displacement 
	Discontinuity 
	Method 
	(PDDM)
	that is used in this work. 
	This method belongs to the class of boundary element methods (BEM),
	which are suitable discretization for problems in which the ratio 
	of volume 
	to the surface is high. 
	Several authors used BEM
	for solving fracture mechanics problems 
	\citep{aliabadi1991boundary,aliabadi2002boundary,cruse2012boundary}.
	A 
	special indirect boundary element method for the case of a line 
	crack in an 
	infinite medium was developed by \cite{crouch1976solution}. 
	This method is based on considering the fracture as a line in 2D 
	(or a 
	surface in 3D)
	along which one defines quantities that take into account the 
	discontinuity 
	in displacements from one side of the fracture to the other. 
	\cite{liu2014revisit} explicitly showed that for problems involving 
	a 
	fracture, both DDM and BEM are equivalent.
	In its original formulation, 
	DDM is based on the fundamental solutions of a point source in an 
	infinite 
	linear elastic medium. 
	These fundamental solutions may be derived from dislocation theory 
	\citep{bobet2005stress}. On the other hand, in a poroelastic medium
	hydraulic fractures may be seen as a manifold 
	across which a discontinuity takes place 
	not only in the rock displacement 
	but also in the fluid flux. We define three \textit{discontinuity 
	fields} 
	with respect to a local coordinate system $ (s,n) $ (Figure 
	\ref{fig:DDM_sample}) as 
	\begin{subequations}
	\begin{equation}
	D_s(s,n,t) = \lim_{n \rightarrow 0-} u_s(s,n,t) - \lim_{n 
	\rightarrow 0+} 
	u_s(s,n,t) \,, \\
 \end{equation}
	\begin{equation}
	D_n(s,n,t) = \lim_{n \rightarrow 0-} u_n(s,n,t) - \lim_{n 
	\rightarrow 0+} 
	u_n(s,n,t) \,, \\
	\end{equation}
	\begin{equation}
	D_q(s,n,t) = \dfrac{q_{0}}{2 \, a \,} \,.
	\end{equation} 
	\end{subequations}
	Here, $D_s $ and $ D_n $ denote the shear and normal displacement 
	discontinuity fields, 
	and $ D_q $ is the flux discontinuity field.
	These fields physically represent the discontinuities in the 
	displacements 
	and the flow determined by a fracture.
	Also, $ u_s $ and $ u_n $ are the shear and normal components of 
	the 
	displacement field $\mathbf{u}$ in the local coordinate system
	and $ q_0 $ is the total flow injection, 
	The above definition implies that a fracture opening corresponds to 
	a 
	negative $D_n$ and counterclockwise movement of the fracture 
	surfaces gives 
	a 
	positive $D_s$.
	Also, the fluid flow is assumed to be negative if its direction is 
	different 
	from the chosen positive direction of total injection.
	Notice that the original displacement discontinuity method 
	developed by 
	\cite{crouch1976solution} does not have the $D_q$ term 
	since it was only developed using elastic solutions with no effect 
	of 
	poroelasticity.
	
	In order to solve for displacement and flux discontinuities, 
	we start from the integral equations relating such discontinuities 
	to 
	stresses and pore pressure in infinite domains 
	\citep{detournay1987poroelastic,vandamme_two}, which read
	for $ i,j=1,2$ (here, summation over repeated indices is used 
	except for 
	the 
	indices $ s,n,q$)
	\begin{subequations}
	\begin{align} \label{eq:int_sigs}
	\sigma_{ij}(\mathbf{x},t) = 
	& 
	\int_{0}^{t} \int_{\Gamma}
	Q_{ik}(\mathbf{\chi}) Q_{jl}(\mathbf{\chi}) 
	S_{s,kl}(\mathbf{x},\mathbf{\chi};t-\tau) D_{s}(\mathbf{\chi},\tau) 
	d\Gamma(\mathbf{\chi}) d\tau 
	\nonumber \\
	+ & \int_{0}^{t} \int_{\Gamma}
	Q_{ik}(\mathbf{\chi}) Q_{jl}(\mathbf{\chi}) 
	S_{n,kl}(\mathbf{x},\mathbf{\chi};t-\tau) D_{n}(\mathbf{\chi},\tau)
	d\Gamma(\mathbf{\chi}) d\tau 
	\nonumber	\\
	+ & \int_{0}^{t} \int_{\Gamma}
	Q_{ik}(\mathbf{\chi}) Q_{jl}(\mathbf{\chi}) 
	S_{q,kl}(\mathbf{x},\mathbf{\chi};t-\tau) D_{q}(\mathbf{\chi},\tau) 
	d\Gamma(\mathbf{\chi}) d\tau \,, 
	\end{align}
	\begin{align} 	\label{eq:int_pp}
	p(\mathbf{x},t) = 
	& 
	\int_{0}^{t}\int_{\Gamma} 
	P_{s}(\mathbf{x},\mathbf{\chi};t-\tau) D_{s}(\mathbf{\chi},\tau)
	d\Gamma(\mathbf{\chi}) d\tau 
	\nonumber \\
	+ &
	\int_{0}^{t}\int_{\Gamma} 
	P_{n}(\mathbf{x},\mathbf{\chi};t-\tau) D_{n}(\mathbf{\chi},\tau)
	d\Gamma(\mathbf{\chi}) d\tau 
	\nonumber \\
	+ &
	\int_{0}^{t}\int_{\Gamma} 
	P_{q}(\mathbf{x},\mathbf{\chi};t-\tau) D_{q}(\mathbf{\chi},\tau)
	d\Gamma(\mathbf{\chi}) d\tau \,.
	\end{align}
	\end{subequations}
	Equations \eqref{eq:int_sigs} and \eqref{eq:int_pp} are analytical 
	solutions
	over the plane, where the \textit{influence functions} 
	$ S_{s,kl} $, $ S_{n,kl} $, $ S_{q,kl} $,
	$ P_{s} $, $ P_{n} $, $ P_{q} $
	are given in \cite{vandamme_two} and 
	\cite{carvalho1991poroelastic}. The 
	influence functions give the solution of a point source or 
	displacement 
	discontinuity from influencing point $\mathbf{\chi}$ and at time 
	$\tau$ on 
	a 
	the influencing point $\mathbf{x}$ at time $t$. Matrix $ \mathbf{Q} 
	$ 
	represents the rotation 
	from the local crack coordinate system to the global coordinate 
	system.
	
	Equations \eqref{eq:int_sigs} and \eqref{eq:int_pp} form a set of 
	integral 
	equations for the unknown fields 
	$D_s$, $D_n$ and $D_q$ 
	with given fields $\sigma_{s}$, $\sigma_{n}$ and $p$. They show 
	that 
	stresses and pore pressures are obtained as a time integral 
	that contains an integral over the fracture $\Gamma$. The 
	discretization of 
	the integral equations using constant spatial and constant temporal 
	elements 
	may be written as 
	\begin{subequations} 
	\begin{align}\label{eq:final_sigs}
	\sum\limits_{\lambda=1}^N 
	A_{xx}^{\beta \lambda}D_{s}^{\lambda,h} +
	\sum\limits_{\lambda=1}^N 
	A_{xy}^{\beta \lambda}D_{n}^{\lambda,h} +
	\sum\limits_{\lambda=1}^N 
	A_{xq}^{\beta \lambda}D_{q}^{\lambda,h} =
	\nonumber \\
	\sigma_{s}^{h}(x^\beta,t) - 
	\sum\limits_{\eta=0}^{h-1}\sum\limits_{\lambda=1}^N
	\bigg(
	A_{xx}^{\beta \lambda,\eta}D_{s}^{\lambda,\eta} + 
	A_{xy}^{\beta \lambda,\eta}D_{n}^{\lambda,\eta} + 
	A_{xq}^{\beta \lambda,\eta}D_{q}^{\lambda,\eta} \bigg), 
	\end{align}
	\begin{align}
	\sum\limits_{\lambda=1}^N 
	A_{yx}^{\beta \lambda}D_{s}^{\lambda,h} +
	\sum\limits_{\lambda=1}^N 
	A_{yy}^{\beta \lambda}D_{n}^{\lambda,h} +
	\sum\limits_{\lambda=1}^N 
	A_{yq}^{\beta \lambda}D_{q}^{\lambda,h} =
	\nonumber \\
	\sigma_{n}^{h}(x^\beta,t) - 
	\sum\limits_{\eta=0}^{h-1}\sum\limits_{\lambda=1}^N
	\bigg(
	A_{yx}^{\beta \lambda,\eta}D_{s}^{\lambda,\eta} + 
	A_{yy}^{\beta \lambda,\eta}D_{n}^{\lambda,\eta} + 
	A_{yq}^{\beta \lambda,\eta}D_{q}^{\lambda,\eta} \bigg), 
	\end{align} 
	\begin{align}\label{eq:final_pp}
	\sum\limits_{\lambda=1}^N 
	A_{px}^{\beta \lambda}D_{s}^{\lambda,h} +
	\sum\limits_{\lambda=1}^N 
	A_{py}^{\beta \lambda}D_{n}^{\lambda,h} +
	\sum\limits_{\lambda=1}^N 
	A_{pq}^{\beta \lambda}D_{q}^{\lambda,h} =
	\nonumber \\ 
	p_p^{h}(x^\beta,t) - 
	\sum\limits_{\eta=0}^{h-1}\sum\limits_{\lambda=1}^N
	\bigg(
	A_{px}^{\beta \lambda,\eta}D_{s}^{\lambda,\eta} + 
	A_{py}^{\beta \lambda,\eta}D_{n}^{\lambda,\eta} + 
	A_{pq}^{\beta \lambda,\eta}D_{q}^{\lambda,\eta} \bigg). 
	\end{align}
	\end{subequations}
	In Equations \eqref{eq:final_sigs}--\eqref{eq:final_pp}, 
	$ A_{ij} $ are the coefficients relating the displacement 
	discontinuities 
	and 
	fluid sources 
	to shear stress, normal stress and pore pressure. For example, $ 
	A_{xx}^{\beta \lambda , \eta} $ is the shear stress that is induced 
	on the 
	observation point $\beta$ from a unit shear displacement 
	discontinuity at 
	the 
	source point $\lambda$ during time $\eta$, where $\eta$ is the time 
	between 
	occurrence of event at the source point and receiving the effect at 
	the 
	observation point. In general, 
	the fracturing fluid pressure is known at the boundary of a 
	hydraulic 
	fracturing problem (i.e. fracture surface), 
	and the fracture surface displacements and flow discontinuity 
	are the unknowns. 
	Therefore, Equations \eqref{eq:final_sigs}--\eqref{eq:final_pp} 
	form a set of $ 3\;N $ linear equations 
	which may be solved for $ 3\;N $ unknowns namely $\sigma_s$, 
	$\sigma_n$, 
	and 
	$p$ 
	for at each time step. Different approaches may be taken for 
	solving above 
	equations. The time marching scheme that is used in this study is 
	explained 
	in previous publications 
	\citep[e.g.,][]{brebbia2012boundary,Rezaei_Journal2018}. In every 
	iteration 
	of 
	the time-marching scheme, the time increments of $ 3\;N $ 
	discontinuity 
	variables are computed. After obtaining the discontinuity fields at 
	any 
	interested time step, Equation \eqref{eq:final_sigs} - 
	\eqref{eq:final_pp} 
	may be used to obtain stresses and pore pressure in any point of 
	the rock 
	body.

	\subsection{Sobol method} \label{sec:glob_var}
	
	Global sensitivity analysis has been used in many applications for 
	purposes 
	such as model verification and understanding, simplifying 
	(reduced-order 
	model), and characterizing the influence of input parameters on the 
	uncertainty of the output \citep[e.g.,][]{archer1997sensitivity, 
	makowski2006global, volkova2008global, lefebvre2010methodological, 
	auder2012screening}. Regression-based and variance-based methods 
	are the 
	main two classes of methods for global sensitivity analysis 
	\citep{arwade2010variance}. Sensitivity measure that characterizes the class 
	is the 
	main 
	cause of this distinction. \cite{sobol1993sensitivity, sobol2001global} 
	introduced a method 
	for global 
	sensitivity analysis that may be used for linear and nonlinear 
	models.
	This method is based on the measurement of the fractional 
	contribution of 
	the input parameters to the variance of the model output. To 
	explain how 
	Sobol technique works, let us assume that a mathematical model is 
	represented by function $f$ such that
	\begin{align}
		y = f(\mathbf{x})
		\label{eq:def}
	\end{align}
	where $\mathbf{x}$ is a set of input parameters on the 
	n-dimensional 
	hypercube such that:
	\begin{align*}
		\Omega^n \coloneqq \left\{\mathbf{x} | 0 \le x_i \le , i= 
		1,\cdots , 
		n \right\}.
	\end{align*}
	The ANOVA representation (abbreviated from Analysis of Variances) 
	of the 
	function $f$ may be written as
	\begin{align}
	\begin{aligned}
	f(x) = f_0 + \displaystyle \sum_{s = 1}^{n}\sum_{i_1 < \cdots < 
	i_s}^{n} 
	f_{i_1 \cdots 
	i_s}(x_{i_{1}}, \cdots , x_{i_{s}}) \,,\quad 1 \le i_1 < \cdots < 
	i_s \le 
	n \,,
	\label{eq:ANOVA}
	\end{aligned}
	\end{align}
	Equation \eqref{eq:ANOVA} may be rearranged to get a 
	series of increasing order Sobol' functions as follows
	\begin{equation}
	\begin{aligned}
	f(x) = f_0 + \sum_{i=1}^{n}f_i(x_i) + \sum_{i=1}^{n} 
	\sum_{j=i+1}^{n} 
	f_{ij}(x_i,x_j) + \cdots + f_{i \cdots n}(x_1,\cdots,x_n) .
	\label{eq:ANOVA_rearang}
	\end{aligned}
	\end{equation}
	For Equation \eqref{eq:ANOVA_rearang} to be true, the following 
	criteria 
	should be satisfied:
	\begin{enumerate}
		\item $f_0$ should be constant
		\item The integral of each member over its own variables should 
		be 
		zero\\
		\begin{align*}
			\int_{0}^{1} f_{i_1 \cdots i_s}(x_{i_{1}}, \dots , 
			x_{i_{s}}) dx_k 
			= 
			0 \hspace{1cm} \forall k = i_1,\dots, i_s 
		\end{align*}
		\item All of the members in Equation \eqref{eq:ANOVA_rearang} 
		are 
		orthogonal, meaning that if $(i_1,\dots, i_s) \ne 
		(j_1,\dots,j_t)$ then
		\begin{align*}
			\int_{\Omega^n} f_{i_1 \dots i_s} \hspace{.1cm} f_{j_1 
			\dots j_t} 
			d\mathbf{x} = 0
		\end{align*} 
	\end{enumerate}
	The individual terms in Equation \eqref{eq:ANOVA_rearang} may be 
	defined as 
	\citep{sobol1993sensitivity,sobol2001global}
	\begin{align}
	\begin{aligned}
	f_0 & = \int_{\Omega^n} f(\mathbf{x}) d\mathbf{x} \\
	f_i(x_i) & = \int_{\Omega^{n-1}} f(x_i,\mathbf{x}_{\sim i}) \; 
	d\mathbf{x}_{\sim i} - f_0 \\
	f_{ij}(x_i,x_j) & = \int_{\Omega^{n-2}} f(x_i,x_j,\mathbf{x}_{\sim 
	ij}) \; 
	d\mathbf{x}_{\sim ij} - f_0 
	- f_i(x_i) - f_j(x_j)
	\end{aligned}
	\label{integ_def}
	\end{align}
	where $\mathbf{x}_{\sim i}$ is the vector corresponding to all 
	variables 
	except $x_i$ in the input set $\mathbf{x}$, and $\mathbf{x}_{\sim 
	ij}$ is 
	the vector corresponding to all variables except $x_i$ and $x_j$in 
	the 
	input set $\mathbf{x}$. Assuming that $f(x)$ is 
	square integrable, total variance of $f$ is given by
	\begin{align}
	\begin{aligned}
	D = V[f] & = \int_{\Omega^n} f^2(\mathbf{x}) \; d\mathbf{x} - f_0^2 
	& = 
	\sum_{s \;= \;1}^{n}\sum_{i_1 < \cdots < i_s}^{n} f_{i_1 \cdots 
	i_s}^2(x_{i_1}, \dots ,x_{i_s}) \; dx_{i_1} \; \cdots \; dx_{i_s}
	\end{aligned}.
	\label{eq:variance}
	\end{align}
	Equation \eqref{eq:variance} can also be written in terms of the 
	partial 
	variances of $f$ as
	\begin{align}
		D & = \sum_{s=1}^{n}\sum_{i_1 < \dots < i_s}^{n}D_{i_1 \dots 
		i_s} = 
		\sum_{i = 1}^{n} D_i + \sum_{i=1}^{n} \sum_{j=i+1}^{n} D_{ij} + 
		\dots + D_{1 \dots n}
	\end{align}
	where $D_i, D_{ij}, \dots, D_{i \dots j}$ can be calculated by 
	integrating 
	the corresponding Sobol functions as follows
	\begin{align}
	\begin{aligned}
	D_i & = \int_{\Omega^1} f_i^2(x_i) \; dx_i \\
	D_{ij} & = \int_{\Omega_2} f_{ij}^2(x_i,x_j) \; dx_i \; dx_j \\
	& . \\
	& . \\
	& . \\
	D_{i_1 \; \cdots i_s} & = \int_{\Omega} f_{i_1 \cdots 
	i_s}^2(x_{i_1}, 
	\cdots, x_{i_s}) \; dx_{i_1} 
	\cdots 	dx_{i_s}
	\end{aligned}
	\label{eq:partial_variance}
	\end{align}
	Using these definitions, one can define Sobol indices that are the 
	ratio of 
	the partial variances to the total variance as
	\begin{align}
	\begin{aligned}
	S_i & = \dfrac{D_i}{D}\\
	S_{ij} & = \dfrac{D_{ij}}{D}\\
	& . \\
	& . \\
	& . \\
	S_{i_1 \dots i_s} & = \dfrac{D_{i_1 \dots i_s}}{D}\\
	\end{aligned}
	\label{eq:sob_div}
	\end{align}
	In this arrangement, greater indices mean a greater impact on the 
	variation 
	of the output parameter. It also should be noted that Sobol indices 
	are 
	non-negative indices that have the following property 
	\begin{align}
	\begin{aligned}
	\sum_{s=1}^{n}\sum_{i_1 < \dots < i_s}^{n}S_{i_1 \dots i_s} = 
	\sum_{i=1}^{n} S_i + \sum_{i = 1}^{n} \sum_{j = 1 + 1}^{n} S_{ij} + 
	\dots + 
	S_{1 \dots n} = 1.
	\end{aligned}
	\end{align}
	Using Sobol indices, one may perform an analysis to order input 
	variables according to their Sobol indices. Numerical examples of 
	such 
	analysis on a polynomial function $f$ may be found in 
	\citep{sobol2001global,saltelli2008global,arwade2010variance}. In 
	the next 
	section, global sensitivity of a simple mathematical model is 
	illustrated 
	to demonstrate the method.

	\subsection{Sobol method for complex functions} 
	\label{subsec:numerical integration}
	For the cases where the function $f$ is not a polynomials such as a 
	numerical 
	simulator, where analytical solution is not available, a Monte 
	Carlo 
	integration is required to perform the integrals that are required 
	by Sobol 
	analysis 
	\citep{sobol1993sensitivity,witarto2018global}. In this approach, 
	Sobol 
	functions $\bar{f}_{i \dots j}$, total variance $\bar{D}$, and 
	partial 
	variances $\bar{D}_{i \dots j}$ can be calculated as follows: 
	\begin{align}
	\bar{f}_0 = \dfrac{1}{N} \sum_{m = 1}^{N} f(x_m)
	\label{eq:f_0} 
	\end{align}
	\begin{align}
	\bar{D} = \dfrac{1}{N} \sum_{m=1}^{N} \; f^2(x_m) - \bar{f}_0^2
	\end{align}
	\begin{align}
	\bar{D}_i = \dfrac{1}{N} \sum_{m=1}^{N} \; f(x_m) \; 
	f(x_{im},\mathbf{x}_{\sim 
	im}^c) - \bar{f}_0^2
	\end{align}
	\begin{align}
	\bar{D_{ij}} = \dfrac{1}{N} \sum_{m=1}^{N}f(x_m) \; 
	f(x_{im},x_{jm}, 
	\mathbf{x}_{\sim 
	ijm}^c) - \bar{D}_i - \bar{D}_j - \bar{f}_0^2.
	\end{align}
	where m is the test number and $N$ is the sample size of the 
	inputs. The 
	bar sign is used to show that the term is numerically integrated. 
	After 
	calculating the $D_{i\cdots j}$, Equation \eqref{eq:sob_div} to 
	calculate 
	the Sobol indices. Once these numerical variables are calculated, 
	Sobol 
	indices 
	may can be obtained using Equation \eqref{eq:sob_div}. 
	
	One advantage of using Sobol analysis is that it can present a 
	complex 
	function 
	with a rather simplified equation. For this purpose, Sobol 
	functions are 
	chosen 
	to a certain degree of accuracy based on the magnitude of the their 
	Sobol 
	index. In Section \ref{Sec:reduced_model}, we present a reduced 
	order model 
	for 
	the 
	model that we described in Subsection \ref{sec_poroddm}. Before that, it is 
	beneficial to show how pore pressure and stresses change around hydraulic 
	fractures as a result of the change in the rock geomechanical properties.

	\section{UNCERTAINTY OF THE GEOMECHANICAL AND OPERATIONAL PARAMETERS} 
	\label{sec:example}
	
\begin{figure}[!h]
		\centering 
		\subfloat[3D view]{
			\includegraphics[width=0.5\textwidth]{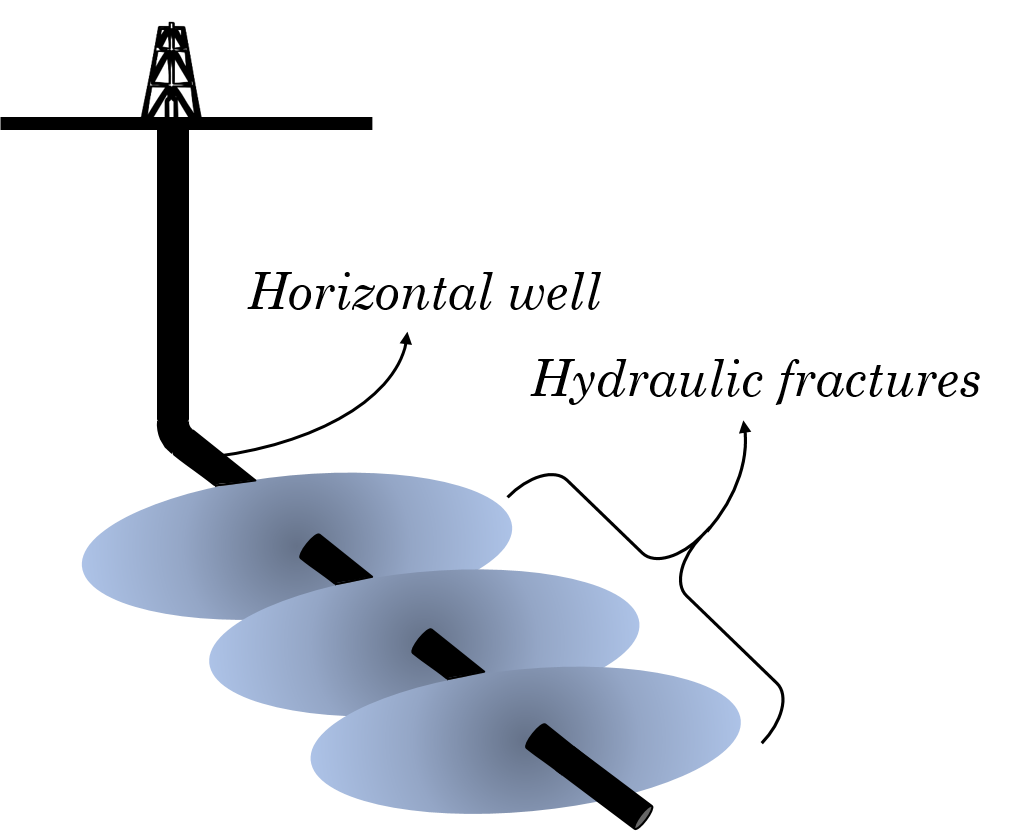}
			\label{fig:hf_3d}
		}
		\subfloat[Top view]{
			\includegraphics[width=0.5\textwidth]{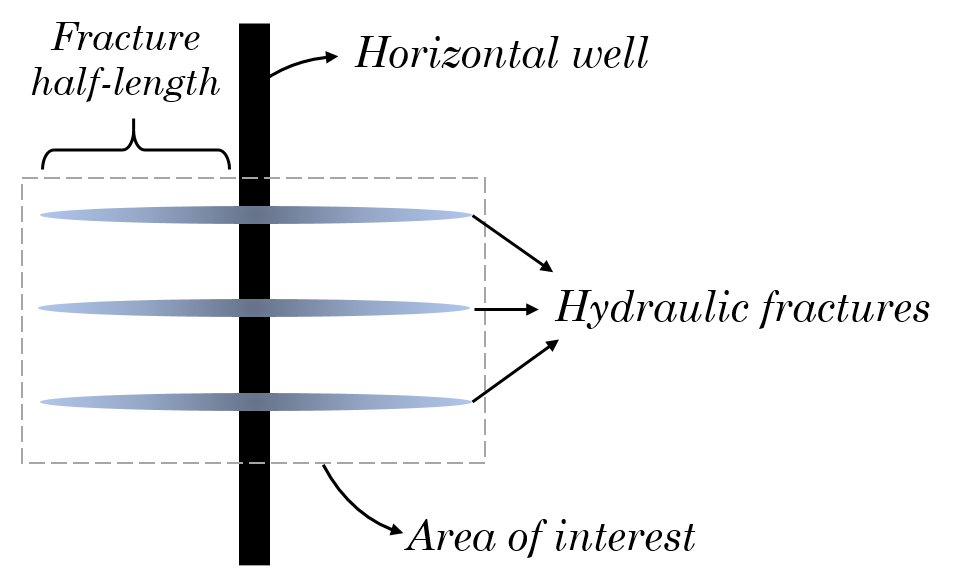}
			\label{fig:hf_2d}
		} 
		\caption{A schematic of hydraulic fracturing design in 
		horizontal wells}
		\label{fig:HF_horizontal}
	\end{figure}

	In this section, the uncertainty of the geomechanical problems is 
	demonstrated using an example of hydraulic fracturing in horizontal 
	wells.
	Hydraulic fracturing is widely used in oil and gas industry 
	as a stimulation technique to increase the hydrocarbon production 
	from 
	tight 
	formations.
	This technique is implemented in horizontal wells through multiple 
	stages.
	Usually, a section of the wellbore is isolated, perforated, and 
	pressurized 
	to open the path for the pressurized fluid to enter the perforation 
	and 
	form 
	a set of clusters \citep{soliman2006fracture, 
	soliman2016fracturing}. 
	Figure \ref{fig:HF_horizontal} shows a schematic of a typical 
	hydraulic 
	fracturing arrangement in horizontal wells. 
	In a normal faulting regime, the preferred direction of the 
	horizontal 
	wellbore 
	is the direction of the minimum horizontal stress for several 
	reasons.
	Firstly, the wellbore is more stable in this direction because 
	loading in 
	this situation is more isotropic.
	Secondly, since the preferred propagation direction of fractures is 
	maximum 
	compressional horizontal stress direction, 
	having the wellbore in this direction helps to have multiple 
	parallel 
	transverse hydraulic fractures (Figure \ref{fig:hf_2d}). 
	
	In the example presented in this section, the production from a 
	horizontal 
	well 
	containing two parallel hydraulic fractures is modeled using 
	different 
	types 
	of rocks.
	The aim is to show the effect of rock properties on the extent and 
	severity 
	of 
	the pore pressure depletion and its subsequent effect on horizontal 
	stresses, 
	even by using the same boundary conditions. 
	
\begin{figure}[!h]
		\centering
		\includegraphics[width = 0.5\textwidth]{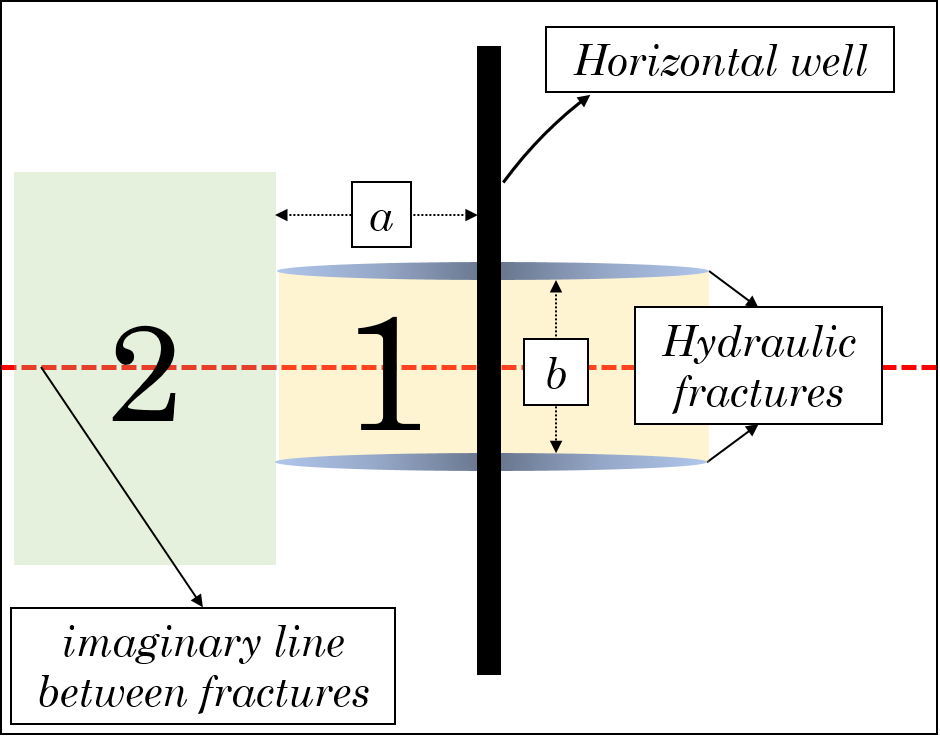}
		\caption{Geometry of the problem that is used for production} 
		\label{fig:cons_prod} 
	\end{figure}

 A set of different rocks are examined by their 
	response to 
	pore 
	pressure depletion. 
	The problem assumptions are as follows. 
	A horizontal well is drilled, and two parallel hydraulic fractures 
	are 
	created orthogonal to it.
	Both hydraulic fractures have reached their final length and effect 
	of 
	stress 
	shadowing on their geometry during their propagation is neglected. 
	Figure \ref{fig:cons_prod} shows the geometry of the problem for this 
	example. 
	The wellbore is put on continuous production for one month, one 
	year, and 
	five years.
	For handling production from fractures, 
	total stress loading mode inside the fracture similar to 
	\cite{mathias2010investigation} is used. 
	After each period of production, pore pressure, maximum horizontal 
	stress 
	($\sigma_H$), 
	minimum horizontal stress ($\sigma_h$), stress anisotropy 
	($\sigma_H-\sigma_h$) are calculated 
	along an imaginary line (dashed red line in Figure 
	\ref{fig:cons_prod}) on 
	the 
	middle point between fractures. 
	Moreover, the changes on these variables are going to be analyzed 
	in two 
	regions along that imaginary line. 
	These two regions are \textit{Region 1} colored by yellow and 
	\textit{Region 
	2} 
	colored by green that respectively represent the area between two 
	fractures 
	and 
	the zone in front of the fracture tips in Figure \ref{fig:cons_prod}. 
	
	In this example, both hydraulic fractures are assumed to have the 
	same 
	lengths, and their half-lengths are equal to $30\; m$.
	Fracture spacing (orthogonal distance between fractures) also is 
	equal to 
	$30\; m$.
	Moreover, maximum horizontal stress, minimum horizontal stress, and 
	reservoir 
	pore pressure are assumed to be $56.53\; MPa$, $55.15\; MPa$, and 
	$48.26\; 
	MPa$ respectively. 
	Furthermore, a constant pressure production is considered for the 
	entire 
	period of the production. 
	Five different rocks are chosen for this purpose. 
	The rocks have chosen somehow to represent a range of 
	low-permeability 
	reservoir rocks 
	in terms of geomechanical properties, although in reality some of 
	them are 
	not reservoir rocks.
	
	Table \ref{tab:rock_prop} presents the geomechanical properties of 
	the 
	rocks 
	that are used for the analysis in this section. 
	Relatively high permeability rocks (Weber Sandstone and Berea 
	Sandstone) 
	are 
	placed deliberately in the analysis 
	to show the difference in their pore pressure depletion compared to 
	ultra 
	tight rocks (Tennessee Marble, Charcoal Granite, and Haynesville 
	Shale).

	\begin{table}[!h]
		\centering
		\caption{Rock poroelastic properties used for the pore pressure 
		depletion 
		example. The table is 
			constructed using \cite{detournay1989poroelasticity} and 
			\cite{rice1976some}. Haynesville shale properties are 
			adopted from 
			\cite{chun2013thermo}.}
		\label{tab:rock_prop}  
		\begin{tabular}{llllllll}
			\hline\noalign{\smallskip}
			Rock type & $G, GPa$ & $\nu$ & $\nu_u$ & $B$ & $c, m^2/s$ & 
			$k, m^2$& 
			$\alpha$\\
			\noalign{\smallskip}\hline\noalign{\smallskip}
			Tennessee Marble & $24$ & $0.25$ & $0.31$ & 
			$0.51$ & $1.3 \times 10^{-5}$ & $1.0 \times 10^{-19}$ & 
			$0.19$ \\ 
			Haynesville Shale & $13.8$ & $0.22$ & $0.46$ & $0.91$ 
			& $4.5 \times 10^{-6}$ & $1.0 \times 10^{-19}$ & $0.96$ \\ 
			Berea Sandstone & $6$ & $0.20$ & $0.25$ & $0.62$ 
			& $4.5 \times 10^{-6}$ & $1.9 \times 10^{-13}$ & $0.96$ \\ 
			Charcoal Granite & 19 & 0.27 & $0.30$ & 0.55 & $7.0 \times 
			10^{-6}$ & 
			$1.0 \times 10^{-19}$ & 0.27 \\
			Weber Sandstone & 12 & 0.15 & 0.29 & 0.73 & $2.1 \times 
			10^{-2}$ & $1.0 
			\times 10^{-15}$ & 0.64 \\
			\noalign{\smallskip}\hline
		\end{tabular}
	\end{table}
	
\begin{figure}[!h]
		\centering
		
		\subfloat[Pore pressure]{
			\includegraphics[width=0.45\textwidth]{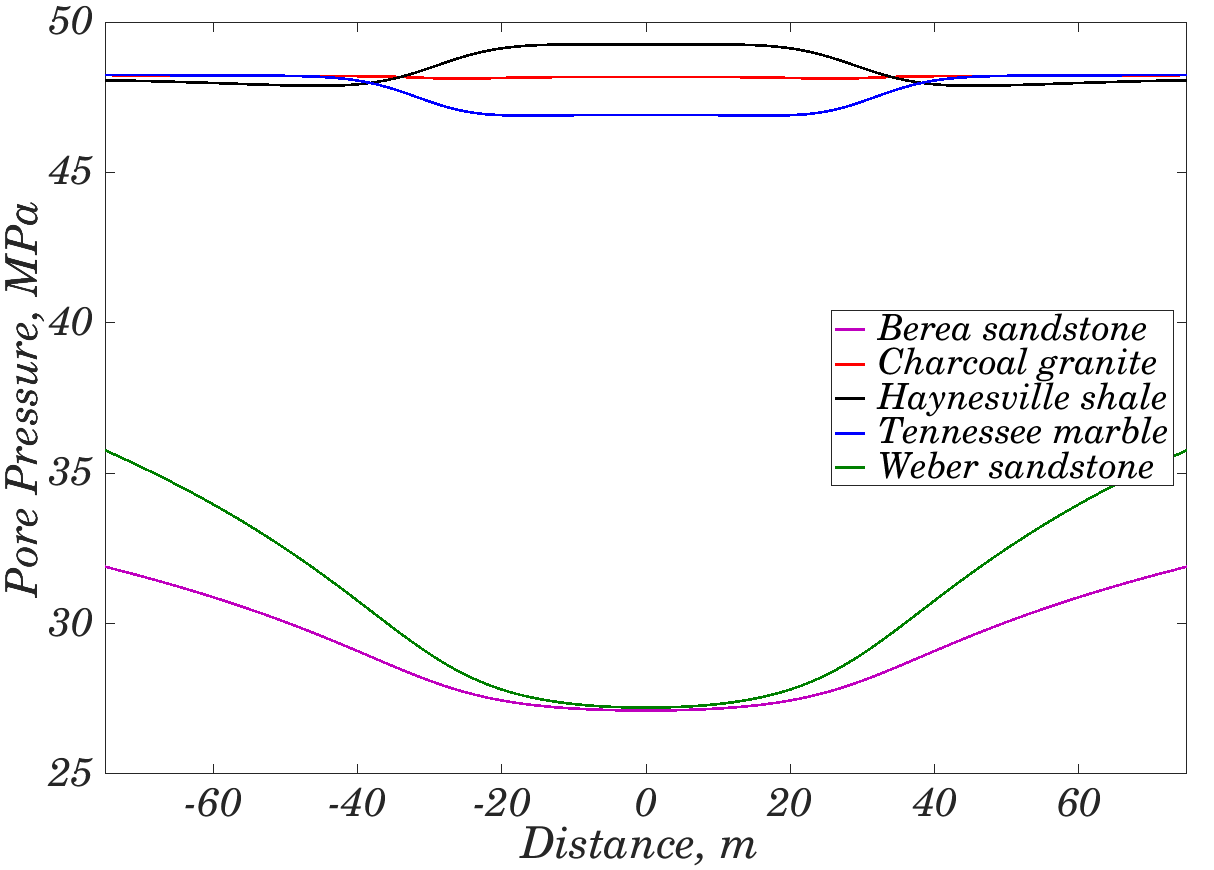}
			\label{fig:diff_pp_1mon}
		}
		\subfloat[Stress anisotropy]{
			\includegraphics[width=0.45\textwidth]{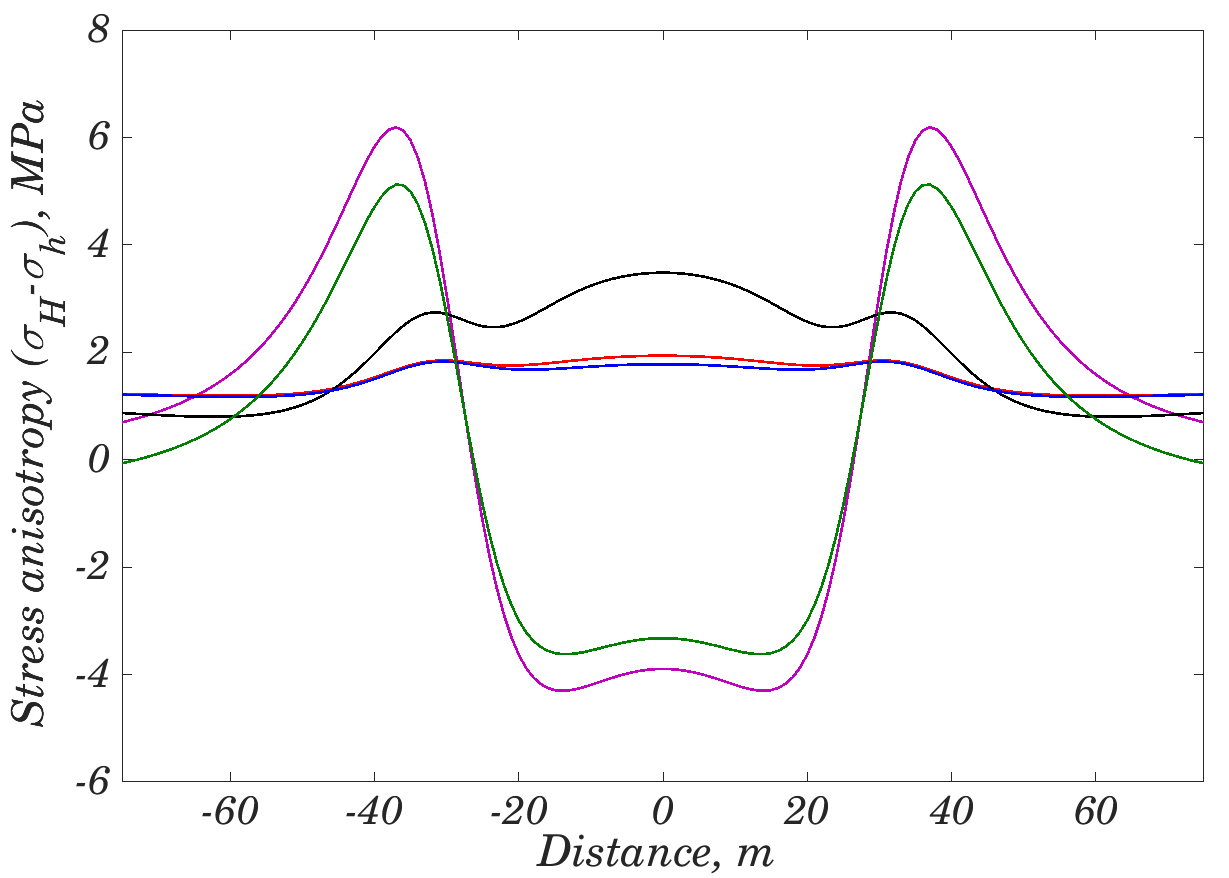}
			\label{fig:diff_aniso_1mon}
		} 
		\qquad
		\subfloat[$\sigma_{xx}$]{
			\includegraphics[width=0.45\textwidth]{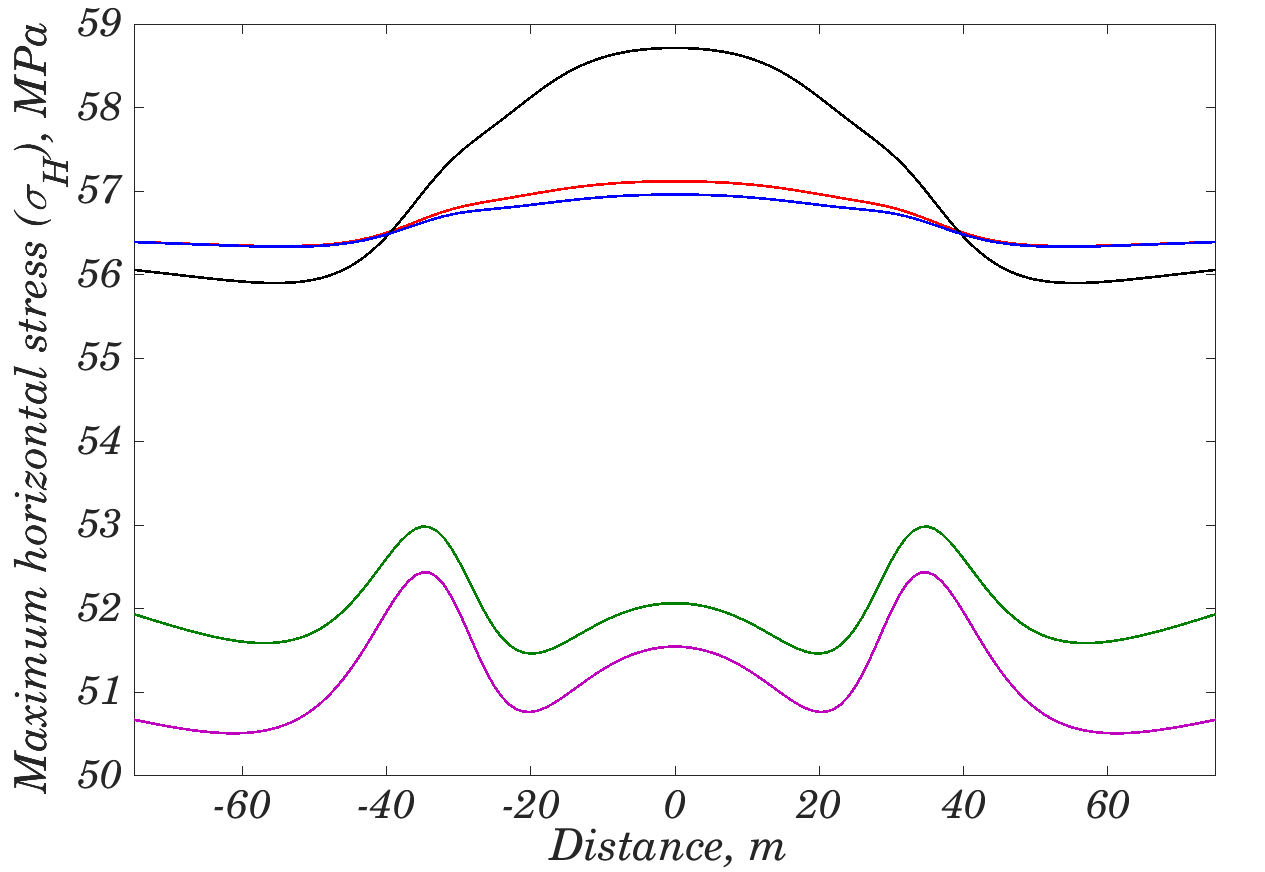}
			\label{fig:diff_sxx_1mon}
		}
		\subfloat[$\sigma_{yy}$]{
			\includegraphics[width=0.45\textwidth]{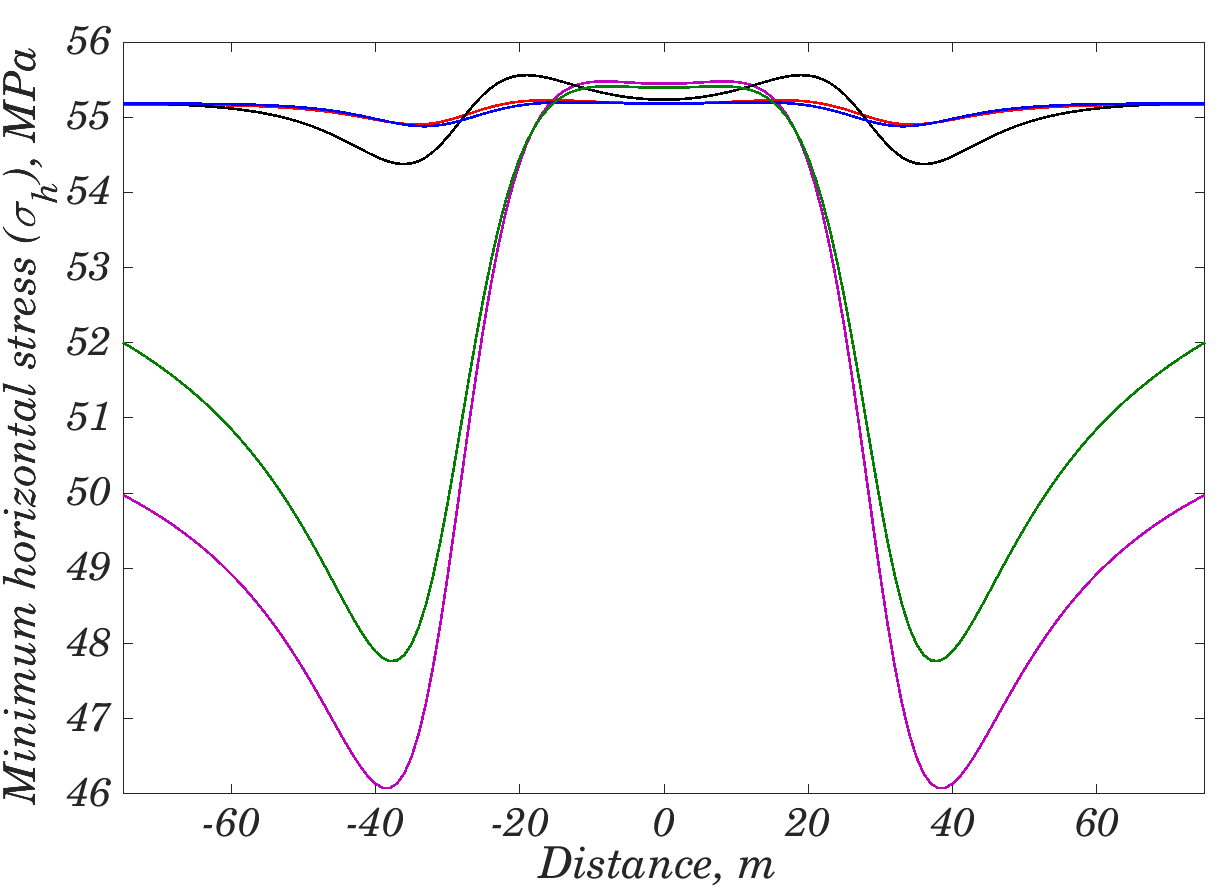}
			\label{fig:diff_syy_1mon}
		} 
		\caption{Response of different rocks after $1$ month of 
		production}
		\label{fig:diff_1mon}
	\end{figure}

	Figure \ref{fig:diff_1mon} shows pore pressure and stresses changes 
	along the 
	line of interest after one month of production. 
	Figure \ref{fig:diff_pp_1mon} shows the pore pressure depletion along 
	that 
	line. 
	Both sandstones experience more than $50\%$ reduction in the pore 
	pressure 
	in 
	Region $1$.
	However, Charcoal Granite shows no change in the pore pressure in 
	that 
	region,
	Tennessee Marble experiences a slightly small depletion, and 
	Haynesville 
	shale shows a slight increase in pore pressure in the same region.
	Maximum horizontal stress shows the same trend as pore pressure, 
	but 
	minimum 
	horizontal stress shows a reverse trend.
	Both sandstones have the highest reduction in the magnitudes of the 
	maximum 
	and minimum horizontal stresses and their stress anisotropies among 
	the 
	analyzed set. 
	Tennessee Marble and Charcoal Granite show a small reduction in the 
	magnitude 
	of these parameters,
	but the Haynesville Shale behavior is slightly different from other 
	rocks 
	in 
	Regions 1 and 2. 
	A slight increase in observed in pore pressure and maximum 
	horizontal 
	stress 
	of Haynesville Shale in Region 1,
	while its behavior falls in between two sandstones and ultra-tight 
	rocks.

\begin{figure}[!h]
		\centering
		
		\subfloat[Pore pressure]{
			\includegraphics[width=0.45\textwidth]{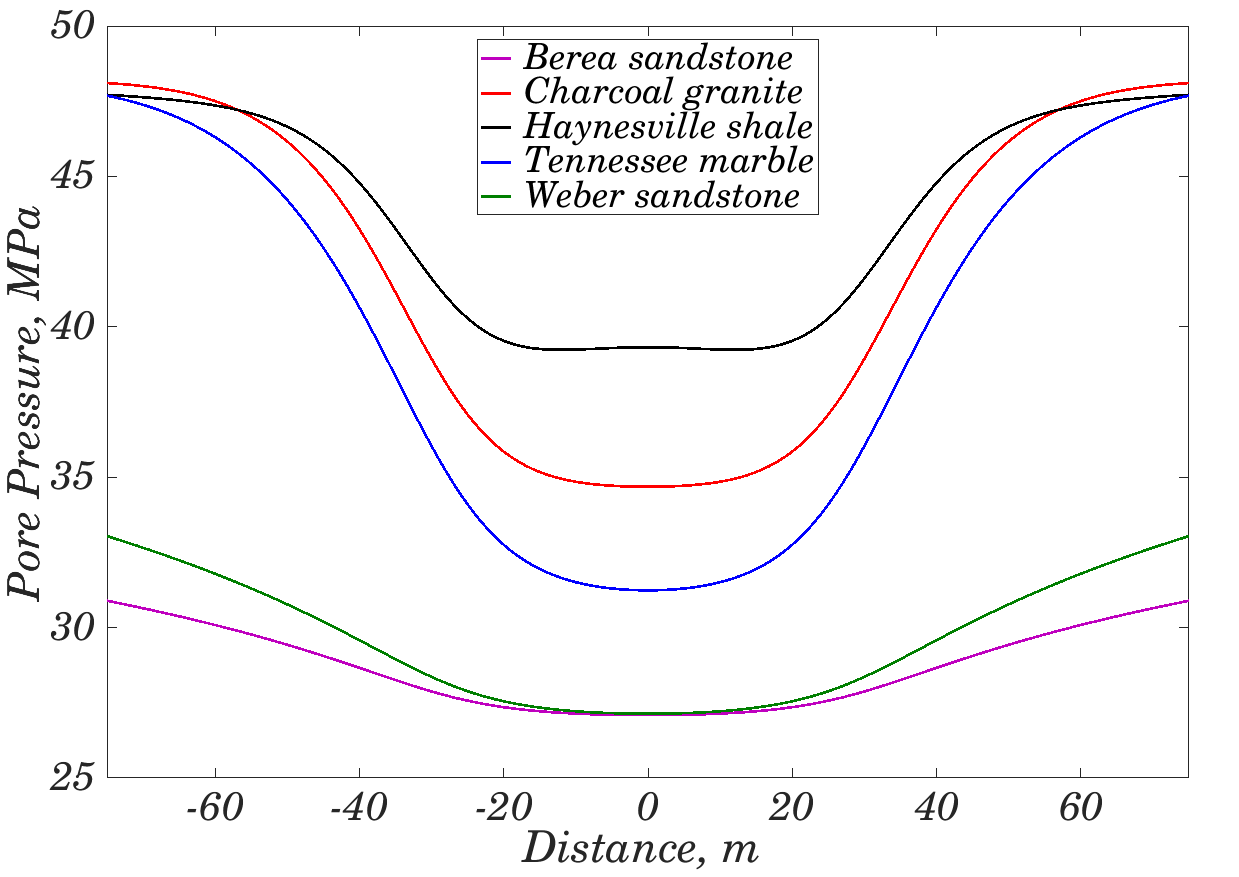}
			\label{fig:diff_pp_1year}
		}
		\subfloat[Stress anisotropy]{
			\includegraphics[width=0.45\textwidth]{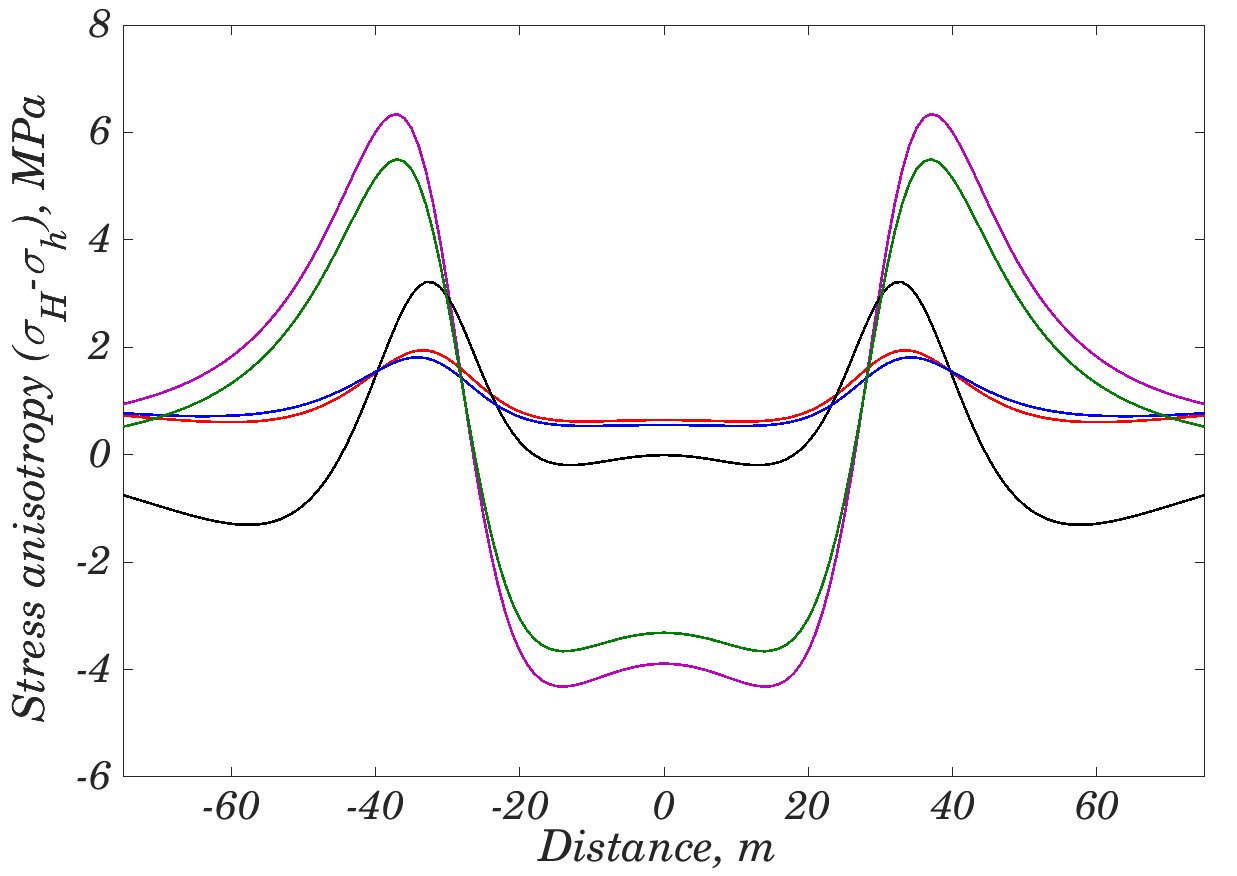}
			\label{fig:diff_aniso_1year}
		} 
		\qquad
		\subfloat[$\sigma_{xx}$]{
			\includegraphics[width=0.45\textwidth]{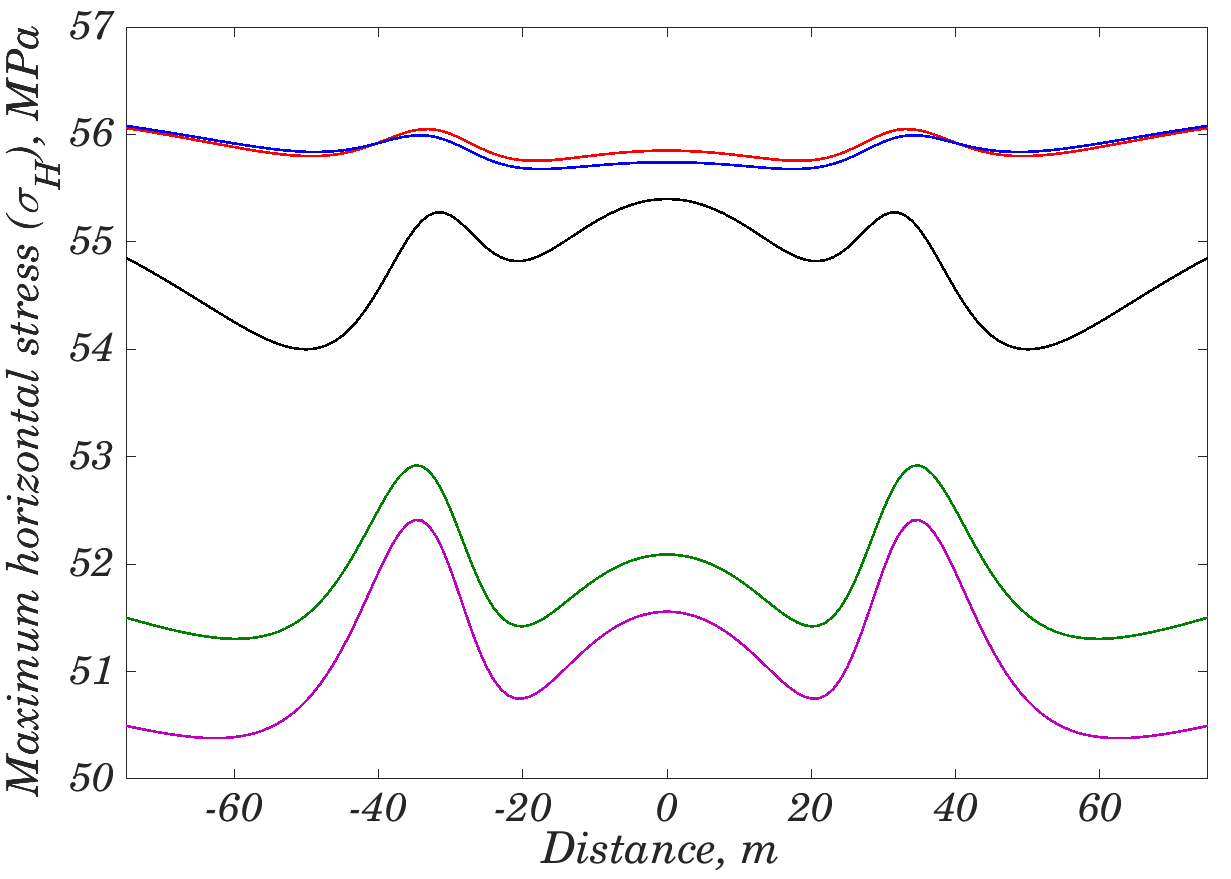}
			\label{fig:diff_sxx_1year}
		}
		\subfloat[$\sigma_{yy}$]{
			\includegraphics[width=0.45\textwidth]{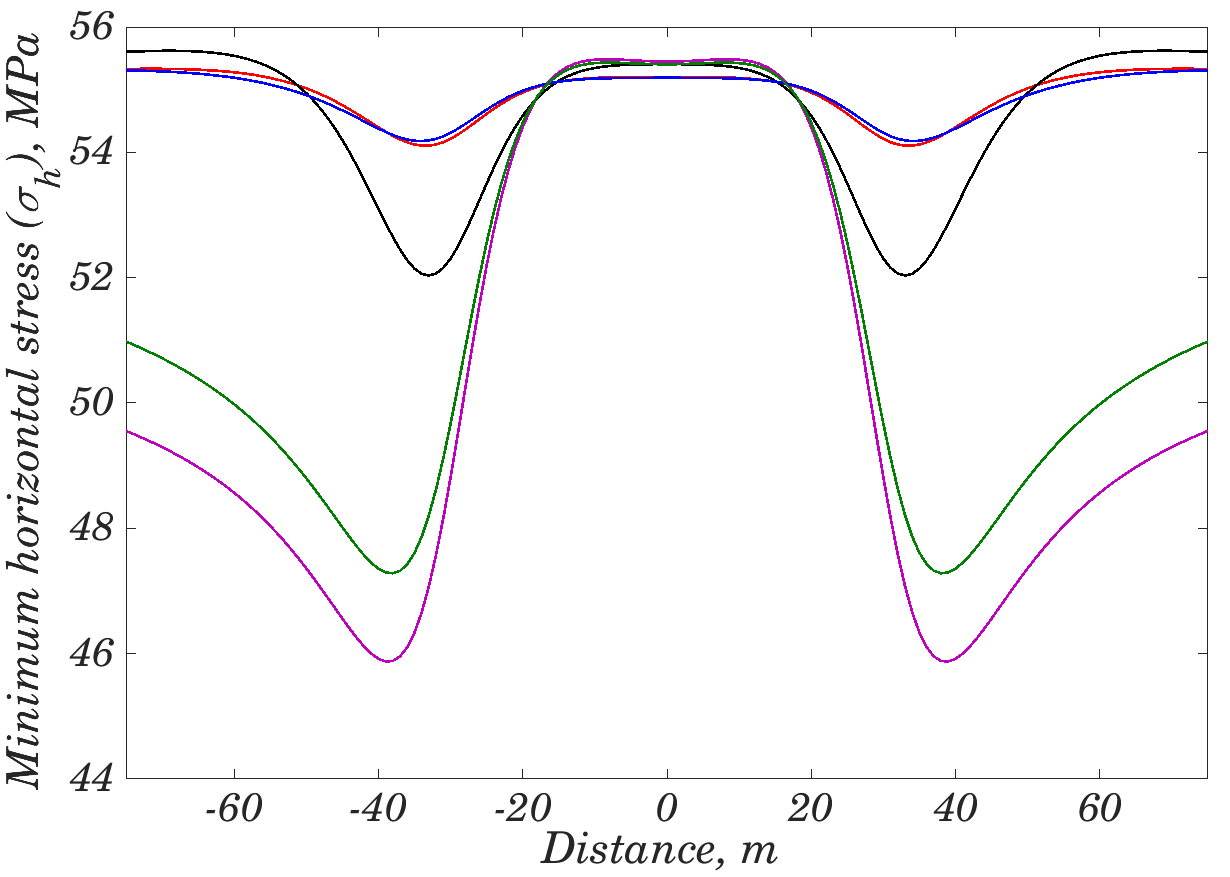}
			\label{fig:diff_syy_1year}
		} 
		\caption{Response of different rocks after $1$ year of 
		production}
		\label{fig:diff_1year}
	\end{figure}

	Another calculation is done on the same line after $1$ year from 
	production. 
	Figure \ref{fig:diff_1year} shows changes of the variables along the 
	interested 
	line. After one year, pore pressure depletion in Region $1$ has 
	reached to 
	the pre-set value ($27 \; MPa$) for both sandstones. Also, pore 
	pressure 
	reduction has traveled more in Region 2. The other three rock 
	samples show 
	partial depletion in Region 1. As is shown in Figure 
	\ref{fig:diff_pp_1year} 
	pore pressure reduction of 
	Haynesville Shale in Region 1 is highest among other rocks after 
	one year. 
	It is also observed that the shale shows a different behavior 
	in terms of the magnitude of both horizontal stresses and 
	anisotropy (Figure 
	\ref{fig:diff_sxx_1year} and \ref{fig:diff_syy_1year}) in both 
	Regions. 
	Moreover, the magnitude of these stresses is decreased and stays 
	between a 
	minimum that belongs to sandstones and a maximum that belongs to 
	Charcoal 
	Granite and Tennessee Marble.
	For the changes in stress anisotropy, after one year it is observed 
	that the stress anisotropy for shale is slightly below of the 
	original 	value	($1.38 \; MPa$). 
	
\begin{figure}[!h]
		\centering
		
		\subfloat[Pore pressure]{
			\includegraphics[width=0.45\textwidth]{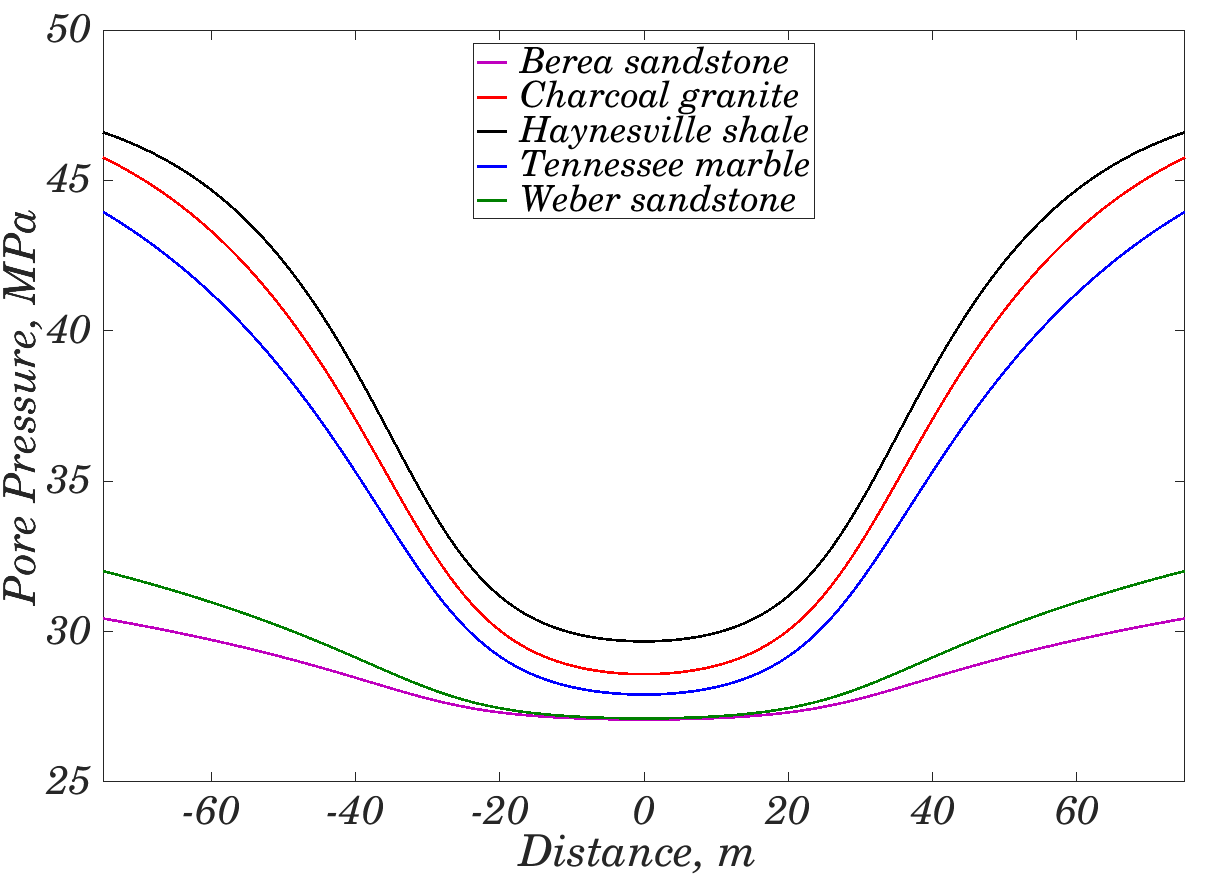}
			\label{fig:diff_pp_5year}
		}
		\subfloat[Stress anisotropy]{
			\includegraphics[width=0.45\textwidth]{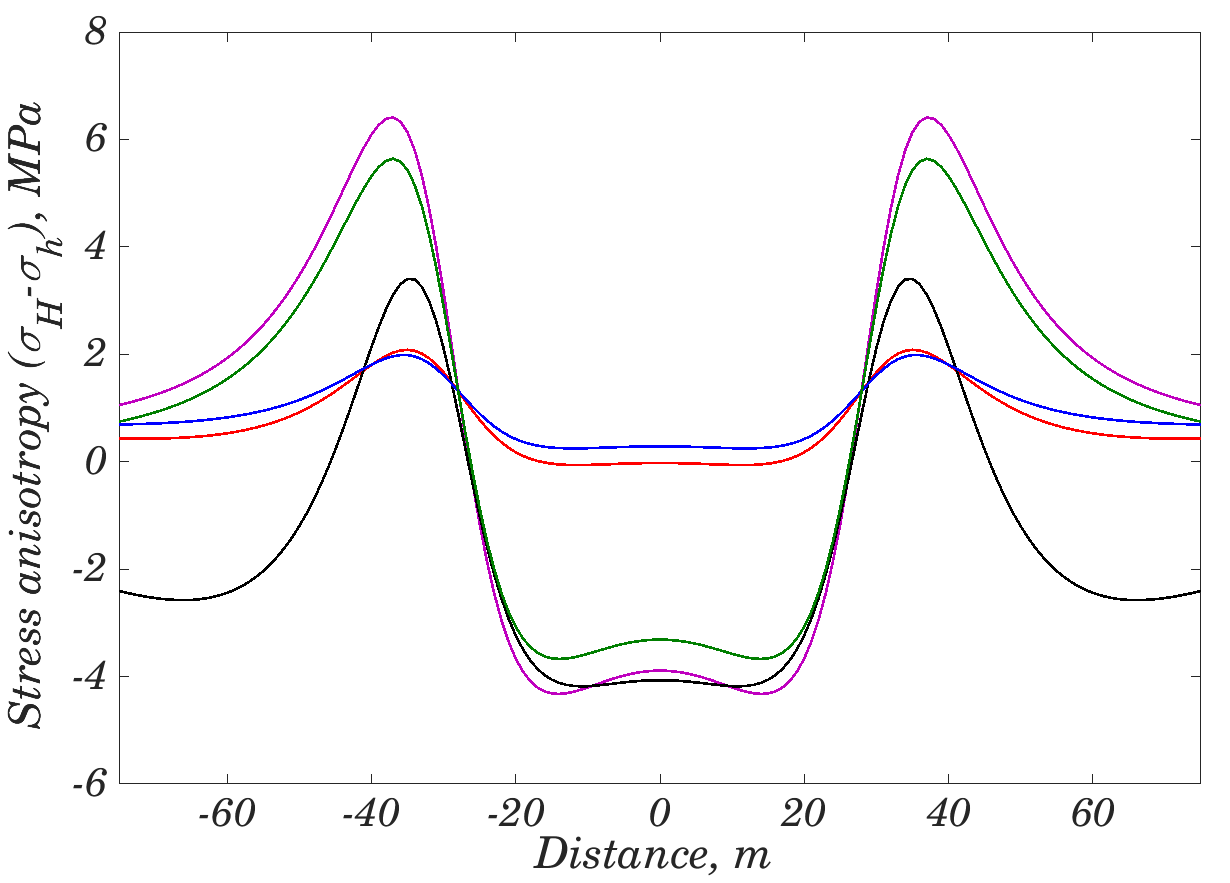}
			\label{fig:diff_aniso_5year}
		} 
		\qquad
		\subfloat[$\sigma_{xx}$]{
			\includegraphics[width=0.45\textwidth]{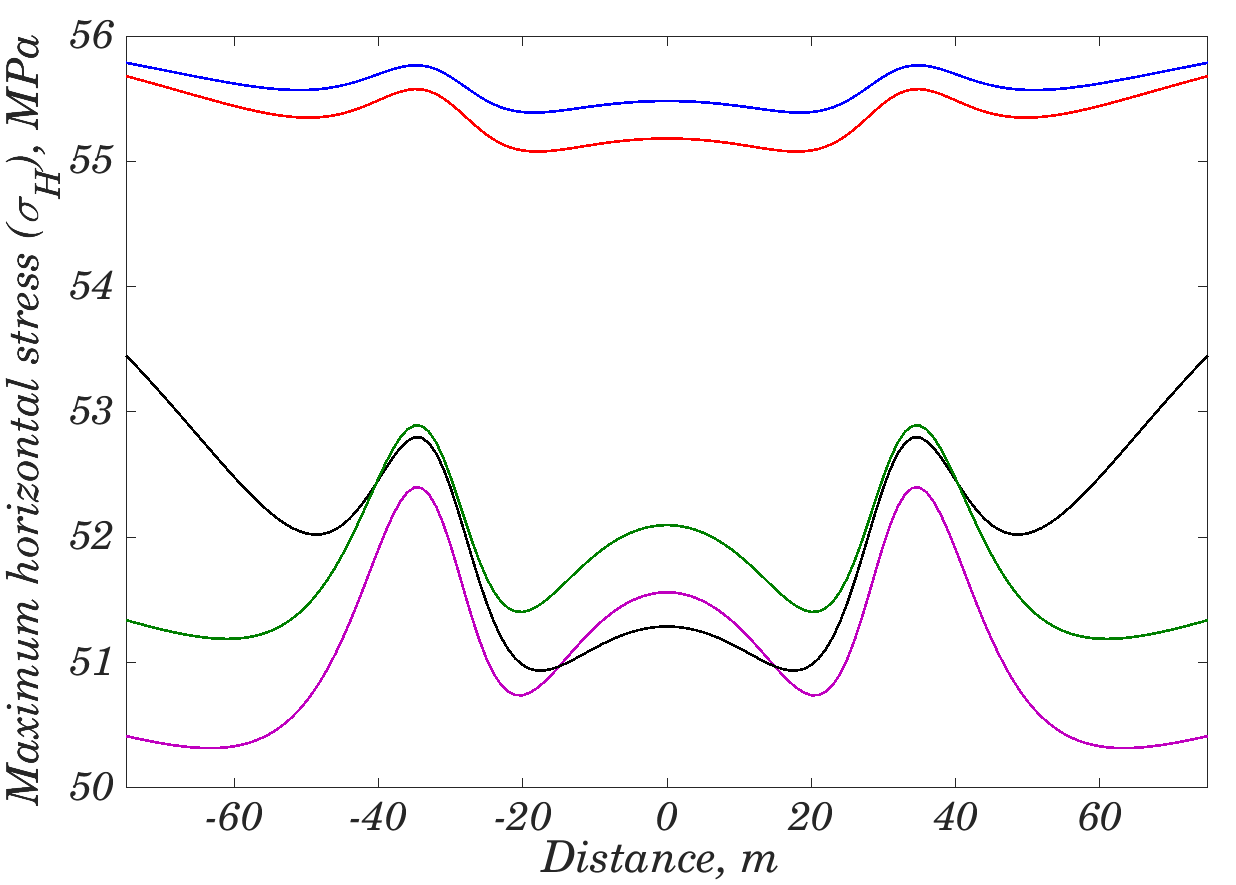}
			\label{fig:diff_sxx_5year}
		}
		\subfloat[$\sigma_{yy}$]{
			\includegraphics[width=0.45\textwidth]{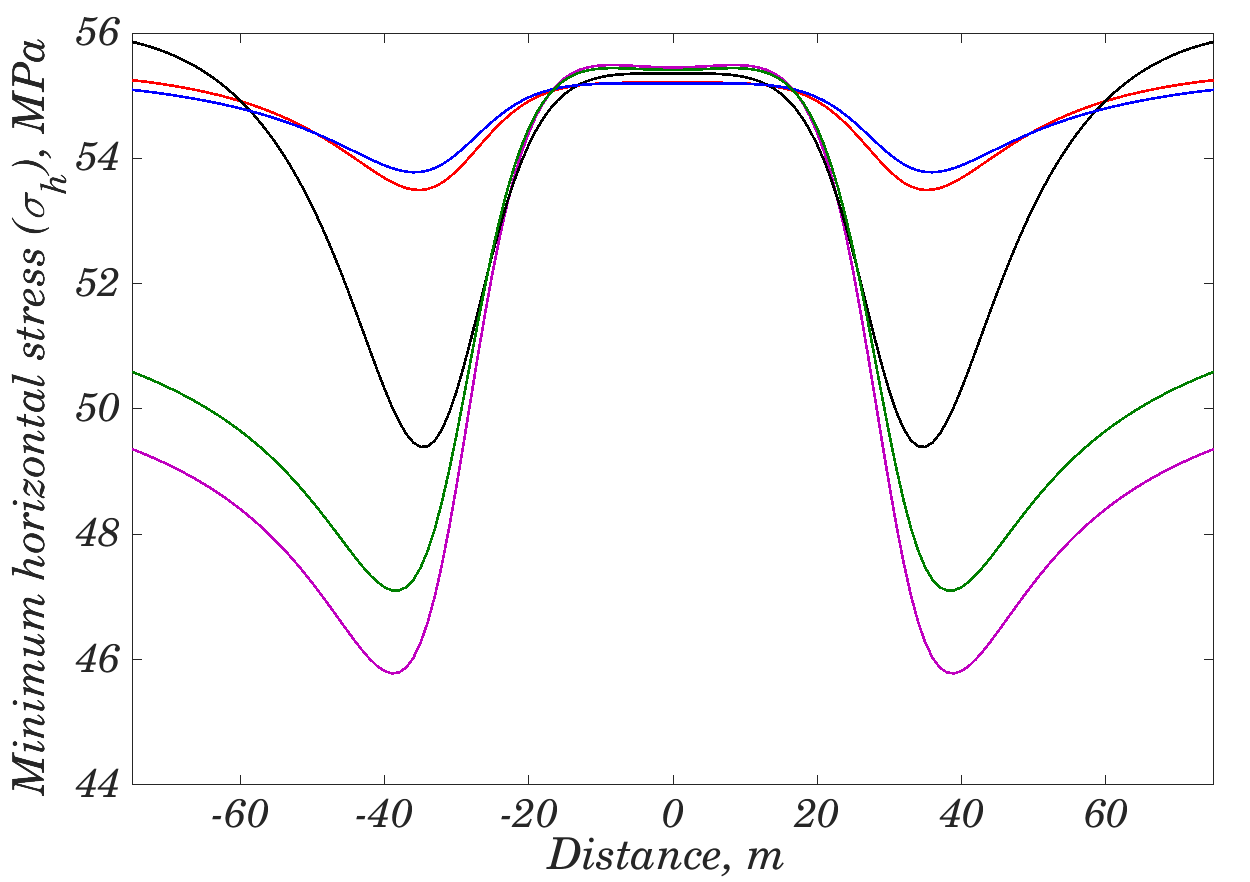}
			\label{fig:diff_syy_5year}
		} 
		\caption{Response of different rocks after $5$ years of 
		production}
		\label{fig:diff_5year}
	\end{figure}

	Plots of the same calculations after five years are shown in Figure 
	\ref{fig:diff_5year}. At this time, pore pressure reduction of the 
	Shale, 
	Charcoal granite, 
	and Tennessee Marble has not reached to the prescribed value yet 
	(this happened in sandstones after one year of production). 
	The pore pressure depletion of sandstones, unlike the other three 
	types of 
	rocks,
	has traveled in Region 2. 
	Comparing pore pressure depletion of the three tight rocks shows 
	that the 
	shale has been depleted less than the other two
	in Region 1. 
	Regarding stress, shale has decreased further down, and at this 
	time its 
	stress magnitudes are at the same level as sandstones. 
	The stress anisotropy plot after five years also shows that stress 
	reversal 
	has happened 
	for shale and two sandstones, although no change or slight change 
	is 
	observed 
	in the other two tight rocks.

	The results of the analysis that was performed using different rock 
	properties 
	may be used to categorize the rock samples that were used in this 
	section 
	into three groups having different behaviors.
	Berea Sandstone and Weber Sandstone formed the first group.
	This group of rocks showed complete depletion in Region 1 as early 
	as $1$ 
	month from production. 
	Pore pressure reduction for these two rocks also travels in more in 
	Region 
	2 
	compared to the other three rocks.
	Haynesville Shale showed a different behavior compared to the other 
	four 
	rocks chosen for this study under the loading condition and problem 
	geometry 
	that was described. The pore pressure depletion plots show that the 
	magnitude 
	of the stresses in both of sandstones stay low
	as a result of depletion during the period of production, and so 
	does their 
	stresses. 
	On the other hand, pore pressure depletion for Charcoal Granite and 
	Tennessee 
	Marble happens very slow, and their stress magnitudes stay the 
	highest 
	during 
	the period
	that was discussed. At early time, shale showed a slight increase 
	in the 
	amount of pore pressure and stresses. Also, after one year its pore 
	pressure 
	was reduced less than any other rock sample in this study,
	its stress magnitude was less than Charcoal Granite and Tennessee 
	Marble 
	and 
	higher than the sandstones. 
	After five years of production, although the area between fractures 
	has not 
	been depleted entirely in shales, the stresses reach their minimum, 
	as low 
	as sandstones, and stress reversal happens.

	\section{GLOBAL SENSITIVITY ANALYSIS OF PORE PRESSURE AND STRESSES} 
	\label{Sec:Global_analysis}
	
\begin{figure}[!h]
		\centering
		\includegraphics[width = 0.5\textwidth]{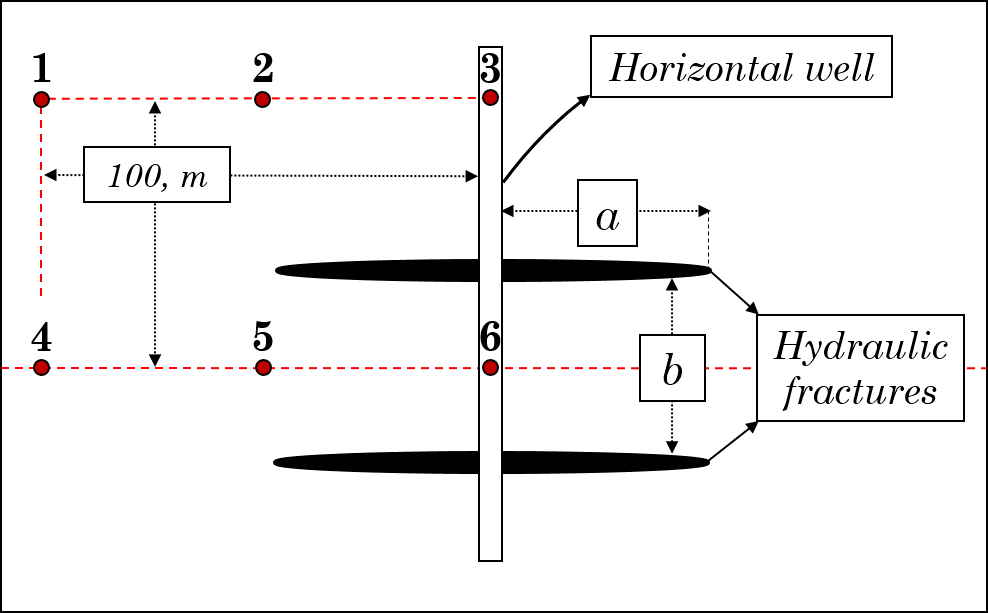}
		\caption{Locations of the six points that are chosen for 
		calculating the 
		changes in QI} 
		\label{fig:sampling} 
	\end{figure}

	In this section, Sobol method is used to analyze the variation of 
	the model 
	output parameters resulting from changes in the inputs. Pore 
	pressure 
	$p_p$, maximum horizontal $\sigma_H$, and minimum horizontal stress 
	$\sigma_h$ are chosen as our quantities of interest (QI). We are 
	precisely 
	interested in tracking the changes of QI at six points around 
	hydraulic 
	fractures and horizontal well as the production time increases. 
	Locations 
	of these six points are shown in Figure \ref{fig:sampling}. Among 
	them, Points $4$ and $1$ are on the path of the infill (off-set) 
	horizontal well and Points $2$, $3$, $5$, and $6$ are on the path 
	of child 
	fracture propagation in refracturing or infill well stimulation. 
	We, 
	intentionally, set the far-field stresses and pore pressure 
	constant during 
	production to limit our analysis to a single well. The prescribed 
	values 
	for these fixed boundary conditions are given in Table 
	\ref{table:model_input}. To track the changes of QI with time, the 
	model is 
	set for three different production periods of $1\;month$, 
	$1\;years$, and 
	$3\;years$. For each of these cases, the simulator is run for a 
	total of 
	$10800$ times to generate the output vectors. Three vectors for 
	pore 
	pressure, maximum horizontal, and minimum horizontal stresses are 
	generated 
	at the end of each period of production. 
	
	\begin{table}[!h]
		\centering
		\begin{tabular}{ll}
			\toprule\noalign{\smallskip}
			Maximum horizontal stress, $\sigma_H$ & $58.60 \text{ MPa}$ 
			\\ \hline
			Minimum horizontal stress, $\sigma_h$ & $55.15 \text{ MPa}$ 
			\\ \hline
			Reservoir pore pressure, $p_r$ & $48.26 \text{ MPa}$ \\ 
			\noalign{\smallskip}\bottomrule
		\end{tabular}
		\caption{Input parameters of for the reservoir far-field 
		properties}
		\label{table:model_input}
	\end{table}
	
	We selected eight independent input variables for the global 
	sensitivity 
	analysis. The input variables are divided into two separate 
	categories, 
	namely 
	the stimulation design variables and rock properties. Table \ref 
	{tab:sens_param} presents these variables and their minimum and 
	maximum 
	values. 
	The design variables include fracture half-length $a$, fracture 
	spacing 
	$b$, 
	and production pressure (fracture pressure) $p_f$. For all of the 
	simulations, 
	fracture pressure was held constant during the period of the 
	production. 
	Rock 
	properties the were selected for sample generation include shear 
	modulus 
	$G$, 
	undrained Poisson's ratio $\nu_u$, drained Poisson's ratio $\nu$, 
	Skemtson's 
	coefficient $B$, and mobility $\kappa$ (rock absolute permeability 
	$k$/fluid 
	viscosity $\mu$).

	\begin{table}[!h]
		\centering
		\caption{Rock poroelastic properties and design variables that 
		are used 
		global sensitivity analysis. The ranges of rock properties are 
		chosen 
		from 
		\cite{cheng2016poroelasticity}. Note that in this section and the following 
		sections of the paper, $\log()$ of the fluid mobility is refereed to as 
		$\kappa$.}
			\label{tab:sens_param}  
		\begin{tabular}{lllllllll}
			\toprule
			& \multicolumn{3}{c}{Design variables} & 
			\multicolumn{5}{c}{Rock properties} \\
			\midrule\noalign{\smallskip}
			Property & $a,\;m$ & $b,\;m$ & $p_f, \; MPa$ & $G, GPa$ & 
			$\nu_u$ & 
			$\nu$ & $B$ & $\log(\kappa), \frac{m^2}{Pa \cdot s}$ \\
			\noalign{\smallskip}\hline\noalign{\smallskip}
			Minimum & $10$ & $10$ & $10$ & 
			$1$ & $0.30$ & $0.1$ & $0.3$ & $-17$ \\ 
			Maximum & $60$ & $30$ & $40$ & $25$ 
			& $0.45$ & $0.29$ & $0.9$ & $-10$ \\ 
			\noalign{\smallskip}\hline\noalign{\smallskip}
			Sobol index & $S_1$ & $S_2$ & $S_3$ & $S_4$ 
			& $S_5$ & $S_6$ & $S_7$ & $S_8$ \\
			\noalign{\smallskip} \bottomrule
		\end{tabular}
	\end{table}
	
	Referring to Equation \eqref{eq:ANOVA_rearang}, number of total terms 
	for the 
	analysis using $8$ input variables will be $2^8\;= 256$, which is 
	also 
	equal 
	to 
	the total number of Sobol indices. Each of the variables in Table 
	\ref 
	{tab:sens_param} are associated with a Sobol index as shown in the 
	table. 
	Using 
	this convention, a first order Sobol index $S_i$ is related to 
	individual 
	contributions of the variables on the output, while higher order 
	Sobol 
	indices 
	$S_{i \dots j}$ represent the interaction between the variables. 
	For 
	example, 
	$S_1$ represents the effect of changes in fracture half-length 
	magnitude to 
	the 
	variation of a specific output variable (e.g., pore pressure), and 
	$S_{37}$ 
	represents this variation due to the simultaneous changes in 
	fracture 
	pressure 
	$p_f$ and Skemptson's coefficient $B$. Using this procedure, three 
	output 
	vectors for pore pressure, maximum stress, and minimum stress are 
	generated 
	for 
	each period, and each of these vectors is analyzed. An open-source 
	library 
	developed by \cite{herman2017salib} for sensitivity analysis is 
	used to 
	perform 
	both the sample generation and Sobol analysis. Results are 
	discussed in the 
	next section.

	\subsection{Analysis of the results after one month} 
	\label{subsec:one_month}
	
\begin{figure}[!h]
		\centering 
		\subfloat[Individual effects]{
			\includegraphics[width=0.45\textwidth]{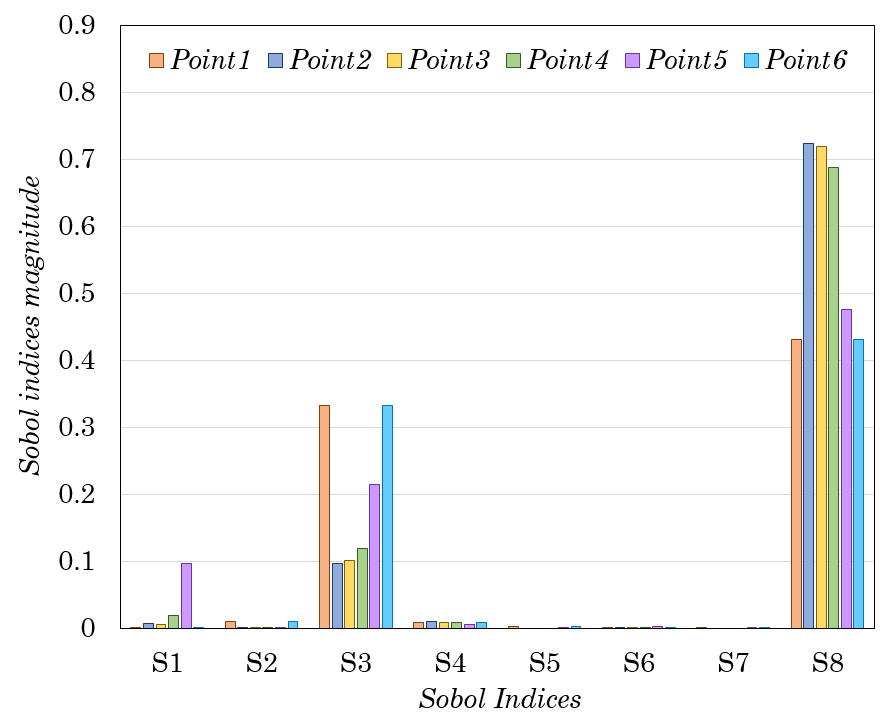}
			\label{fig:press_Si_1month}
		}
		\subfloat[Interaction effects ($>0.01$)]{
			\includegraphics[width=0.45\textwidth]{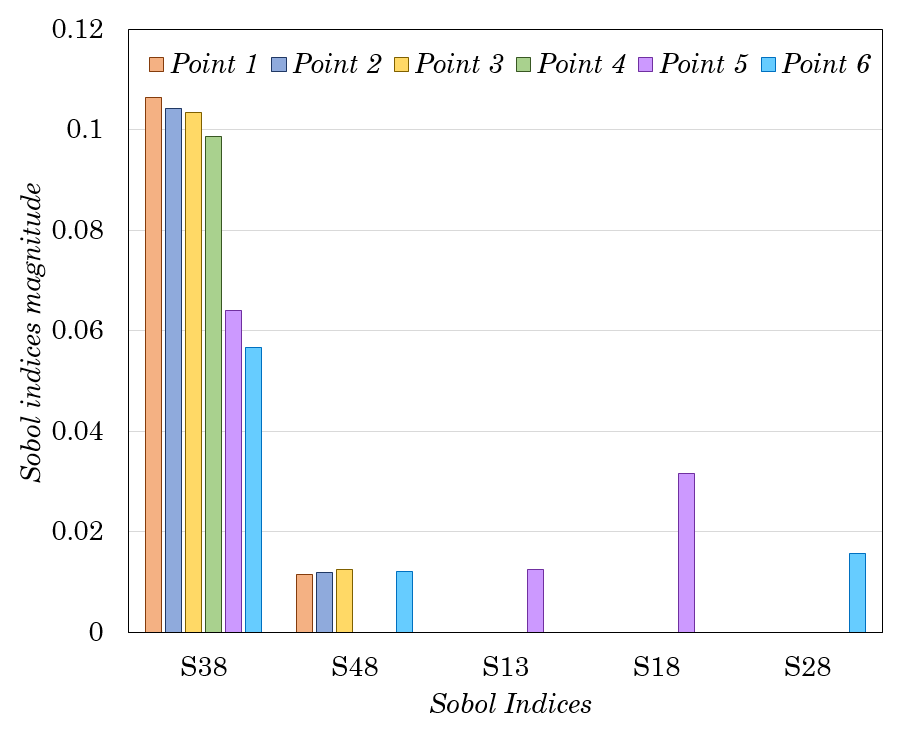}
			\label{fig:press_Sij_1month}
		} 
		\caption{Sobol indices of the input sample for pore pressure 
		after 1 
		month from production}
		\label{fig:press_1month}
	\end{figure}

\begin{figure}[!h]
		\centering 
		\subfloat[Individual effects]{
			\includegraphics[width=0.45\textwidth]{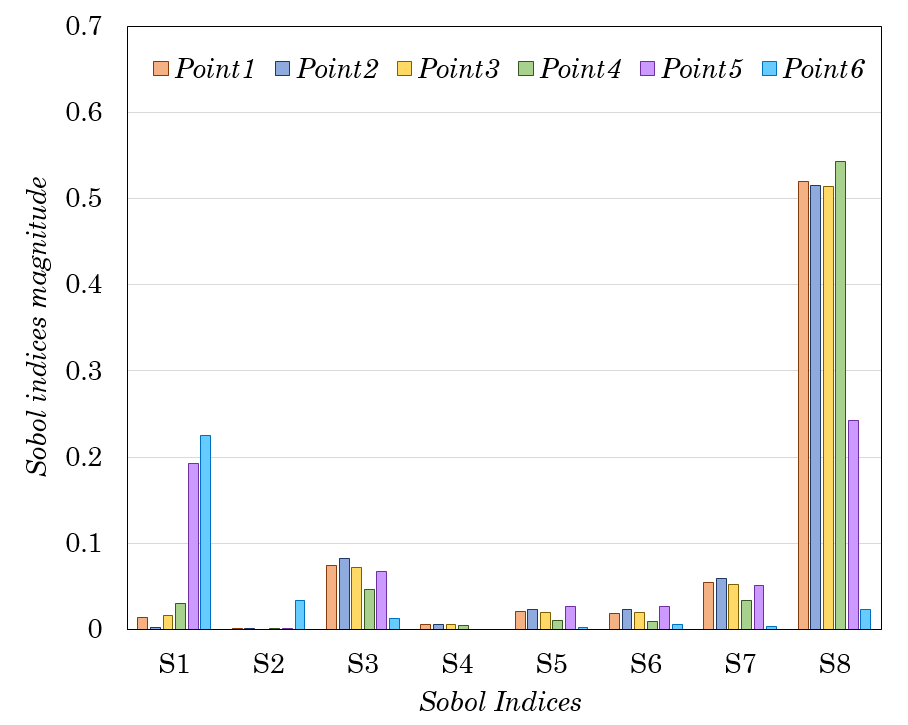}
			\label{fig:MinStress_Si_1month}
		}
		\subfloat[Interaction effects ($>0.03$)]{
			\includegraphics[width=0.45\textwidth]{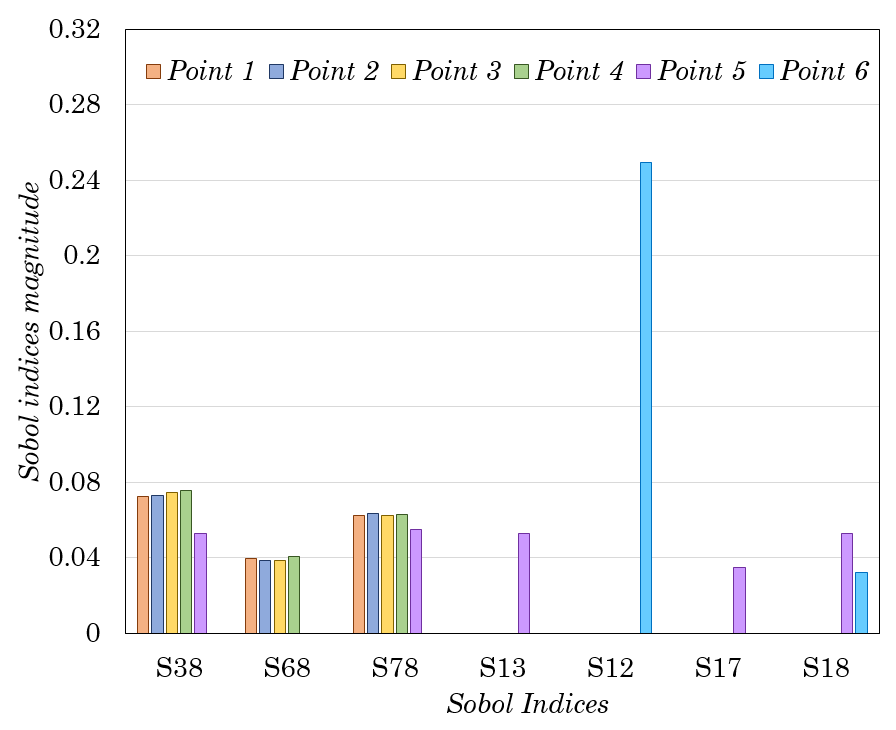}
			\label{fig:MinStress_Sij_1month}
		} 
		\caption{Sobol indices of the input sample for minimum 
		horizontal 
		stress after 1 month from production}
		\label{fig:MinStress_1month}
	\end{figure}

\begin{figure}[!h]
		\centering 
		\subfloat[Individual effects]{
			\includegraphics[width=0.45\textwidth]{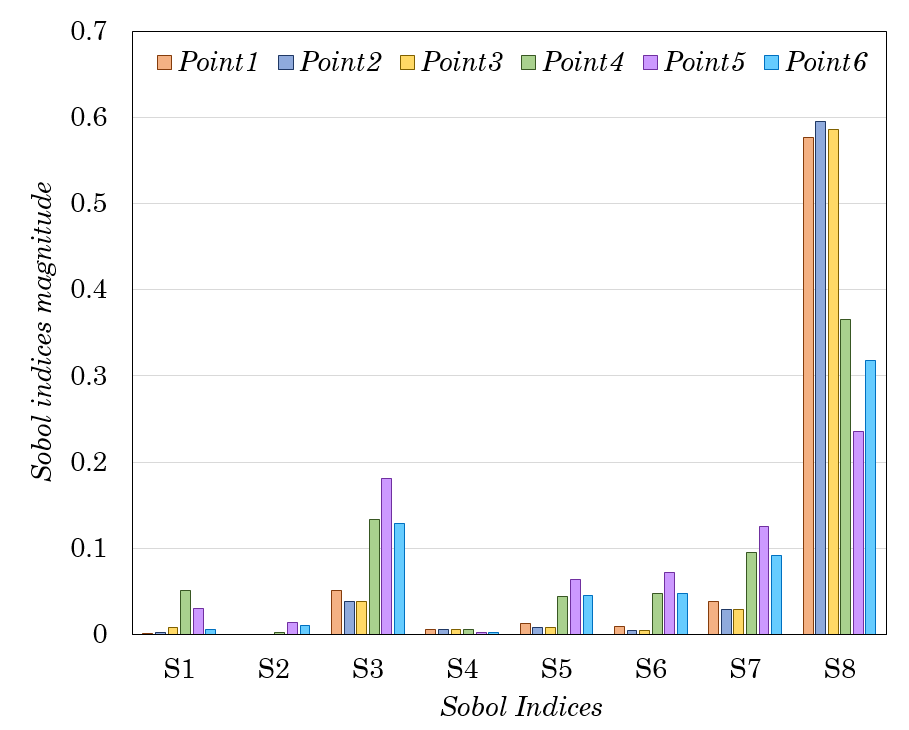}
			\label{fig:MaxStress_Si_1month}
		}
		\subfloat[Interaction effects ($>0.03$)]{
			\includegraphics[width=0.45\textwidth]{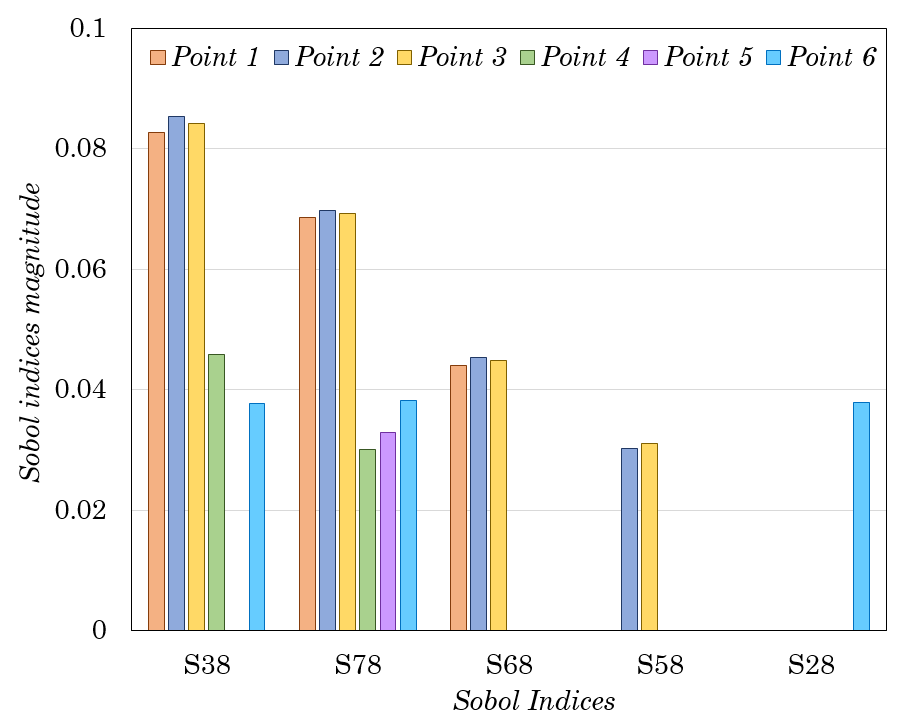}
			\label{fig:MaxStress_Sij_1month}
		} 
		\caption{Sobol indices of the input sample for maximum 
		horizontal 
		stress after 1 month from production}
		\label{fig:MaxStress_1month}
	\end{figure}

	In this section, a global sensitivity of the pore pressure and 
	stresses 
	after 
	$1\;month$ of production is presented. After one month, it not 
	expected 
	from 
	pore pressure depletion extent to move far away from the wellbore. 
	Figures 
	\ref{fig:press_1month}--\ref{fig:MaxStress_1month} show the result 
	of SA 
	for 
	this production period. As presented in Figure 
	\ref{fig:press_Si_1month} for 
	pore pressure, the main individual contributions to changes of the 
	pore 
	pressure are $S_8$ and $S_3$, which correspond to mobility $\kappa$ 
	and 
	production pressure $p_f$ respectively. These two variables 
	contribute to 
	more 
	than $80\%$ of the changes in pore pressure changes around 
	hydraulic 
	fractures 
	and horizontal well. Also, an interesting observation is that $S_1$ 
	(fracture 
	half-length) contributes to $10\%$ of the changes at Point $5$. 
	This is 
	important because Point $5$ is the desired path for propagation of 
	the 
	child 
	fracture in refracturing process. Figure \ref{fig:press_Sij_1month} 
	shows the 
	interaction effects $S_{ij}$ between the selected input parameters. 
	Any 
	Sobol 
	index that has a value less than $0.01$ is excluded. This does 
	not put any limitation on our analysis since $\Sigma S_{i\dots j} 
	\approx 
	0.85-0.9$ at each point. As expected, the interaction between 
	production 
	pressure 
	and mobility have the greatest contribution among other 
	interactions. Also, 
	a 
	small contribution is observed from $S_{48}$ at Points $1$, $2$, 
	$3$, and 
	$6$. 
	
	Figure \ref{fig:MinStress_Si_1month} shows the individual 
	contributions of 
	inputs on the minimum horizontal stress. It is observed that $B$, 
	$p_f$, 
	and 
	drained Poisson's ratio $\nu$ are the main contributors to the 
	changes of 
	the 
	minimum horizontal stress at the points that are located outside 
	the 
	fracture 
	spacing (i.e., Points $1$, $2$, $3$, and $4$) in this case. 
	However, for 
	the 
	points close to the fracture tip area, and inside the spacing 
	(i.e., Points 
	$5$ 
	and $6$), different results are observed. At Point $5$, fracture 
	half-length 
	$a$ and $B$ are significant contributors, while at Point $6$ 
	fracture 
	half-length $a$ and fracture spacing $b$ are the dominating 
	contributors 
	among 
	other individual terms. Also, it is observed that changes of 
	fracture 
	half-length and fracture spacing ($S_{12}$) have the most 
	significant 
	impact 
	($\approx\;24\%$) on the minimum horizontal stress at Point $6$, 
	compared 
	to 
	all other variables. This shows the importance of these two 
	variables on 
	the 
	minimum horizontal stress changes. Note that Point $6$ is the 
	possible path 
	of 
	child fracture initiation in refracturing. At the points outside 
	the 
	fracture 
	spacing area, however, $S_{38}$, $S_{68}$, and $S_{78}$ have the 
	main 
	contributions to changes of minimum horizontal stress. 
	
	Figure \ref{fig:MaxStress_Si_1month} shows the first order Sobol 
	indices 
	$S_i$ 
	for maximum horizontal stress. In this case, main contributions are 
	$S_8$, 
	$S_3$, $S_7$, $S_6$, and $S_5$ respectively. No significant 
	contribution 
	from 
	design variables $a$ and $b$ is observed. Also, The main effect for 
	interaction 
	terms are due to $S_{38}$, $S_{78}$, $S_{68}$. From this, it can be 
	concluded 
	that the main contributions to the changes in maximum horizontal 
	stress are 
	due 
	to the rock properties rather than the design variables. These 
	variables 
	contribute to about $70\%$ of the changes in maximum horizontal 
	stress. It 
	should be noted that although $p_f$ is among design variables, the 
	main 
	source 
	of the changes due to this variable is the difference between this 
	variable 
	and 
	far-field pore pressure ($\Delta p$).  
	
	\subsection{Analysis of the results after one year} 
	\label{subsec:OneYear}
	
\begin{figure}[!h]
		\centering 
		\subfloat[Individual effects]{
			\includegraphics[width=0.45\textwidth]{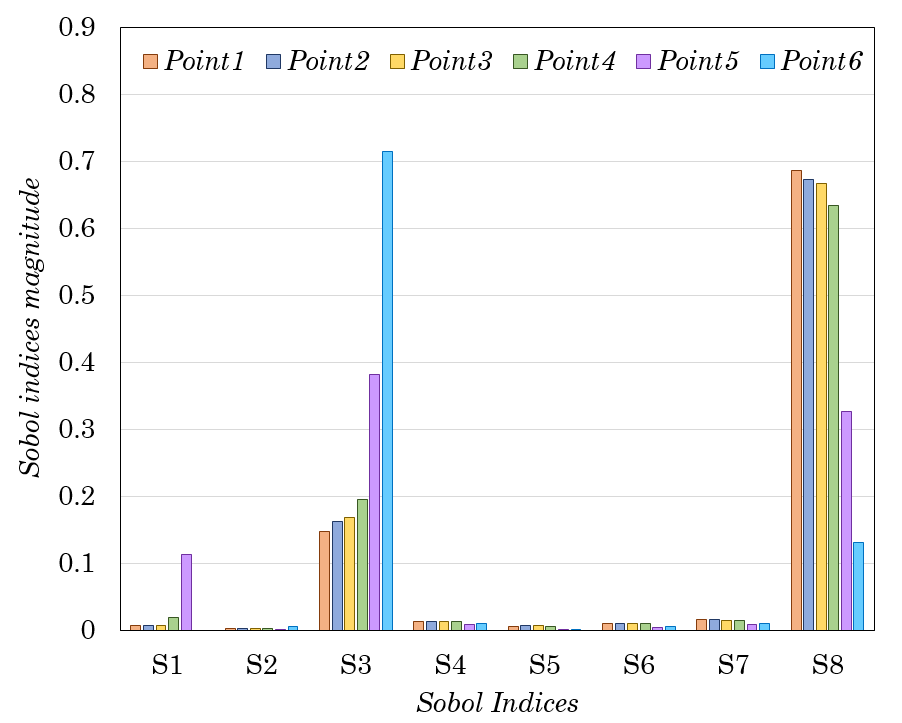}
			\label{fig:press_Si_1year}
		}
		\subfloat[Interaction effects ($>0.03$)]{
			\includegraphics[width=0.45\textwidth]{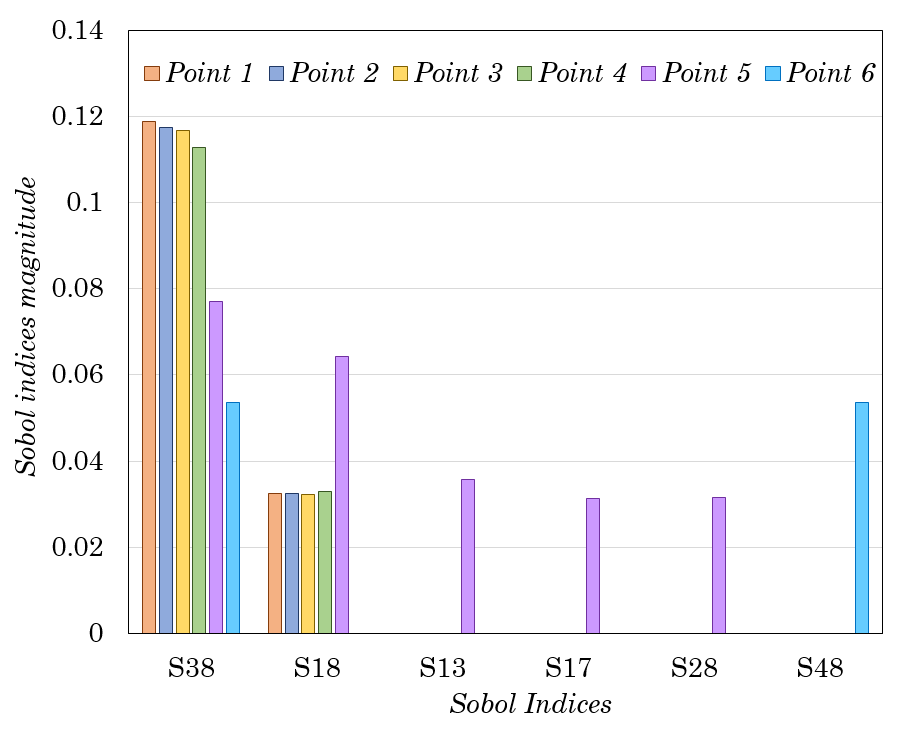}
			\label{fig:press_Sij_1year}
		} 
		\caption{Sobol indices of the input sample for pore pressure 
		after 1 
		year from production}
		\label{fig:press_1year}
	\end{figure}

\begin{figure}[!h]
		\centering 
		\subfloat[Individual effects]{
			\includegraphics[width=0.45\textwidth]{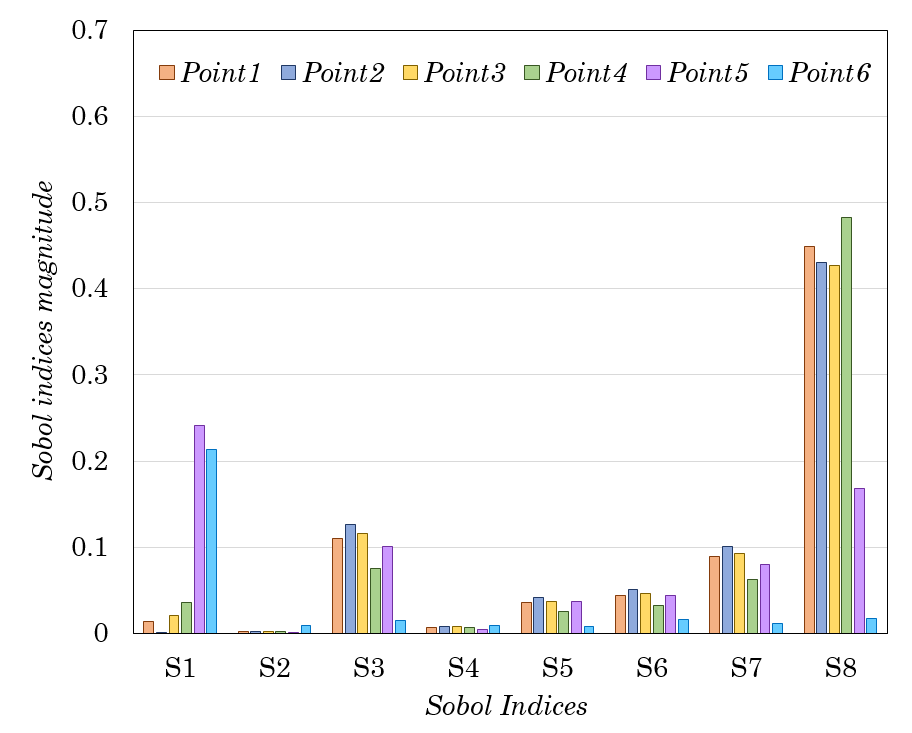}
			\label{fig:MinStress_Si_1year}
		}
		\subfloat[Interaction effects ($>0.03$)]{
			\includegraphics[width=0.45\textwidth]{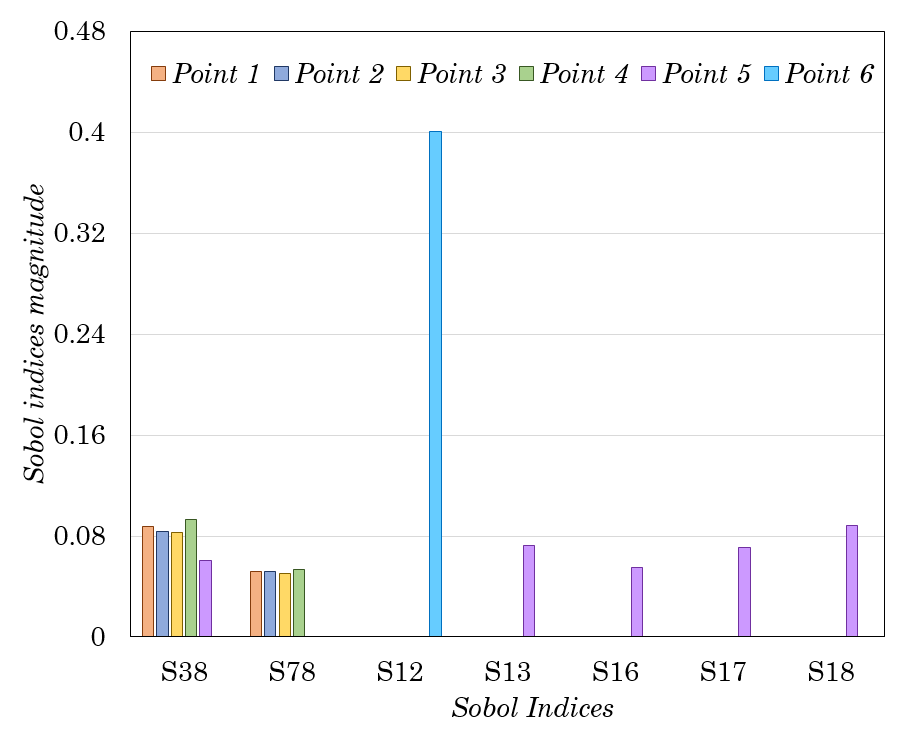}
			\label{fig:MinStress_Sij_1year}
		} 
		\caption{Sobol indices of the input sample for minimum 
		horizontal 
		stress after 1 year from production}
		\label{fig:MinStress_1year}
	\end{figure}

\begin{figure}[!h]
		\centering 
		\subfloat[Individual effects]{
			\includegraphics[width=0.45\textwidth]{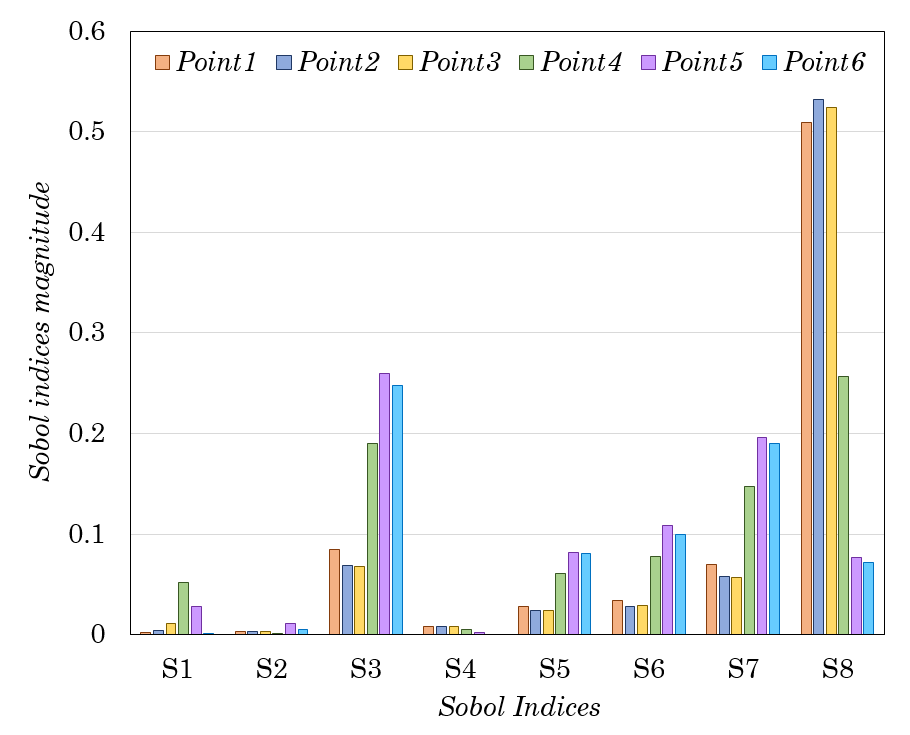}
			\label{fig:MaxStress_Si_1year}
		}
		\subfloat[Interaction effects ($>0.03$)]{
			\includegraphics[width=0.45\textwidth]{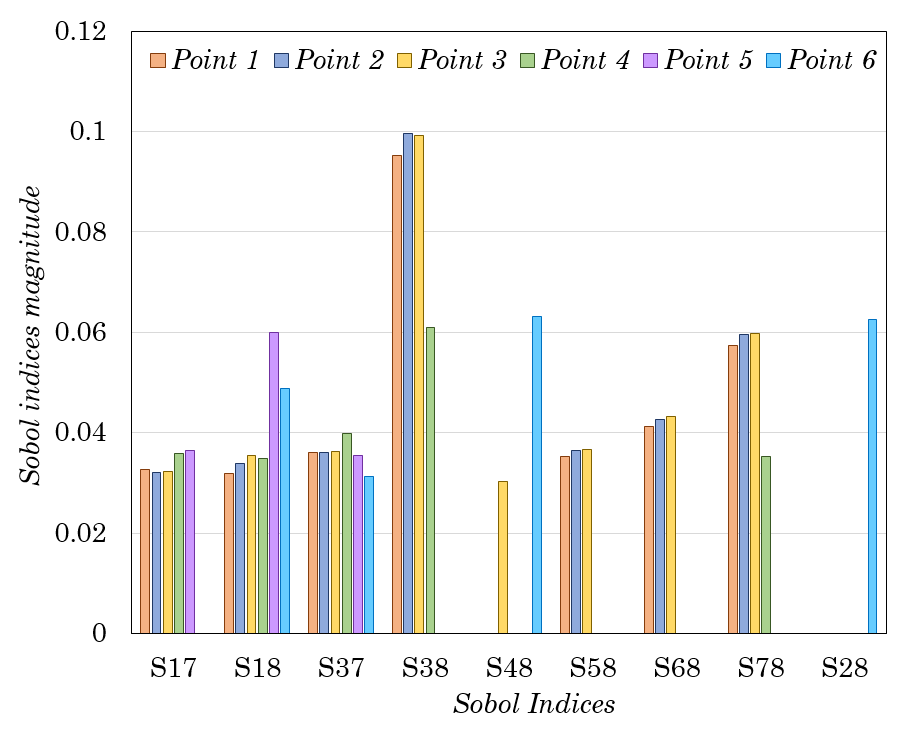}
			\label{fig:MaxStress_Sij_1year}
		} 
		\caption{Sobol indices of the input sample for maximum 
		horizontal 
		stress after 1 year from production}
		\label{fig:MaxStress_1year}
	\end{figure}

	Figures \ref{fig:press_1year}--\ref{fig:MaxStress_1year} show the 
	results of 
	Sobol analysis after one year from production. There are some 
	similarities 
	and several differences in the dominant contributors compared to 
	the 
	results after $1\;month$. For example, similar to the $1$ month 
	case, 
	$S_{3}$ and $S_8$ are still dominating individual contributors to 
	the 
	variance of the pore pressure as shown in Figure 
	\ref{fig:press_Si_1year}. In 
	contrary, unlike the $1\;month$ case, the effect of $S_3$ is 
	greater than 
	$S_8$ at Point $6$ after $1\;year$. The main reason for this is 
	that the 
	effect of pore pressure depletion has reached to the half-way 
	between 
	fractures after this period. Among higher order Sobol indices, 
	$S_{38}$ has 
	the greatest impact on the pore pressure changes. The second 
	greatest 
	impact from interaction terms is due to $S_{18}$, while after a 
	month 
	$S_{48}$ was the second dominant contributor among higher-order 
	indices. 
	
	Sobol indices for minimum horizontal stress are presented in Figure 
	\ref{fig:MinStress_1year}. The dominant individual Sobol indices 
	are the 
	same as $1\;month$ in this case. The only difference is that all of 
	the 
	rock properties variables have a more significant Sobol index after 
	one 
	year, indicating advancement of the pore pressure depletion in the 
	rock. 
	Among interaction terms, $S_{38}$ and $S_{78}$ are the dominating 
	interactions. Unlike the $1\;month$ case, there is no significant 
	impact 
	from $S_{68}$. Also, $S_{12}$ is still the greatest Sobol index at 
	Point 
	$6$, but its value has increased from $24\%$ to $40\%$.
	
	Moreover, the combination of the fracture length and rock 
	properties have a 
	considerable impact on the changes in minimum horizontal stress at 
	Point 
	$5$. 
	This makes sense because this point is in front of the fracture 
	tip. The 
	total contribution of the individual and interaction terms that 
	include 
	fracture length on the variance of minimum horizontal stress 
	changes at 
	Point 
	$5$ is almost $65\%$. 
	
	Figure \ref{fig:MaxStress_1year} represents the Sobol indices for 
	maximum 
	horizontal stress after one year. It can be seen that after one 
	year, 
	mobility has a smaller effect on the changes in maximum horizontal 
	stress 
	compared to the one month case. But, the contribution from $S_3$, 
	$S_5$, 
	$S_6$, and $S_7$ have increased. It is also observed that $S_3$ and 
	$S_7$ 
	are the dominant Sobol indices after one year. Among higher-order 
	indices, 
	$S_{17}$, $S_{18}$, and $S_37$ are the dominant terms. Also, other 
	terms 
	such as $S_{38}$, $S_{48}$, $S_{58}$, $S_{68}$, $S_{78}$, and 
	$S_{28}$ have 
	considerable effects, but they do not affect all of the points 
	equally 
	
	\subsection{Analysis of the results after three years} 
	\label{subsec:three_years}
	
\begin{figure}[!h]
		\centering 
		\subfloat[Individual effects]{
			\includegraphics[width=0.45\textwidth]{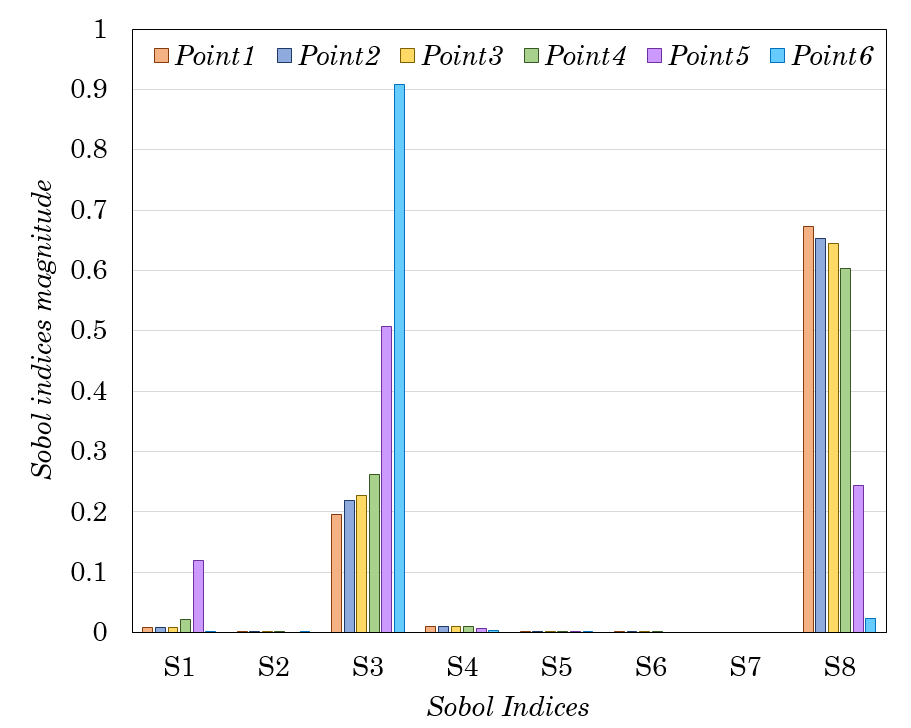}
			\label{fig:press_Si_3years}
		}
		\subfloat[Interaction effects ($>0.03$)]{
			\includegraphics[width=0.45\textwidth]{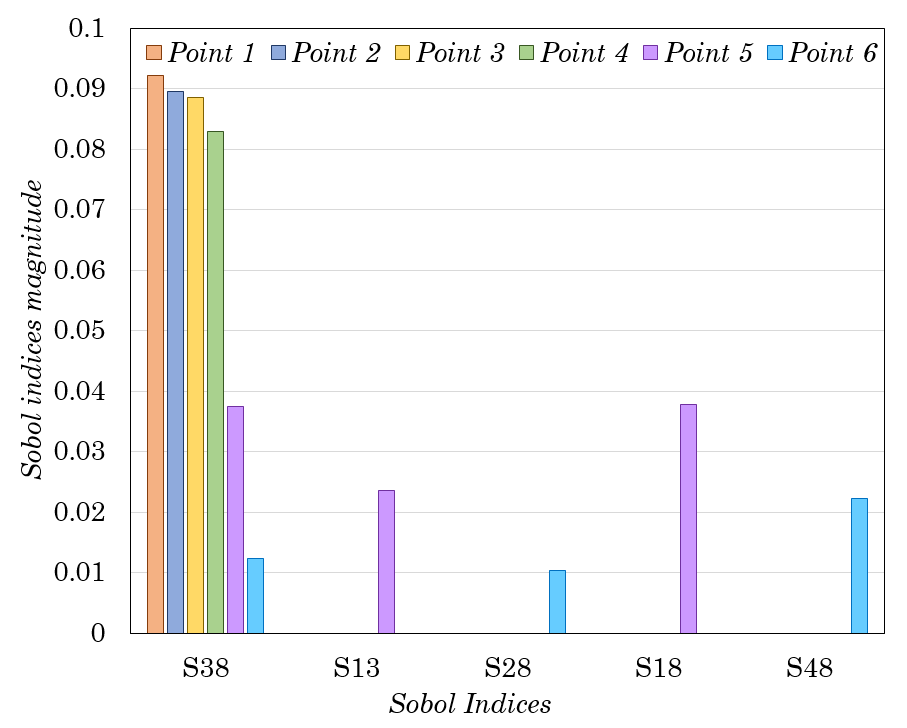}
			\label{fig:press_Sij_3years}
		} 
		\caption{Sobol indices of the input sample for pore pressure 
		after 3 
		years from production}
		\label{fig:press_3years}
	\end{figure}

\begin{figure}[!h]
		\centering 
		\subfloat[Individual effects]{
			\includegraphics[width=0.45\textwidth]{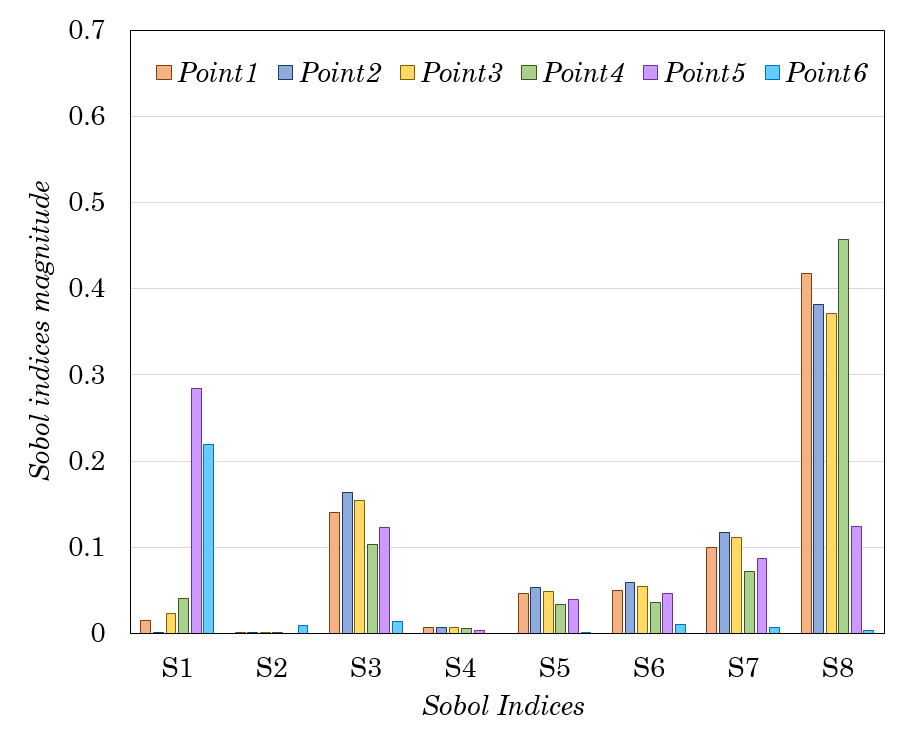}
			\label{fig:MinStress_Si_3years}
		}
		\subfloat[Interaction effects ($>0.03$)]{
			\includegraphics[width=0.45\textwidth]{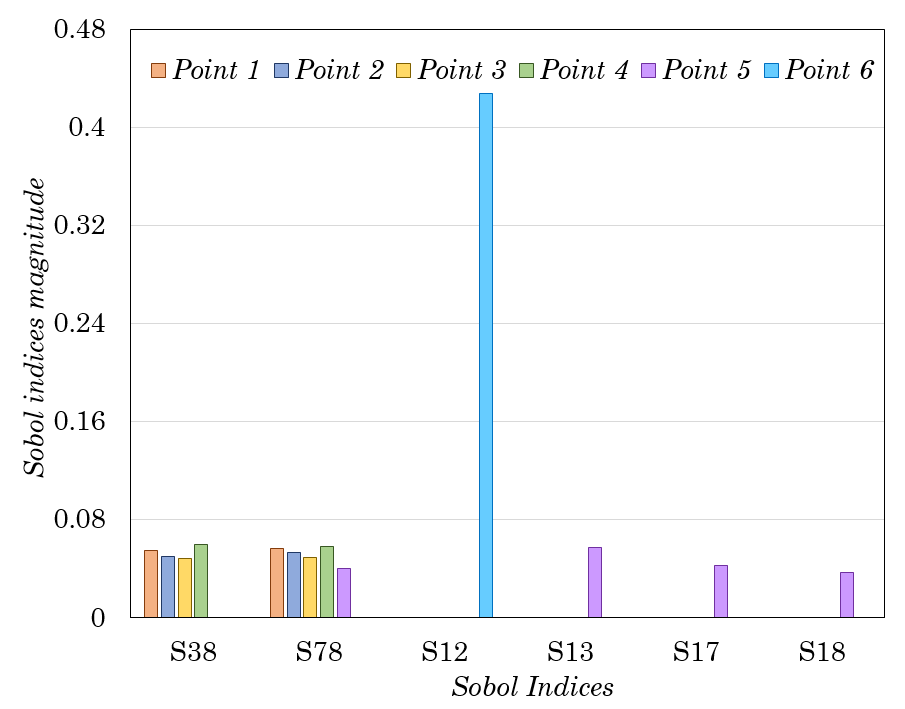}
			\label{fig:MinStress_Sij_3years}
		} 
		\caption{Sobol indices of the input sample for minimum 
		horizontal 
		stress after 3 years from production}
		\label{fig:MinStress_3years}
	\end{figure}

\begin{figure}[!h]
		\centering 
		\subfloat[Individual effects]{
			\includegraphics[width=0.45\textwidth]{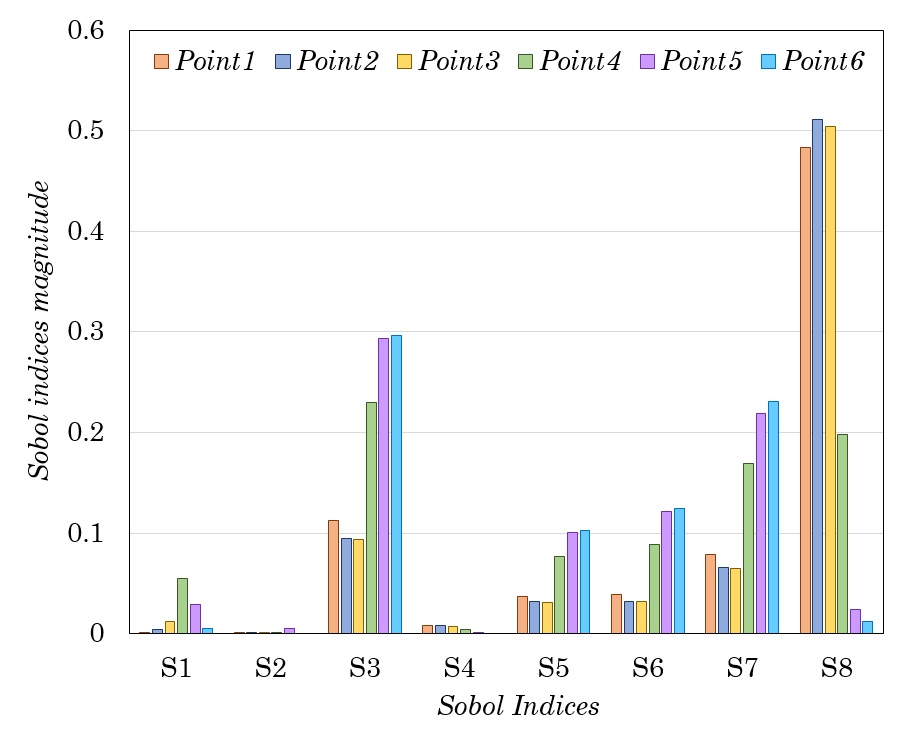}
			\label{fig:MaxStress_Si_3years}
		}
		\subfloat[Interaction effects ($>0.03$)]{
			\includegraphics[width=0.45\textwidth]{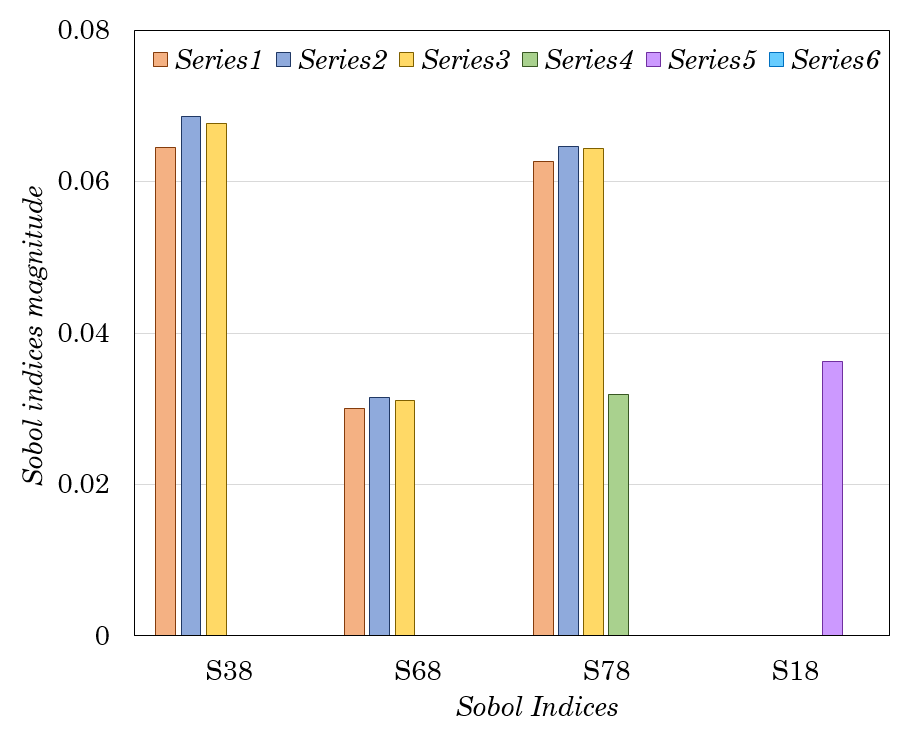}
			\label{fig:MaxStress_Sij_3years}
		} 
		\caption{Sobol indices of the input sample for maximum 
		horizontal 
		stress after 3 years from production}
		\label{fig:MaxStress_3years}
	\end{figure}

	A similar analysis is performed after three years of production. 
	Results 
	are 
	presented in Figures \ref{fig:press_3years}--\ref{fig:MaxStress_3years}. 
	Analysis of pore pressure (Figure \ref{fig:press_3years}) shows that 
	as time 
	progresses, $S_8$ decreases at all of the points. In contrary, 
	$S_3$ 
	increases 
	as time progresses. The increase of $S_{3}$ at Points $5$ and $6$ 
	is much 
	more 
	than the other points since they are closer to the 
	spacing area 
	between fractures. At Point $6$, almost $90\%$ of the variation in 
	pore 
	pressure is due to the changes in production pressure $p_f$, while 
	this 
	value 
	is close to $50\%$ at Point $5$ and less than $30\%$ in the rest of 
	points. 
	Also, fracture half-length plays an important role in the variation 
	of the 
	pore 
	pressure in the neighborhood of the tip region (Point $5$). The 
	contribution 
	of 
	fracture half-length and its interactions with other terms on the 
	variation 
	of 
	pore pressure at this point is close to $20\%$. 
	
	Figure \ref{fig:MinStress_3years} represents the Sobol indices for 
	minimum 
	horizontal stress after three years. After this period, a slight 
	decrease 
	in 
	$S_8$ is observed, while other individual indices stay almost 
	constant. Two 
	Points $5$ and $6$ are of particular interest for minimum 
	horizontal 
	stress. 
	The dominant variable that causes most of the variation in minimum 
	horizontal 
	stress at these two points is fracture half-length. Also, $S_{12}$ 
	contributes 
	to about $45\%$ of the changes at Point $6$. Among interaction 
	terms, 
	$S_{38}$ 
	and $S_{78}$ are the dominant interaction indices at all point 
	except Point 
	$6$, and $S_{12}$ at Point $6$ increases slightly compared to the 
	case of 
	$1\;year$.
	
	Figure \ref{fig:MaxStress_3years} presents the Sobol analysis results 
	for 
	maximum horizontal stress after $3\;years$. Similar to changes that 
	were 
	observed by comparing $1\;year$ and $1\;month$ results, it is 
	observed in 
	this 
	case that the effect of $S_{8}$ decreases further compared to 
	$1\;year$. 
	However, the rate of changes is less than the changes that were 
	observed 
	between two previous changes. Also, it observed that $S_{3}$ 
	increases at 
	all 
	of the points. Moreover, the dominant interaction indices that are greater 
	than 
	$0.03$ 
	have decreased from nine to four.

	\section{REDUCED ORDER MODEL FOR PORE PRESSURE AND STRESSES} 
	\label{Sec:reduced_model}
	
	In this section, we utilize one of the greatest benefits of the	Sobol	analysis 
	which is its ability to present a reduced order (mathematically simple) 
	model (ROM) for a relatively	complex function such as the one that we 
	discussed in Section \ref{sec_poroddm}. For this purpose, we consider the case 
	of production from the	hydraulically-fractured well for one year (Section 
	\ref{subsec:OneYear}) and present a ROM for for pore pressure, maximum 
	horizontal stress, and minimum horizontal stress at different points around 
	fractures. In our approach, we use	those Sobol indices that contribute to 
	$80\%$ to  $90\%$ changes in the results.
	
	In order to be avoid repetition, we group the points at which pore pressure 
	and stresses reduces similarly. It	was observed that the points located 
	outside the	spacing area	(i.e., Points 1 - 4) show the same behavior in terms 
	of 	the dominant Sobol indices, and, correspondingly,	their corresponding	
	dominant	Sobol	functions are similar. Therefore, as a	representative of these 
	points we	only present the reduced order model for Point 1. Next, we present 
	the ROM	for	Points 5 and 6 because these two pints showed different 
	behavior, and their representative Sobol functions are different. 

	\subsection{Sobol functions and ROM at Point 1}
\begin{figure}[!h]
		\centering
		\subfloat[$f_3 $]{
			\graphicspath{{figures/ReducedOrder/fi/Pp/Point1/}}
			\includegraphics[width=0.33\textwidth]{f3.png}
			\label{fig:f3@point1_Pp}
		}
		\subfloat[$f8 $]{
			\graphicspath{{figures/ReducedOrder/fi/Pp/Point1/}}
			\includegraphics[width=0.33\textwidth]{f8.png}
			\label{fig:f8@point1_Pp}
		} 
		\subfloat[$f_{38}$]{
			\graphicspath{{figures/ReducedOrder/fij/Pp/}}
			\includegraphics[width=0.33\textwidth]{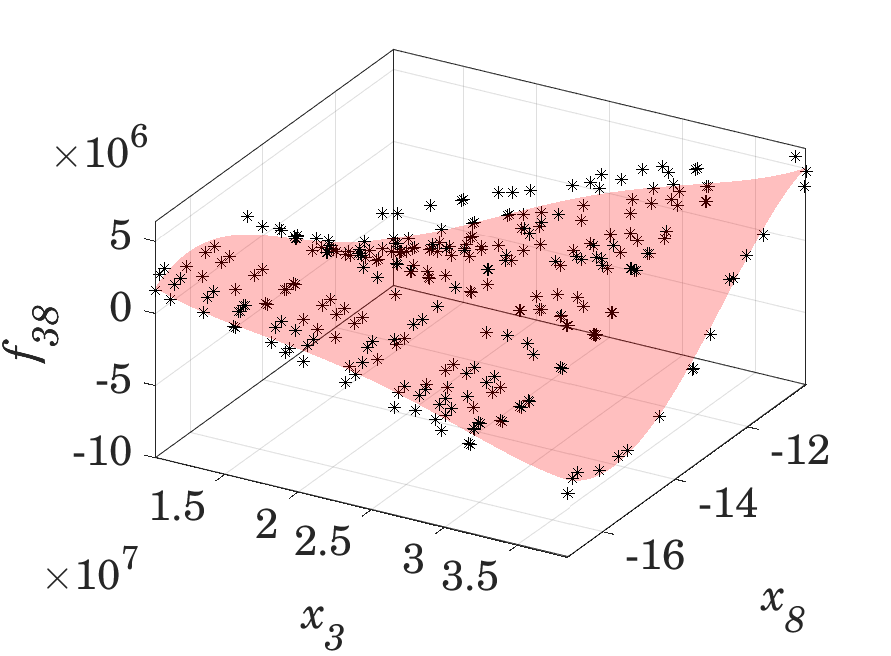}
			\label{fig:f38@point1_Pp}
		}
		\caption{The dominant Sobol functions for pore pressure at 
		Point 1.}
		\label{fig:SobolFunctions@Point1_Pp}
	\end{figure}

	\paragraph{\emph{Pore pressure.}} As shown in Figure \ref{fig:press_1year}, 
	production pressure	$P_p$, mobility $\kappa$, and their combinations (i.e., 
	$S_3, \; S_8, \;	S_{38}$) accounts for more than $90\%$	of	the pore pressure 
	changes at Points 1--4. Therefore, a reduced order model 
	(ROM) may be presented using these variables at these points. For example, 
	Equation \eqref{eq:ANOVA_rearang} for pore pressure at 
	Point 1 
	with $90\%$ accuracy may 
	be written as:
	\begin{align}
		Pp^{p_1}_{1year} = f_0 + f_3 + f_8 + f_{38}
		\label{eq:ROM_Pp_P1_summary}
	\end{align}
	Figure \ref{fig:SobolFunctions@Point1_Pp} shows these three functions 
	correspondingly. It can be seen that the pore pressure is directly related to 
	production pressure $p_f$. It is also has an inverse relationship with 
	$\log()$ of the fluid mobility. As shown in Figure \ref{fig:f8@point1_Pp}, 
	extremely low fluid mobilities (i.e., $> 10^{-13} \; \frac{m^2}{Pa.s}$), the 
	function is almost constant and changes in fluid mobility will not affect the 
	pore pressure.   
	%

\begin{figure}[H]
	 	\centering
	\graphicspath{{figures/ReducedOrder/fi/SigMin/Point1/}}
	 	\subfloat[$f_3$]{
	 		\includegraphics[width=0.33\textwidth]{f3.png}
	 		\label{fig:f3@point1_SigMin}
	 	}
	 	\subfloat[$f_5$]{
	 		\includegraphics[width=0.33\textwidth]{f5.png}
	 		\label{fig:f5@point1_SigMin}
	 	} 
	 	\subfloat[$f_6$]{
	 		\includegraphics[width=0.33\textwidth]{f6.png}
	 		\label{fig:f6@point1_SigMin}
	 	}
	 	\quad
	 	\subfloat[$f_7$]{
	 		\includegraphics[width=0.33\textwidth]{f7.png}
	 		\label{fig:f7@point1_SigMin}
	 	}
	 	\subfloat[$f_8$]{
	 		\includegraphics[width=0.33\textwidth]{f8.png}
	 		\label{fig:f8@point1_SigMin}
	 	} 
	 	\subfloat[$f_{38} $]{
	 		\graphicspath{{figures/ReducedOrder/fij/SigMin/}}
	 		\includegraphics[width=0.33\textwidth]{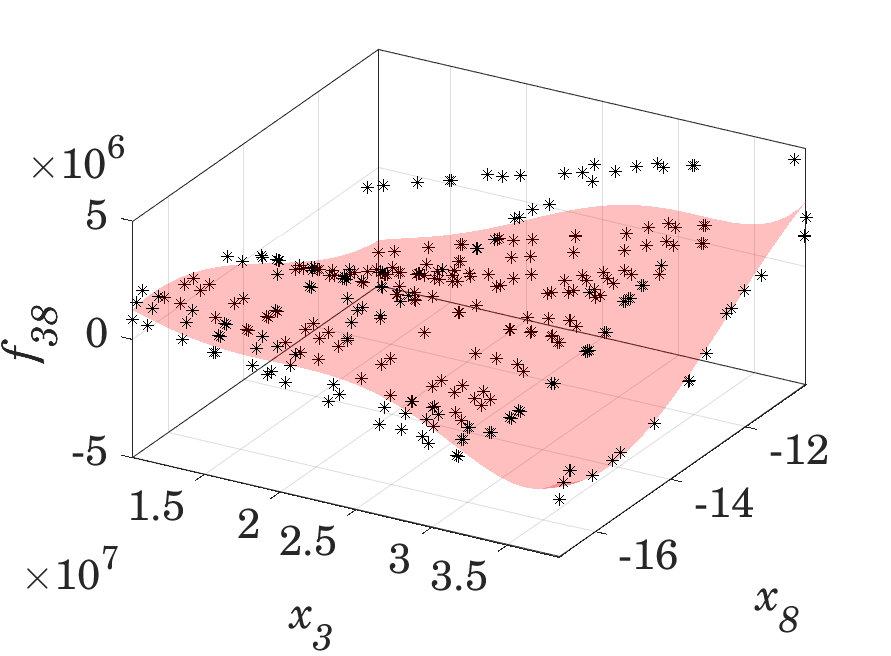}
	 		\label{fig:f38@point1_SigMin}
	 	}
	 	\quad
	 	\subfloat[$f_{78}$]{
	 		\graphicspath{{figures/ReducedOrder/fij/SigMin/}}
	 		\includegraphics[width=0.33\textwidth]{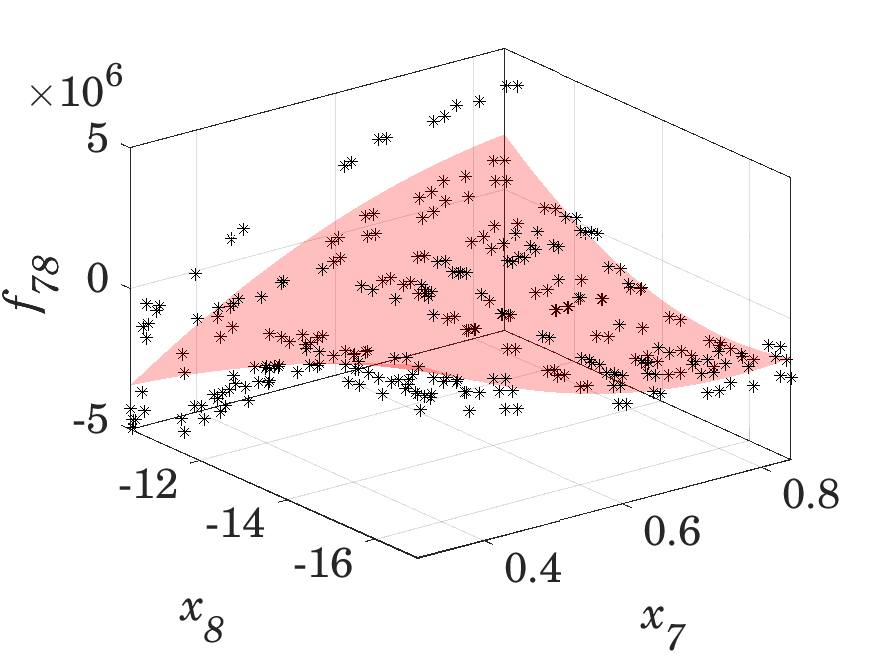}
	 		\label{fig:f78@point1_SigMin}
	 	}
	 	\caption{The dominant Sobol functions for minimum horizontal 
	 	stress at 
	 	Point 1.}
	 	\label{fig:SobolFunctions@Point1_SigMin}
	\end{figure}

	Using Equation \eqref{eq:ROM_Pp_P1_summary} and Figure 
	\ref{fig:SobolFunctions@Point1_Pp} one may construct the reduced order model 
	with the following items: 
	
	\begin{subequations}

	\begin{align}
	f_0 = & \; A_0 
		\label{eq:ROM_Pp_P1_f0} \\
	f_3 = & \; A_0 \; p_f + A_1 
	\label{eq:ROM_Pp_P1_f3} \\
	f_8 = & \; A_0 \; \sin(A_1 \; \kappa + A_2) + A_3 \; \sin(A_4 \; 
	\kappa - A_5)
	\label{eq:ROM_Pp_P1_f8}\\
	f_{38} = & - A_0 - A_1 \; p_f - A_2 \; \kappa - 
	A_3 \; p_f^2 + A_4 \; p_f \; \kappa + A_5 \; \kappa^2 + 
	A_6 \; p_f^3 - A_7 \; p_f^2 \; \kappa - A_8 \; p_f \; \kappa^2 - 
	A_9 \; \kappa^3 
	\label{eq:ROM_Pp_P1_f38}
	\end{align}
	\end{subequations}
	where the coefficients for each function is presented in Table 
	\ref{App:Tab:ROM_Coeff_Pp_P1} in the Appendix. 
	\paragraph{\emph{Minimum horizontal stress.}} Similar analysis can be done 
	for minimum horizontal stress at Point 1. As	shown in Figure 
	\ref{fig:MinStress_1year}, the dominant Sobol indices	for	minimum horizontal 
	stress at Point 1 are $S_3$, $S_5$, $S_6$, $S_7$, $S_{8}$, $S_{38}$, and 
	$S_{78}$. Figure \ref{fig:SobolFunctions@Point1_SigMin} shows th plots of 
	these 
	variables and their corresponding Sobol function for minimum horizontal stress 
	after	one	year from production. It can be observed that pore pressure, drained 
	Poisson's ratio, and Skemptson's coefficient have a directly effect on the 
	minimum horizontal stress, and undrained Poisson's ratio and mobility have an 
	inverse effect on the minimum horizontal stress at Point 1.
	For the minimum horizontal stress at Point 1, the following reduced order 
	function may be constructed:		
	\begin{subequations}
	\begin{align}
	f_0 = & \; A_0 
 	\label{eq:ROM_SigMin_P1_f0}\\
	f_3 = & \; A_0 \; p_p   + A_1 
	 \label{eq:ROM_SigMin_P1_f3}\\
	f_5 = & \; A_0 \; \nu_u + A_1 
		\label{eq:ROM_SigMin_P1_f_5}\\
	f_6 = & \; A_0 \; \nu - A_1 
		\label{eq:ROM_SigMin_P1_f6}\\
	f_7 = & \; A_0 \; B^3 - A_1 \; B^2 + A_2 \; B - A_3 
		\label{eq:ROM_SigMin_P1_f7}\\
	f_8 = & \; A_0 \; \sin(A_1 \; \kappa + A_2) + A_3 \; \sin(A_4 \; \kappa - A_5) 
	+ A_6 \; \sin(A_7 \; \kappa + A_8) 
		\label{eq:ROM_SigMin_P1_f8}\\
	f_{38} = & \; - A_0 + A_1 \; p_f - A_2 \; \kappa	- A_3 \; p_f^2 - A_4 \; p_f 
	\; \kappa + A_5 \; \kappa^2 	- A_6 \; p_f^3 + A_7 \; p_f^2 \; \kappa -A_8 \; 
	p_f \; \kappa^2 + A_9 \;	\kappa^3
		\label{eq:ROM_SigMin_P1_f38}
	\end{align}
	\end{subequations}
	
\begin{figure}[!h]
		\centering
		\subfloat[$f_{3}$]{
			\includegraphics[width=0.33\textwidth]{f3.png}
			\label{fig:f3@point1_SigMax}
		}
		\subfloat[$f_{5}$]{
			\includegraphics[width=0.33\textwidth]{f5.png}
			\label{fig:f5@point1_SigMax}
		} 
		\subfloat[$f_{6}$]{
			\includegraphics[width=0.33\textwidth]{f6.png}
			\label{fig:f6@point1_SigMax}
		}
	\quad
		\subfloat[$f_{7}$]{
			\includegraphics[width=0.33\textwidth]{f7.png}
			\label{fig:f7@point1_SigMax}
		}
		\subfloat[$f_{8}$]{
			\includegraphics[width=0.33\textwidth]{f8.png}
			\label{fig:f8@point1_SigMax}
		} 
		\subfloat[$f_{38}$]{
			\graphicspath{{figures/ReducedOrder/fij/SigMax/}}
			\includegraphics[width=0.33\textwidth]{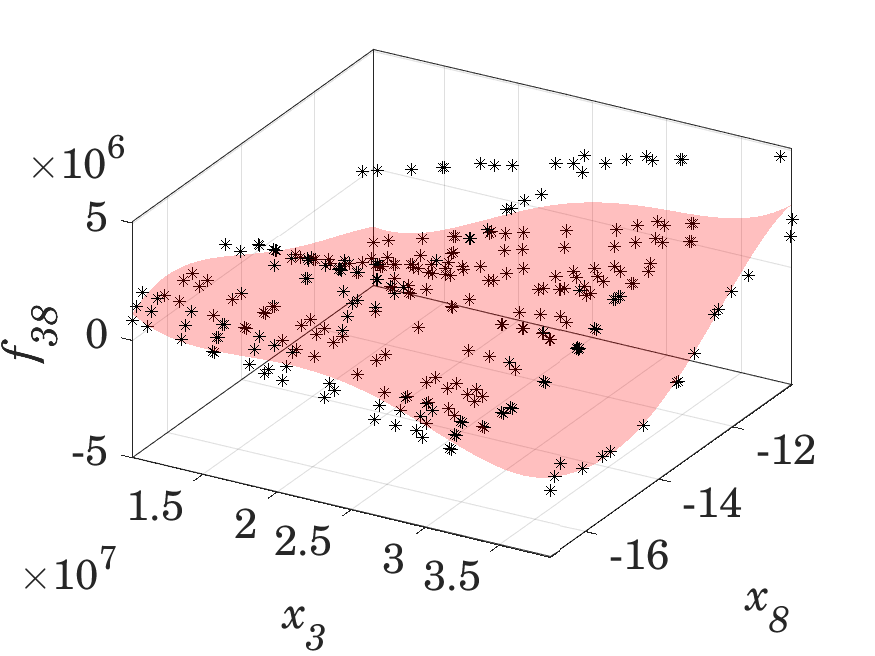}
			\label{fig:f38@point1_SigMax}
		}
		\caption{The dominant Sobol functions for maximum horizontal 
		stress at 
		Point 1}
		\label{fig:SobolFunctions@Point1}
	\end{figure}

	Coefficients of the functions in Equations 
	\eqref{eq:ROM_SigMin_P1_f0}--\eqref{eq:ROM_SigMin_P1_f38} are 
	presented in Table \ref{App:Tab:ROM_Coeff_SigMin_P1}. 
	\paragraph{\emph{Maximum horizontal stress.}} It can be verified that $f_3$, 
	$f_5$, $f_6$, $f_7$, $f_8$, and $f_{38}$ are the dominant Sobol indices that 
	affect the maximum horizontal stress at Point 6. As shown in Figure 
	\ref{fig:MaxStress_1year}, 
	these variables generate about $90\%$ of the calculated model's maximum 
	horizontal	stress. Therefore, one may construct the ROM using these 
	terms as:  
	\begin{align}
		\bar{\sigma}_{H_{1 year}}^{P_1} = f_0 + f_3 + f_5 + f_6 + f_7 
		+ f_8 + 	f_{38}
	\end{align} 
	The corresponding Sobol functions for maximum horizontal stress are plotted in 
	Figure 
	\ref{fig:SobolFunctions@Point1}. It was observed that $p_f$ and $\kappa$ show 
	the same linear 
	relationship 
	with the reduced order function as in pore pressure and minimum horizontal 
	stress. Also, similar to minimum horizontal stress, undrained Poisson's ratio 
	and mobility have effect on the minimum horizontal stress at Point 1.  
	\graphicspath{{figures/ReducedOrder/fi/SigMax/Point1/}}
	
	For maximum horizontal stress, the ROM can be represented as:
	
	\begin{subequations}
	\begin{align}
	f_0 = & \;  A_0 
	\label{eq:ROM_SigMax_P1_f0} \\
	f_3 = & \;  A_0 \; p_f - A_1 
	\label{eq:ROM_SigMax_P1_f3} \\
	f_5 = & \; -A_0 \; \nu_u + A_1 
	\label{eq:ROM_SigMax_P1_f5} \\
	f_6 = & \;  A_0 \; \nu^2 - A_1 \; \nu + A_2 
	\label{eq:ROM_SigMax_P1_f6} \\
	f_7 = & \;  A_0 \; B^3 - A_1 \; B^2 + A_2 \; B - A_3 
	\label{eq:ROM_SigMax_P1_f7} \\
	f_8 = & \;  A_0 \; \sin(A_1 \; \kappa + A_2) + 
		        A_3 \; \sin(A_4 \; \kappa - A_5) + 
		        A_6 \; \sin(A_7 \; \kappa + A_8) 
	\label{eq:ROM_SigMax_P1_f8} \\
	f_{38} = & \; - A_0 + A_1 \; p_f - A_2 \; \kappa 
	+ A_3 \; p_f \; \kappa - A_4 \; \kappa^2 + 
	A_5 \; \kappa^3 
	\label{eq:ROM_SigMax_P1_f38}
	\end{align}
	\end{subequations}	
	
	Coefficients of the function in Equations 
	\eqref{eq:ROM_SigMax_P1_f0}--\eqref{eq:ROM_SigMax_P1_f38} are 
	presented 
	in Table \ref{App:Tab:ROM_Coeff_SigMax_P1}. Similar analysis may be done for 
	other points that are located at the outside of the spacing area. Also, it 
	should be mentioned that the choice of the individual function in the reduced 
	order model is arbitrary and by the best possible match. In the next section, 
	Sobol functions and reduced order models for Point 5 is presented. 
	
	\subsection{Sobol functions and ROM at Point 5}
\begin{figure}[!h]
	\centering
	\subfloat[$f_{1}$]{
		\includegraphics[width=0.33\textwidth]{f1.png}
		\label{fig:f1@point5_Pp}
	}
	\subfloat[$f_{3}$]{
		\includegraphics[width=0.33\textwidth]{f3.png}
		\label{fig:f3@point5_Pp}
	} 
	\quad
	\subfloat[$f_{8}$]{
		\includegraphics[width=0.33\textwidth]{f8.png}
		\label{fig:f8@point5_Pp}
	}
	\subfloat[$f_{38}$]{
		\graphicspath{{figures/ReducedOrder/fij/Pp/}}
		\includegraphics[width=0.33\textwidth]{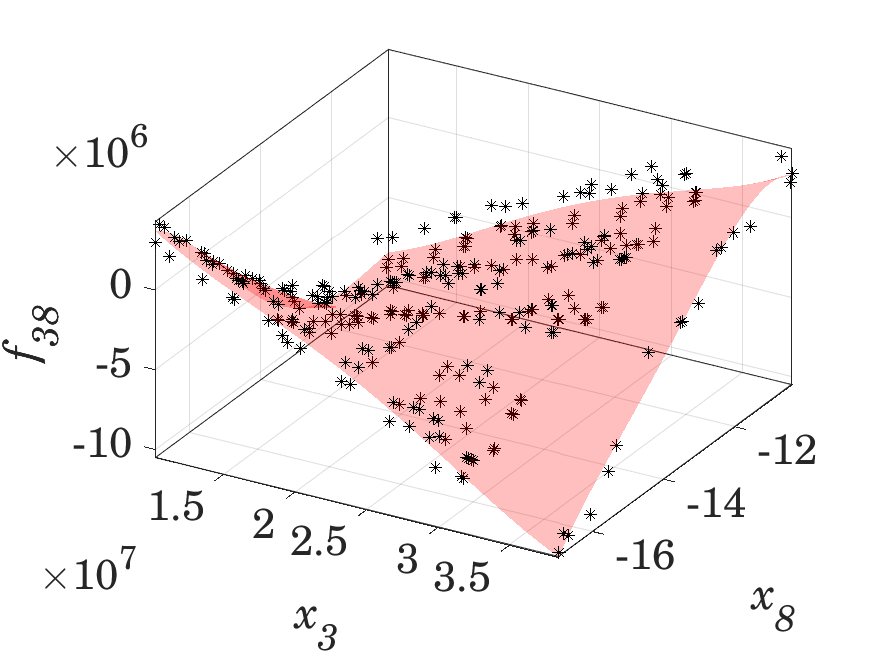}
		\label{fig:f38@point5_Pp}
	}
	\caption{Dominant Sobol functions for pore pressure at Point 5}
	\label{fig:SobolFunctions@Point5_Pp}
	\end{figure}

	\paragraph{\emph{Pore pressure.}} Point 5 is located at the point between the tips of 
	the pre-existing 
	fractures. 
	Therefore, understanding the changes in pore pressure, the maximum and minimum 
	horizontal stress at this point are extremely important 
	for refracturing and infill drilling applications. The dominant variables 
	affecting the pore pressure at Point 5 are 
	fracture 
	half-length ($a$), production pressure ($p_f$), mobility ($\kappa$) 
	and the 
	interaction between production pressure and mobility (cf.,	Figure 
	\ref{fig:press_1month}). Figure \ref{fig:SobolFunctions@Point5_Pp} 
	shows the	plot of the corresponding Sobol functions for these variables. 	
	\graphicspath{{figures/ReducedOrder/fi/Pp/Point5/}}
	
	As it can be seen, fracture half-length has a considerable 
	effect on the pore pressure at Point 5. Change of fracture half-length (with 
	the same spacing) from $57 \; m$ to $10 \; m$ will cause a reduction of $11 \; 
	Mpa$ in the pore pressure at Point 5. Using the fitted functions in Figure 
	\ref{fig:SobolFunctions@Point5_Pp}, the ROM for pore pressure at Point 5 after 
	a year from production can be constructed as follows:
	 
	\begin{subequations}
	\begin{align} 
	f_0 = & \; A_0
	\label{eq:ROM_Pp_P5_f0} \\ 
	f_1 = & \;-A_0 \; a^2 - A_1 \; a + A_2
	\label{eq:ROM_Pp_P5_f1} \\
	f_3 = & \; A_0 \; p_f - A_1
	\label{eq:ROM_Pp_P5_f3} \\ 
	f_8 = & \; A_0 \; \sin(A_1 \; \kappa - A_2) + 
	           A_3 \; \sin(A_4 \; \kappa - A_5) + 
	           A_6 \; \sin(A_7 \; \kappa + A_8)
	\label{eq:ROM_Pp_P5_f8} \\ 
	f_{38} = & \; -A_0 + A_1 \; p_f -A_2 \; \kappa -A_3 \; p_f^2 + A_4 \; p_f \; 
	\kappa + A_5 \; \kappa^2 + A_6 \; p_f^3 - A_7 \; p_f^2 \; \kappa - A_8 \; p_f 
	\; \kappa^2 + A_9 \; \kappa^3
	\label{eq:ROM_Pp_P5_f38} 
	\end{align}
	\end{subequations} 
	
\begin{figure}[!h]
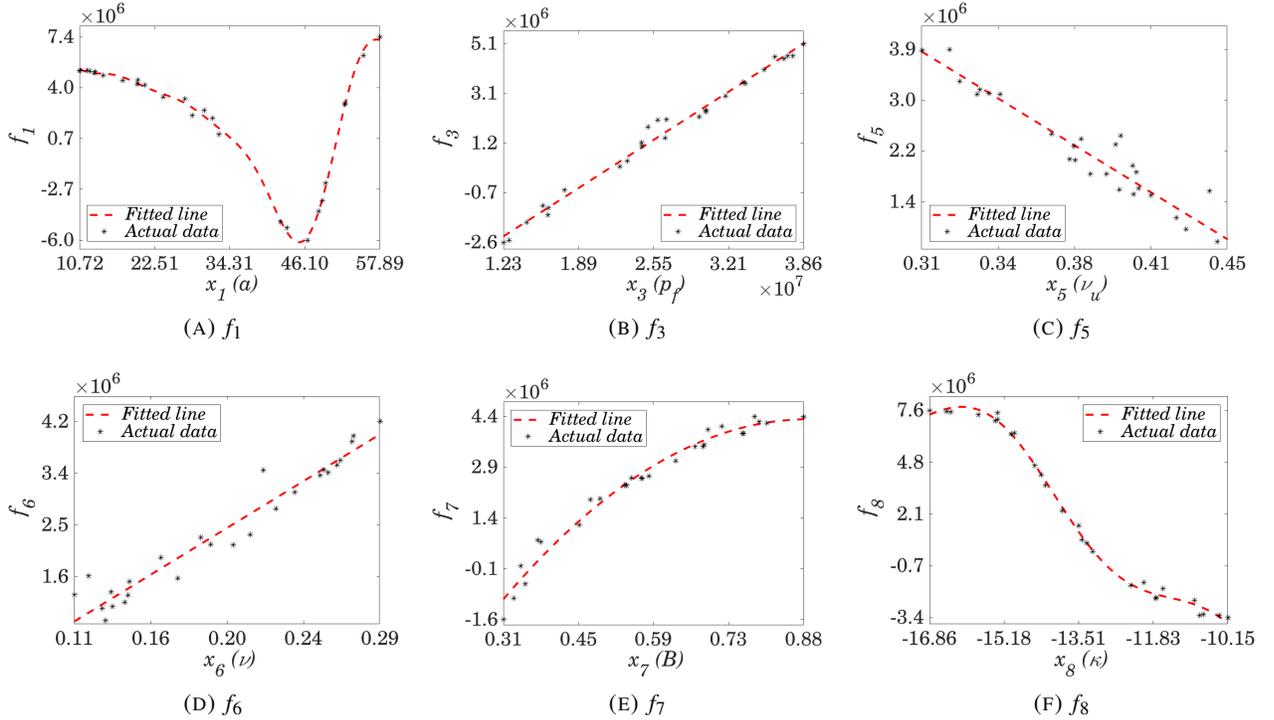

		\centering
		\subfloat[$f_{1}$]{
			\includegraphics[width=0.33\textwidth]{f1.png}
			\label{fig:f1@point5_SigMin}
		}
		\subfloat[$f_{3}$]{
			\includegraphics[width=0.33\textwidth]{f3.png}
			\label{fig:f3@point5_SigMin}
		} 
		\subfloat[$f_{5}$]{
			\includegraphics[width=0.33\textwidth]{f5.png}
			\label{fig:f5@point5_SigMin}
		}
	\quad
	\subfloat[$f_{6}$]{
	 	\includegraphics[width=0.33\textwidth]{f6.png}
	 	\label{fig:f6@point5_SigMin}
	}
	\subfloat[$f_{7}$]{
	 	\includegraphics[width=0.33\textwidth]{f7.png}
	 	\label{fig:f7@point5_SigMin}
	} 
	\subfloat[$f_{8}$]{
	 	\includegraphics[width=0.33\textwidth]{f8.png}
	 	\label{fig:f8@point5_SigMin}
	}
	\caption{Dominant first-order Sobol functions for minimum horizontal stress	at 
	Point 5}
	\label{fig:SobolFunctions@Point5_SigMin_a}
    \end{figure}

	The coefficients of the ROM's functions in Equations 
	\eqref{eq:ROM_Pp_P5_f0}--\eqref{eq:ROM_Pp_P5_f38} are 
	presented in Table \ref{App:Tab:ROM_Coeff_Pp_P5}. Among all of the function 
	that are presented for pore pressure, $f_8$ has the greatest impact on the 
	changes of the pore pressure, followed by production pressure, and the fracture 
	half-length. 
	\paragraph{\emph{Minimum horizontal stress.}} The dominant Sobol indices affecting 
	minimum	
	horizontal stress at Point 6 are $S_1$, $S_3$, $S_5$, $S_6$, $S_7$, $S_8$, 
	$S_{38}$, $S_{13}$, $S_{16}$, $S_{17}$, and $S_{18}$. Although some of these 
	variables' contributions such as $S_{13}$ is only about $8\%$, their total 
	contribution from them can capture more than $90\%$ of the model output for 
	minimum horizontal stress. Since many of the variables contribute to the 
	changes of the minimum horizontal stress at Point 5, the functions of these 
	variables are grouped into first-order function and second-order function 
	categories. Figure \ref{fig:SobolFunctions@Point5_SigMin_a} shows the plots of 
	the first-order functions for these variables. 
	
	\graphicspath{{figures/ReducedOrder/fi/SigMin/Point5/}}

\begin{figure}[!h]
		\centering
	\subfloat[$f_{38}$]{		
	 	\includegraphics[width=0.33\textwidth]{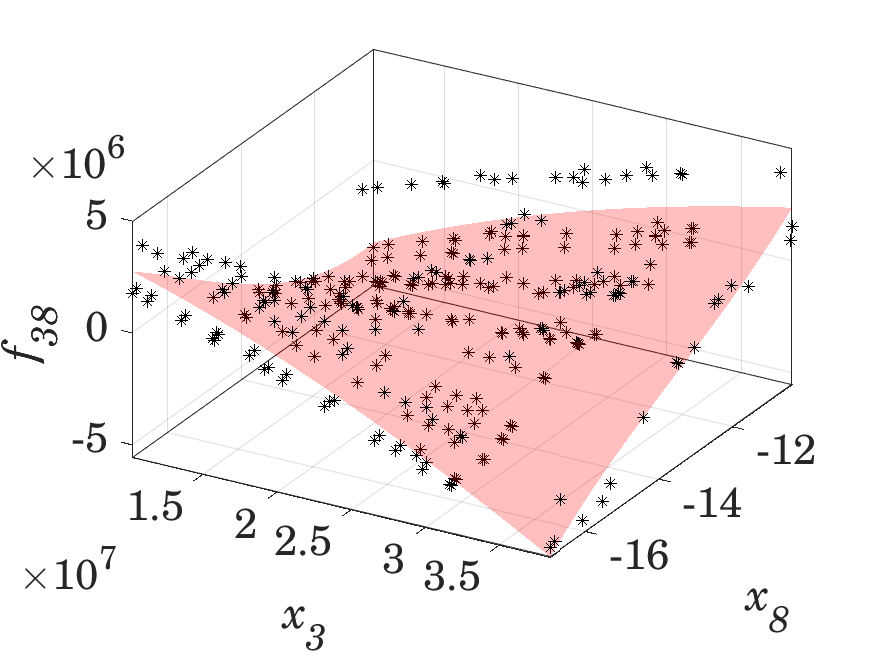}
	 	\label{fig:f38@point5_SigMin}
	}
	\subfloat[$f_{13}$]{
	 	\includegraphics[width=0.33\textwidth]{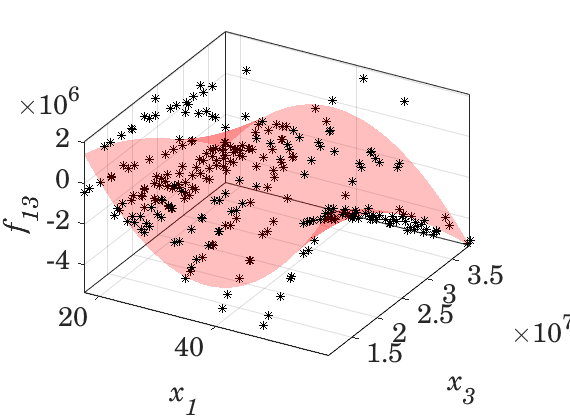}
	 	\label{fig:f13@point5_SigMin}
	} 
	\subfloat[$f_{16}$]{
	 	\includegraphics[width=0.33\textwidth]{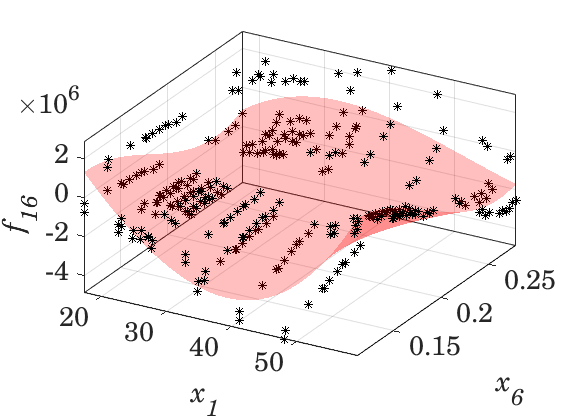}
	 	\label{fig:f16@point5_SigMin}
	}
	\quad
	\subfloat[$f_{17}$]{
	 	\includegraphics[width=0.33\textwidth]{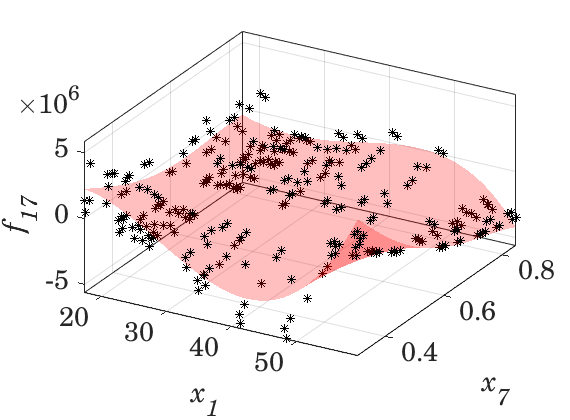}
	 	\label{fig:f17@point5_SigMin}
	}
	\subfloat[$f_{18}$]{
	 	\includegraphics[width=0.33\textwidth]{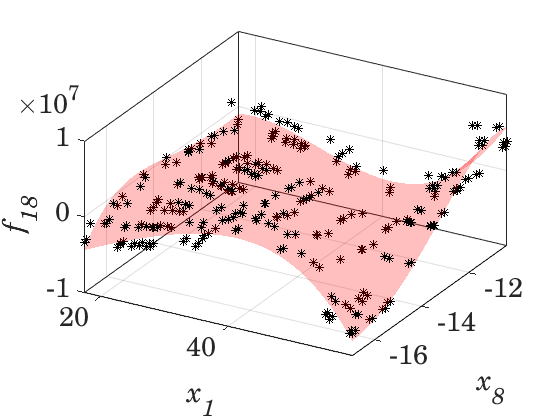}
	 	\label{fig:f18@point5_SigMin}
	} 
		\caption{Dominant second-order Sobol functions for minimum horizontal stress	
		at Point 5}
		\label{fig:SobolFunctions@Point5_SigMin_b}
	\end{figure}

	As shown in the figure, as fracture half-length increases, it initially has an 
	inverse effect on the minimum horizontal stress up to certain point. After 
	passing that point, its effect is directly proportional to the minimum 
	horizontal stress.
	Figure \ref{fig:SobolFunctions@Point5_SigMin_b} represents the dominant 
	second-order term functions. As shown in the figure, one of the individual 
	terms in most of these functions is fracture half-length. Fracture half-length 
	accounts for $50\%$ of the changes in minimum horizontal stress at Point 5.
	\graphicspath{{figures/ReducedOrder/fij/SigMin/}}
	
	Among the dominant variables, production pressure and fracture half-length have 
	a considerable impact on the minimum horizontal stress at this point. Using the 
	functions in Figures \ref{fig:SobolFunctions@Point5_SigMin_a} and 
	\ref{fig:SobolFunctions@Point5_SigMin_b}, the ROM for minimum horizontal 
	stress at Point 5 can be presented as follows:
	
	\begin{subequations}
	\begin{align}
	f_0 = & \; A_0
	\label{eq:ROM_SigMin_P5_f0}\\
	f_1 = & \; A_1 - A_2 \; \cos(A_0 \; a) + A_3 \; \sin(A_0 \; a) + 
	           A_4 \; \cos(2 \; A_0 \; a) + A_5 \; \sin(2 \; A_0 \; a) + 
	           A_6 \; \cos(3 \; A_0 \; a) + \nonumber \\ &  A_7 \; \sin(3 \; A_0 
	           \; 
	           a) + 
	           A_8 \; \cos(4 \; A_0 \; a) + A_9 \; \sin(4 \; A_0 \; a) 
	\label{eq:ROM_SigMin_P5_f1}\\
	f_3 = & \;  A_0 \; p_f - A_1 
	\label{eq:ROM_SigMin_P5_f3}\\
	f_5 = & \; -A_0 \; \nu_u + A_1 
	\label{eq:ROM_SigMin_P5_f5}\\
	f_6 = & \;  A_0 \; \nu - A_1 
	\label{eq:ROM_SigMin_P5_f6}\\
	f_7 = & \; -A_0 \; B^2 + A_1 \; B - A_2
	\label{eq:ROM_SigMin_P5_f7}\\
	f_8 = & \;  A_0 \; \sin(A_1 \; \kappa + A_2) + 
	            A_3 \; \sin(A_4 \; \kappa + A_5)
	\label{eq:ROM_SigMin_P5_f8}\\
	f_{38} = & \; - A_0 + A_1 \; p_f - A_2 \; \kappa - A_3 \; p_f^2 + 
	             A_4 \; p_f \kappa + A_5 \; \kappa^2 + A_6 \; p_f^3 + 
	           A_7 \; p_f^2 \; \kappa - A_8 \; p_f \; \kappa^2 + A_9 \; \kappa^3 
	\label{eq:ROM_SigMin_P5_f38}\\
	f_{13} = & - A_0 - A_1 \; a + A_2 \; p_f + A_3 \; a^2 + A_4 \; a \; p_f 
	       +   - A_5 \; p_f^2 + A_6 \; a^3 - A-7 \; a^2 \; p_f - A_8 \; a \; p_f^2 
	        \nonumber \\ & - A_9 \; p_f^3 - A_{10} \;	a^4	A_{11} \; a^3 \; p_f + 
	        A_{12} \; 
	        a^2 \; p_f^2 - A_{13} \; a \; p_f^3 + A_14 \; p_f^4  
	\label{eq:ROM_SigMin_P5_f13}\\
	f_{16} = & A_0 - A_1 \; a + A_2 \; \nu + A_3 \;	a^2 + A_4 \; a \; \nu - A_5 \; 
	      \nu^2 + A_6 \; a^3 - A_7 \; a^2 \; \nu - A_8 \; a \; \nu^2 + A_9 \; \nu^3
	\label{eq:ROM_SigMin_P5_16}\\
	f_{17} = & -A_0 - A_1 \;a + A_2 \;B + A_3 \; a^2 + A_4 \;a \;B + A_5 \; B^2 
	         + A_6 \; a^3 - A_7	\; a^2 \; B	- A_8 \;a \; B^2 - A_9 \; B^3 
	         \nonumber \\ &  - 
	         A_{10}\; a^4 - A_{11} \; a^3 \; B + A_{12} \; a^2 \; B^2 - A_{13} \; 
	         a \; B^3 - A_{14} \; B^4	
	\label{eq:ROM_SigMin_P5_17}\\	
	f_{18} = & - A_0 - A_1 \; a - A_2 \; \kappa + A_3 \; a^2 + 
	            A_4 \; a \; \kappa + A_5 \; \kappa^2 + A_6 \; a^3 +	A_7 \; a^2 \; 
	            \kappa + A_8 \; a \; \kappa^2 + A_9 \; \kappa^3 - \nonumber \\ & 
	            A_{10} \; a^4 + 
            	 A_{11} \; a^3 \; \kappa -A_{12} \; a^2 \; \kappa^2	
            	- A_{13} 
            	\; a \; \kappa^3 - A_{14} \; \kappa^4
 \label{eq:ROM_SigMin_P5_f18} 
	\end{align}
	\end{subequations}
	
	The coefficients of the functions in Equations 
	\eqref{eq:ROM_SigMin_P5_f0}--\eqref{eq:ROM_SigMin_P5_f18} are 
	presented in Tables \ref{App:Tab:ROM_Coeff_SigMin_P5_a} and 
	\ref{App:Tab:ROM_Coeff_SigMin_P5_b}. Depending on the desired degree of 
	confidence, the 
	ROM in Equations \eqref{eq:ROM_SigMin_P5_f0}--\eqref{eq:ROM_SigMin_P5_f18} may 
	be truncated to have smaller 
	number of variables.
	
\begin{figure}[!h]
		\centering
		\subfloat[$f3$]{
			\includegraphics[width=0.33\textwidth]{f3.png}
			\label{fig:f3@point5_SigMax}
		}
		\subfloat[$f5$]{
			\includegraphics[width=0.33\textwidth]{f5.png}
			\label{fig:f5@point5_SigMax}
		} 
		\subfloat[$f6$]{
			\includegraphics[width=0.33\textwidth]{f6.png}
			\label{fig:f6@point5_SigMax}
		}
	\quad
		\subfloat[$f7$]{
			\includegraphics[width=0.33\textwidth]{f7.png}
			\label{fig:f7@point5_SigMax}
		} 
		\subfloat[$f8$]{
			\includegraphics[width=0.33\textwidth]{f8.png}
			\label{fig:f8@point5_SigMax}
		}
	 	\subfloat[$f_{18}$]{
	\graphicspath{{figures/ReducedOrder/fij/SigMax/}}
	 	\includegraphics[width=0.33\textwidth]{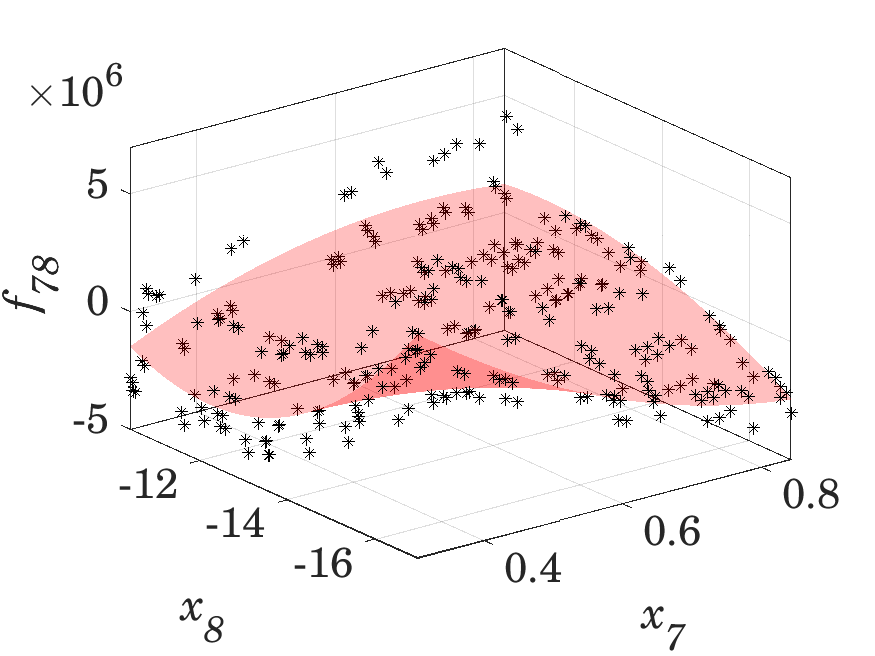}
	 	\label{fig:f18@point5_SigMax}
	} 
		\caption{Dominant Sobol functions for maximum horizontal stress at Point 5}
		\label{fig:SobolFunctions@Point5_SigMax}
	\end{figure}

	\paragraph{\emph{Maximum horizontal stress.}} The dominant affecting variables to the 
	changes in maximum horizontal stress at point 5 are $p_f$, $\nu_u$, 
	$\nu$, $B$, $\kappa$, and $f_{18}$. Figure 
	\ref{fig:SobolFunctions@Point5_SigMax} shows the plot of these variables and 
	their corresponding Sobol function.
	\graphicspath{{figures/ReducedOrder/fi/SigMax/Point5/}}
	
	Using these functions, the following reduced order model may be constructed 
	for maximum horizontal stress:
	
	\begin{subequations}
	\begin{align}
	f_0 = & \; A_0 
	\label{eq:ROM_SigMax_P5_f0} \\
	f_3 = & \; A_0 \; p_f   - A_1
	\label{eq:ROM_SigMax_P5_f3} \\ 
	f_5 = & \;-A_0 \; \nu_u + A_1
	\label{eq:ROM_SigMax_P5_f_5} \\
	f_6 = & \; A_0 \; \nu - A_1
	\label{eq:ROM_SigMax_P5_f6} \\
	f_7 = & \; A_0 \; B^3 - A_1 \; B^2 + A_2 \;	B - A_3
	\label{eq:ROM_SigMax_P5_f7} \\  
	f_8 = & \; A_1 - A_2 \; \cos(A_0 \; \kappa) + A_3 \; \sin(A_0 \; \kappa) + 
	A_4  
	\; \cos(2 \; A_0 \; \kappa) - A_5 \; \sin(2	\; A_0 \; \kappa)
	\label{eq:ROM_SigMax_P5_f8} \\
	f_{18} = & \; -A_0 -A_1 \; a + A_2 \; \kappa + 
	A_3 \; a^2 + 
	A_4 \; a \; \kappa + A_5 \; \kappa^2 + A_6 \; 
	a^3 + A_7 \; a^2 \; \kappa + A_8 \; a \; \kappa^2 - A_9 \; 
	\kappa^3 - A_{10} \; a^4 \nonumber \\ & - A_{11} \; a^3 \; \kappa + A{12} \; 
	a^2 \; 
	\kappa^2 
	-A_{13} \; a \; \kappa^3 -A_{14} \; \kappa^4
	\label{eq:ROM_SigMax_P5_f18} 
	\end{align}
	\end{subequations}
	
	Coefficients of the functions in Equations 
	\eqref{eq:ROM_SigMax_P5_f0}--\eqref{eq:ROM_SigMax_P5_f18} are
	presented in Table \ref{App:Tab:ROM_Coeff_SigMax_P5}. Next point that we 
	investigate in this paper is Point 6 which is located close to horizontal well 
	and in the middle of pre-existing fracture's spacing.
		
	\subsection{Sobol functions and ROM at Point 6}
	Point 6 is the most important point among all other points for refracturing 
	process because it is	the	point were the refracture will be placed. Therefore, 
	it is important to	keep	track of the changes in pore pressure, maximum 
	horizontal, and minimum	horizontal stresses at this point. 
	
\begin{figure}[!h]
		\centering
	
		\subfloat[$f_{3}$]{
			\includegraphics[width=0.33\textwidth]{f3.png}
			\label{fig:f2@point6_Pp}
		} 
		\subfloat[$f_{8}$]{
			\includegraphics[width=0.33\textwidth]{f8.png}
			\label{fig:f3@point6_Pp}
		}
	\quad \graphicspath{{figures/ReducedOrder/fij/Pp/}}
	\subfloat[$f_{38}$]{
	 	\includegraphics[width=0.33\textwidth]{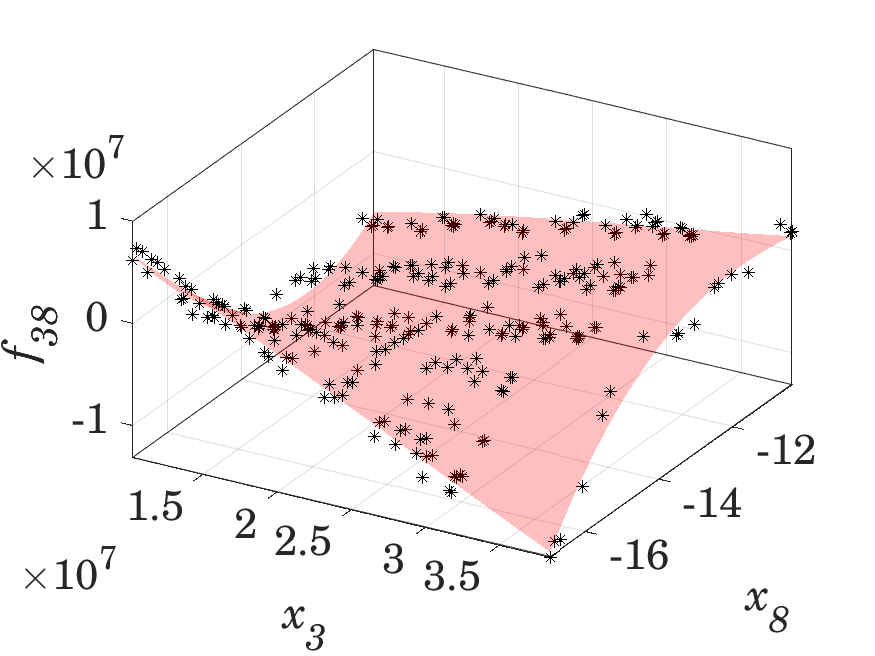}
	 	\label{fig:f38@point6_Pp}
	} 
	\subfloat[$f_{48}$]{
	 	\includegraphics[width=0.33\textwidth]{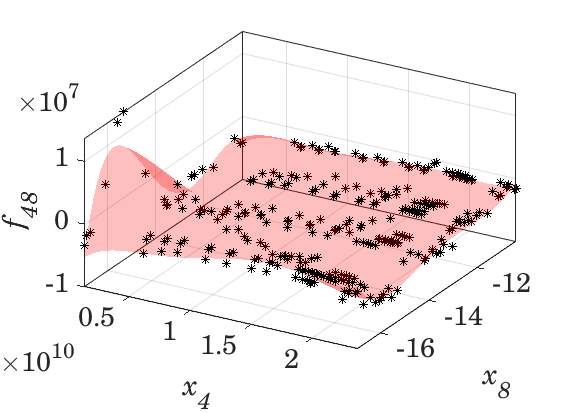}
	 	\label{fig:f47@point6_Pp}
	}
		\caption{Dominant Sobol functions for pore pressure at Point 6}
		\label{fig:SobolFunctions@Point6_Pp}
	\end{figure}

	\paragraph{\emph{Pore pressure.}} There are four main Sobol indices that control the 
	changes of the pore	pressure	at Point 6. These inputs are $S_3$, $S_8$ 
	$S_{38}$, and $S_{48}$. Figure \ref{fig:SobolFunctions@Point6_Pp} shows the 
	dominant Sobol functions	for pore	pressure at Point 6.

	\graphicspath{{figures/ReducedOrder/fi/Pp/Point6/}}
	
	The reduced order model for pore pressure at Point 6 can be represented 
	as:
	
	\begin{subequations}
	\begin{align}
	f_0 = & \; A_0 
	\label{eq:ROM_Pp_P6_f0} \\
	f_1 = & \; A_0 \; p_f - A_1 
	\label{eq:ROM_Pp_P6_f1} \\
	f_8 = & \; A_0 \; \sin( A_1 \; \kappa + A_2) + 
	           A_3	\;	\sin( A_4 \; \kappa - A_5)
	\label{eq:ROM_Pp_P6_f_8} \\
	f_{38} = & \; - A_0 + A_1 \; a - A_2 \; \kappa -	A_3 \; a^2 + A_4 \; a \; 
	           \kappa + A_5 \;	\kappa^2 +	A_6 \; a^3 - 	A_7 \; a^2 \; \kappa - A_8 
	           \; a \; \kappa^2 +	A_9 \;	\kappa^3 
	\label{eq:ROM_Pp_P6_f38} \\           
	f_{48} = & \; - A_0 + A_1 \; G + A_2 \; \kappa -	A_3 \; G^2 + A_4 \; G \; 
	         \kappa + A_5 \;	\kappa^2 -	A_6 \; G^3  -	A_7 \; G^2 \; \kappa - 
	         A_8 \; 
	         G \; 	\kappa^2 \\& - A_{9}	\;	\kappa^3 + A_{10} \; G^4 + A_{11} \; 
	         G^3 \; 
	         \kappa -	A_{12}	\;	G^2 \; \kappa^2	- A_{13} \; G \; \kappa^3 - 
	         A_{14} \; \kappa^4 +	A_{15}	\; G^5	 \\ & -	A_{16} \; G^4 \; \kappa + 
	         A_{17} 
	         \; G^3 \; \kappa^2 	+ A_{18} \; G^2 \; \kappa^3 + A_{19} \; G \; 
	         \kappa^4 -	A_{20}	\; \kappa^5
	\label{eq:ROM_Pp_P6_f48} 
	\end{align}
	\end{subequations}
	
	Coefficients of the functions in Equations 
	\eqref{eq:ROM_Pp_P6_f0}--\eqref{eq:ROM_Pp_P6_f48} are presented 
	in Table \ref{App:Tab:ROM_Coeff_Pp_P6}. 
	
\begin{figure}[!h]
		\centering
		\subfloat[$f_{1}$]{
			\includegraphics[width=0.5\textwidth]{f1.png}
			\label{fig:f1@point6_SigMin}
		} \\
		\subfloat[$f_{12}$]{
			\graphicspath{{figures/ReducedOrder/fij/SigMin/}}
			\includegraphics[width=0.5\textwidth]{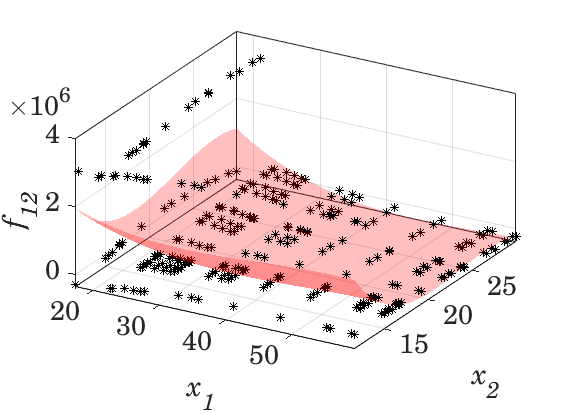}
			\label{fig:f8@point6_SigMin}
		}
		\caption{Dominant Sobol functions for minimum horizontal stress at Point	6.}
		\label{fig:SobolFunctions@Point6_SigMin}
	\end{figure}

	\paragraph{\emph{Minimum horizontal stress.}} There are two dominant input variables 
	that contribute to more than $65\%$ of the changes in minimum horizontal 
	stress at Point 6. These two variables are fracture half-length and the 
	interaction effect of the fracture half-length and fracture spacing. About half
	of the changes in the minimum horizontal stress is due to changes of these two 
	variables (i.e., combination effect). Figure 
	\ref{fig:SobolFunctions@Point6_SigMin} shows the plots of these two variables 
	and their corresponding functions. 
	
	\graphicspath{{figures/ReducedOrder/fi/SigMin/Point6/}}

 Based on the plots of the dominant Sobol function for minimum horizontal 
 stress at Point 6, the reduced order model for $\sigma_h$ at this point may be 
 obtained using: 
	
	\begin{subequations}
	\begin{align}
	f_0 = & \; 5.53 \times 10^7
	\label{eq:ROM_SigMin_P6_f0} \\
	f_1 = & \; - A_0 \; a^8 + A_1 \; a^7 - A_2 \; a^6 + A_3 \; a^5 + 
	         A_4 \; a^4 - A_5 \; a^3 - A_6 \; a^2 - A_7 \; a - A_8 
	\label{eq:ROM_SigMin_P6_f1}\\        
	f_{12} = & \; - A_0 + A_1 \; a - A_2 \; b +	A_3	\;	a^2 - A_4 \; a \; b + A_5 
	\; b^2	- A_6 \; a^3 + A_7 \; a^2 \; b + A_8 \; a \;	b^2 -	A_9 \; b^3
	\label{eq:ROM_SigMin_P6_f12}  
	\end{align}
	\end{subequations}

Coefficients of the functions in Equations 
\eqref{eq:ROM_SigMin_P6_f0}--\eqref{eq:ROM_SigMin_P6_f12} are 
presented in Table \ref{App:Tab:ROM_Coeff_SigMin_P6}	
	
\begin{figure}[H]
		\centering
	
		\subfloat[$f_{1}$]{
			\includegraphics[width=0.33\textwidth]{f1.png}
			\label{fig:f1@point6_SigMax}
		} 
		\subfloat[$f_{3}$]{
			\includegraphics[width=0.33\textwidth]{f3.png}
			\label{fig:f3@point6_SigMax}
	} 
	\subfloat[$f_{5}$]{
	 	\includegraphics[width=0.33\textwidth]{f5.png}
	 	\label{fig:f5@point6_SigMax}
	}
	\quad
	\subfloat[$f_{6}$]{
	 	\includegraphics[width=0.33\textwidth]{f6.png}
	 	\label{fig:f6@point6_SigMax}
	}
	\subfloat[$f_{7}$]{
	 	\includegraphics[width=0.33\textwidth]{f7.png}
	 	\label{fig:f7@point6_SigMax}
	}
	\subfloat[$f_{8}$]{
	 	\includegraphics[width=0.33\textwidth]{f8.png}
	 	\label{fig:f8@point6_SigMax}
	}
	\quad
		\subfloat[$f_{18}$]{
			\graphicspath{{figures/ReducedOrder/fij/SigMax/}}
			\includegraphics[width=0.33\textwidth]{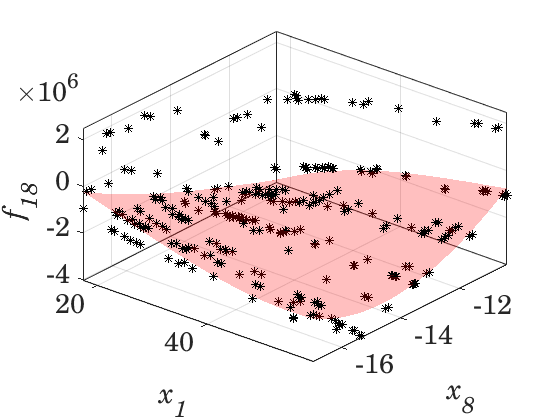}
			\label{fig:f18@point6_SigMax}
		}
		\caption{Dominant Sobol functions for maximum horizontal stress at Point 6}
		\label{fig:SobolFunctions@Point6_SigMax}
	\end{figure}

\paragraph{\emph{Maximum horizontal stress.}} There are seven dominant contributors to 
the changes in maximum horizontal stress at Point 6. These contributors are 
fracture half-length, production pressure, undrained Poisson's ratio, drained 
Poisson's ratio, Skemptson's coefficient, mobility, and interaction of fracture 
half-length and mobility. Figure \ref{fig:SobolFunctions@Point6_SigMax} shows 
plots of these contributors and their corresponding Sobol functions.
	
	\graphicspath{{figures/ReducedOrder/fi/SigMax/Point6/}}
	
	Using the functions in Figure \ref{fig:SobolFunctions@Point6_SigMax}, the 
	following ROM can be used to calculate the maximum horizontal stress at Point 
	6:	
	
	\begin{subequations}
	\begin{align}
	f_0 = & \;  A_0 
	\label{eq:ROM_SigMax_P6_f0} \\
	f_1 = & \; -A_0 \; a^2 + A_1 \; a - A_2	
	\label{eq:ROM_SigMax_P6_f1} \\
	f_3 = & \;  A_0 \; p_f - A_1 
	\label{eq:ROM_SigMax_P6_f3} \\
	f_5 = & \; -A_0 \; \nu_u + A_1 
	\label{eq:ROM_SigMax_P6_f5} \\
	f_6 = & \;  A_0 \; \nu - A_1 
	\label{eq:ROM_SigMax_P6_f6} \\
	f_7 = & \; -A_0 \; B^2 + A_1 \; B - A_2
	\label{eq:ROM_SigMax_P6_f7} \\ 
	f_8 = & \; A_0 \; \sin(A_1 \; \kappa + A_2) + A_3 \;	\sin(A_4 \; \kappa - A_5) 
	+ A_6 \; \sin(A_7 \; \kappa +	A_8)	
	\label{eq:ROM_SigMax_P6_f8} \\
	f_{18} = & \;  - A_0 - A_1 \; a - A_2 \; \kappa -	A_3 \; a^2 + A_4 \; a \; 
	\kappa + A_5 \;	\kappa^2 +	A_6 \; a^3 -	A_7 \; a^2 \; \kappa + A_8 \; a \; 
	\kappa^2 - A_9 \;	\kappa^3 
	\label{eq:ROM_SigMax_P6_f18}  
	\end{align}
	\end{subequations}
	
	The coefficients of the functions in Equations 
	\eqref{eq:ROM_SigMax_P6_f0}--\eqref{eq:ROM_SigMax_P6_f18}are 
	presented in Table \ref{App:Tab:ROM_Coeff_SigMin_P6}. In this section, 
	different reduced order models were presented for calculating the pore 
	pressure, maximum horizontal stress and minimum horizontal stress at different 
	points around a hydraulically-fractured horizontal well. Presented ROMs are 
	valid for the case of one year production with constant production pressure. 
	for any other point or a different time, a set of new ROM may be constructed 
	using the same procedure. It also showed that depending on the location of the 
	desired point with respect to the horizontal well and pre-existing hydraulic 
	fractures, different set of dominant variables contribute to the changes of 
	the desired quantity (i.e., pore pressure or stresses). As an example, we 
	presented the ROM for the points that are located outside the spacing, at 
	the tip in the hypothetical line connecting the tips of fractures, and in the 
	middle line between fractures.

	\section{CONCLUDING REMARKS}
	\label{Sec:Sobol_CR}
	
	In this study, using a fully-coupled geomechanical model, it was 
	shown that 
	different hydraulically-fractured rocks have dissimilar 
	pore 
	pressure depletion and stress change under the same boundary 
	conditions. 
	This 
	dissimilarity is illustrated using two parallel hydraulic fractures after 
	varying 
	production time periods. Results show that, besides the different 
	behaviors 
	that were observed for different rocks, the changes in the 
	mentioned 
	variables were different from point to point around hydraulic 
	fractures for 
	the same rock type. To find the most influencing inputs on the 
	model 
	outputs, 
	a global sensitivity analysis based on Sobol method is used. We 
	chose eight 
	parameters as our set of input parameters. Pore pressure, maximum 
	horizontal 
	stress, and minimum horizontal stress were chosen as the quantities of 
	interests.
	
	The results showed that mobility $\kappa$ and production pressure 
	$p_f$ 
	(which 
	needs to be seen as $\Delta p$ between reservoir and fracture) and 
	their 
	interactions are the dominant properties that cause most of the 
	pore 
	pressure 
	changes. These two variables are also dominant for the changes in 
	minimum 
	and 
	maximum horizontal stresses. It was observed that as the point of 
	observation 
	gets closer to the spacing between fractures, fracture half-length 
	and its 
	interactions contribute to $20\%$ of the changes in pore pressure 
	depletion. 
	Moreover, it was observed that interaction between fracture 
	half-length and 
	fracture spacing ($S_{12}$) has the most significant impact on the 
	minimum 
	horizontal stress at a point inside the spacing. Thus, 
	selecting 
	these 
	two variables appropriately will increase the chance of success in 
	operations 
	such as refracturing. 
	
	The following inferences are valuable for 
	practical applications. Different input variables 
	are dominant responsible elements for variation 
	of different output variables. The dominant factors 
	are not the same with time and 
	location. Thus, in order to perform any further operation around a 
	hydraulically-fractured well-bore, time and location have to be 
	considered carefully. The Sobol method offers a nice 
	framework for a systematic parametric study, as it can capture not only 
	the influence of individual parameters but also the 
	influence of interactions among the parameters. 
	
	Finally, we suggest three plausible future works. The 
	first research effort can be towards combining fracture 
	modeling with double porosity/permeability models 
	\citep{nakshatrala2018modeling,joodat2018modeling}. 
	The second effort is to incorporate inertial effects (e.g., 
	Forchheimer-type models) and pressure-dependence 
	viscosity \citep{chang2017modification,mapakshi2018scalable} 
	on the hydraulic fracture propagation. 
	The third effort can be towards utilizing phase modeling 
	for hydraulic fracture propagation on the lines similar 
	to \citep{miehe2010thermodynamically}.

	\bibliographystyle{spbasic}  
	\bibliography{biblio} 
	
\newpage
\appendix
\section{COEFFICIENTS OF THE REDUCED ORDER MODELS}

	\begin{table}[h]
	\centering
	\caption{Coefficients of the reduced order model for pore pressure 
	at Point 
	1 (Equations \eqref{eq:ROM_Pp_P1_f0}--\eqref{eq:ROM_Pp_P1_f38})}
	\label{App:Tab:ROM_Coeff_Pp_P1}  
	\begin{tabular}{lllll}
		\hline\noalign{\smallskip}
		Coeff & $f_0$ & $f_3$ & $f_8$ & $f_{38}$ \\ 
		\noalign{\smallskip}\hline\noalign{\smallskip}
	$A_0$ & $4.01 \times 10^7$ & $0.2348$ & $8.9 \times 10^6$ &	$1.161	\times 
	10^6$ \\
	$A_1$ & $-$ & $2.93 \times 10^6$ & $0.43$ & $4.734 \times 
	10^5$  \\
	$A_2$ & $-$ & $-$ & $2.165$ & $6718$ \\
	$A_3$ & $-$ & $-$ & $7.8 \times 10^5$ & $1.856 \times 10^5$ \\
	$A_4$ & $-$ & $-$ & $1.7$ & $1.995 \times 10^6$ \\
	$A_5$ & $-$ & $-$ & $0.14$ & $1.667 \times 10^5$ \\
	$A_6$ & $-$ & $-$ & $-$ & $1408$ \\
	$A_7$ & $-$ & $-$ & $-$ & $1664$ \\
	$A_8$ & $-$ & $-$ & $-$ & $5.835 \times 10^5$ \\
	$A_9$ & $-$ & $-$ & $-$ & $4.631 \times 10^4$ \\    
		\noalign{\smallskip}\hline
	\end{tabular}
	\end{table}

	\begin{table}[!h]
	\centering
	\caption{Coefficients of the reduced order model for the minimum 
	horizontal 
	stress 
	at Point 1 (Equations 
	\eqref{eq:ROM_SigMin_P1_f0}--\eqref{eq:ROM_SigMin_P1_f38})}
	\label{App:Tab:ROM_Coeff_SigMin_P1}  
	\begin{tabular}{llllllll}
		\hline\noalign{\smallskip}
	Coeff & $f_0$ & $f_3$ & $f_5$ & $f_6$ & $f_7$ & $f_8$ & $f_{38}$ \\ 
	\noalign{\smallskip}\hline\noalign{\smallskip}
	$A_0$ & $5.13 \times 10^7$ & $0.1284$ & $1.117 \times 10^7$ & $8.878 \times 
	10^6$ & $1.508 \times 10^7$ & $1.236 \times 10^6$ & $4.47 \times 10^5$ \\
	$A_1$ & $-$ & $2.248 \times 10^6$ & $5.529 \times 10^6$ & $2.898 \times 10^5$ 
	& $3.305 \times 10^7$ & $0.9195$ & $1.404 \times 10^5$ \\
	$A_2$ & $-$ & $-$ & $-$ & $-$ & $2.656$ & $2.074$ & $1.442 \times 10^5$ \\
	$A_3$ & $-$ & $-$ & $-$ & $-$ & $5.05 \times 10^6$ & $6.497 \times 10^6$ & 
	$1.806 \times 10^5$ \\
	$A_4$ & $-$ & $-$ & $-$ & $-$ & $-$ & $0.1707$ & $9.587 \times 10^5$ \\
	$A_5$ & $-$ & $-$ & $-$ & $-$ & $-$ & $1.042$ & $2.08 \times 10^5$ \\
	$A_6$ & $-$ & $-$ & $-$ & $-$ & $-$ & $1.97 \times 10^5$ & $6396$ \\
	$A_7$ & $-$ & $-$ & $-$ & $-$ & $-$ & $2.26$ & $5117$ \\
	$A_8$ & $-$ & $-$ & $-$ & $-$ & $-$ & $7.68$ & $2.532 \times 10^5$ \\
	$A_9$ & $-$ & $-$ & $-$ & $-$ & $-$ & $-$ & $1.505 \times 10^4$ \\
		\noalign{\smallskip}\hline
		\end{tabular}
		\end{table}

	\begin{table}[!h]
		\centering
		\caption{Coefficients of the reduced order model for the maximum horizontal 
		stress at Point 1 (Equations 
		\eqref{eq:ROM_SigMax_P1_f0}--\eqref{eq:ROM_SigMax_P1_f38})}
		\label{App:Tab:ROM_Coeff_SigMax_P1}  
		\begin{tabular}{llllllll}
		\hline\noalign{\smallskip}
		Coeff & $f_0$ & $f_3$ & $f_5$ & $f_6$ & $f_7$ & $f_8$ & $f_{38}$\\ 
		\noalign{\smallskip}\hline\noalign{\smallskip}
		$A_0$ & $5.55 \times 10^7$ & $1$ & $8.78 \times 10^6$ & $2.933 \times 10^7$ 
		& $1.392 \times 10^7$ & $1.433 \times 10^6$ & $3.951 \times 10^7$ \\
		$A_1$ & $-$ & $1.604 \times 10^6$ & $4.468 \times 10^6$ & $4.687 \times 
		10^6$ & $2.93 \times 10^7$ & $0.8957$ & $2.278$  \\
		$A_2$ & $-$ & $-$ & $-$ & $9.515 \times 10^5$ & $2.233 
		\times 10^7$ & $1.726$ & $4.13 \times 10^6$ \\
		$A_3$ & $-$ & $-$ & $-$ & $-$ & $4.276 \times 10^6$ & $5.315 \times 10^6$ & 
		$0.2681$ \\
		$A_4$ & $-$ & $-$ & $-$ & $-$ & $-$ & $0.1887$ & $7.187 \times 10^4$ \\
		$A_5$ & $-$ & $-$ & $-$ & $-$ & $-$ & $0.804$ & $1797$ \\
		$A_6$ & $-$ & $-$ & $-$ & $-$ & $-$ & $2.289 \times 10^5$ & $-$ \\
		$A_7$ & $-$ & $-$ & $-$ & $-$ & $-$ & $2.185$ & $-$ \\
		$A_8$ & $-$ & $-$ & $-$ & $-$ & $-$ & $6.766$ & $-$ \\
			\noalign{\smallskip}\hline
		\end{tabular}
	\end{table}

	\begin{table}[!h]
	\centering
	\caption{Coefficients of the reduced order model for pore pressure 
	at Point 5 (Equations \eqref{eq:ROM_Pp_P5_f0}--\eqref{eq:ROM_Pp_P5_f38})}
	\label{App:Tab:ROM_Coeff_Pp_P5}  
	\begin{tabular}{llllll}
		\hline\noalign{\smallskip}
		Coeff & $f_0$ & $f_1$ & $f_3$ & $f_8$ & $f_{38}$ \\ 
		\noalign{\smallskip}\hline\noalign{\smallskip}
		$A_0$ & $3.23 \times 10^7$ & $3439$ & $0.5258$ & $2.939 \times 10^7$ & 
		$3.526 \times 10^6$ \\
		$A_1$ & $-$ & $5929$ & $9.196 \times 10^6$ & $0.107$ & $1.341 \times 10^6$ 
		\\
		$A_2$ & $-$ & $9.585 \times 10^6$ & $-$ & $1.871$ & $5.341 \times 10^5$ \\
		$A_3$ & $-$ & $-$ & $-$ & $2.294 \times 10^6$ & $2.164 \times 10^5$ \\
		$A_4$ & $-$ & $-$ & $-$ & $1.085$ & $2.553 \times 10^6$ \\
		$A_5$ & $-$ & $-$ & $-$ & $0.06293$ & $5.576 \times 10^5$ \\
		$A_6$ & $-$ & $-$ & $-$ & $4.435 \times 10^5$ & $3049$ \\
		$A_7$ & $-$ & $-$ & $-$ & $2.071$ & $1.674 \times 10^4$ \\
		$A_8$ & $-$ & $-$ & $-$ & $2.674$ & $7.557 \times 10^5$ \\
		$A_9$ & $-$ & $-$ & $-$ & $-$ & $1.106 \times 10^5$ \\
		\noalign{\smallskip}\hline
		\end{tabular}
		\end{table}

		\begin{table}[!h]
			\centering
			\caption{Coefficients of the first order Sobol functions for the 
			minimum horizontal stress at Point 5 (Equations 
			\eqref{eq:ROM_SigMin_P5_f0}--\eqref{eq:ROM_SigMin_P5_f8}).}
			\label{App:Tab:ROM_Coeff_SigMin_P5_a}  
			\begin{tabular}{llllllll}
				\hline\noalign{\smallskip}
				Coeff & $f_0$ & $f_1$ & $f_3$ & $f_5$ & $f_6$ & $f_7$ & $f_8$\\ 
				\noalign{\smallskip}\hline\noalign{\smallskip}
				$A_0$ & $4.75 \times 10^7$ & $0.1193$ & $0.2848$ & $2.227 \times 10^7$ & 
				$1.8 \times 10^7$ & $1.58 \times 10^7$ & $1.38 \times 10^6$\\
				$A_1$ & $-$ & $1.95 \times 10^6$ & $5.914 \times 10^6$ & $1.066 \times 
				10^7$ & $1.146 \times 10^6$ & $2.802 \times 10^7$ & $1.191$ \\
				$A_2$ & $-$ & $8.279 \times 10^5$ &  $-$ & $-$ & $-$ & $8.152 \times 10^6$ 
				& $1.187$ \\
				$A_3$ & $-$ & $4.604 \times 10^6$ & $-$ &  $-$ & $-$ & $-$ & $7.738 \times 
				10^6$ \\
				$A_4$ & $-$ & $9.287 \times 10^5$ & $-$ &  $-$ & $-$ & $-$ & $0.298$ \\
				$A_5$ & $-$ & $2.434 \times 10^6$ & $-$ &  $-$ & $-$ & $-$ & $0.6055$ \\
				$A_6$ & $-$ & $8.653 \times 10^5$ & $-$ &  $-$ & $-$ & $-$ & $-$ \\
				$A_7$ & $-$ & $1.009 \times 10^6$ & $-$ &  $-$ & $-$ & $-$ & $-$ \\
				$A_8$ & $-$ & $3.608 \times 10^5$ & $-$ &  $-$ & $-$ & $-$ & $-$ \\
				$A_9$ & $-$ & $2.868 \times 10^5$ & $-$ &  $-$ & $-$ & $-$ & $-$ \\
				\noalign{\smallskip}\hline
			\end{tabular}
		\end{table}

	\begin{table}[!h]
		\centering
		\caption{Coefficients of the second order Sobol functions for the 
			minimum horizontal stress at Point 5 (Equations 
			\eqref{eq:ROM_SigMin_P5_f38}--\eqref{eq:ROM_SigMin_P5_f18}).}
		\label{App:Tab:ROM_Coeff_SigMin_P5_b}  
		\begin{tabular}{llllll}
			\hline\noalign{\smallskip}
			Coeff & $f_{38}$ & $f_{13}$ & $f_{16}$ & $f_{17}$ & $f_{18}$ \\ 
			\noalign{\smallskip}\hline\noalign{\smallskip}
			$A_0$ & $9.637 \times 10^5$ & $1.123 \times 10^6$ & $1.06 \times 10^7$ 
			& $1.154 \times 10^6$ & $3.744 \times 10^6$ \\
			$A_1$ & $9.514 \times 10^4$ & $4.399 \times 10^5$ & $1.284 \times 10^6$ 
			& $7.664 \times 10^5$ & $1.954 \times 10^6$ \\
			$A_2$ & $4.621 \times 10^5$ & $1.114 \times 10^6$ & $9.594 \times 10^7$ 
			& $1.28 \times 10^6$ & $2.229 \times 10^6$ \\
			$A_3$ & $2.127 \times 10^5$ & $5.819 \times 10^5$ & $1.082 \times 10^4$ 
			& $3.21 \times 10^5$ & $1.827 \times 10^6$ \\
			$A_4$ & $1.347 \times 10^6$ & $1.24  \times 10^6$ & $7.181 \times 10^6$ 
			& $1.673 \times 10^6$ & $1.111 \times 10^5$ \\
			$A_5$ & $3.006 \times 10^5$ & $9.329 \times 10^5$ & $1.228 \times 10^9$ 
			& $1.64 \times 10^5$ & $1.796 \times 10^6$ \\
			$A_6$ & $1.411 \times 10^4$ & $9.442 \times 10^4$ & $99.63$ & $2.329 
			\times 10^5$ & $4.715 \times 10^5$ \\
			$A_7$ & $4099$              & $8.217 \times 10^5$ & $9.425 \times 10^4$ 
			& $9.159 \times 10^5$ & $1.76 \times 10^6$ \\
			$A_8$ & $1.667 \times 10^5$ & $1.604 \times 10^4$ & $3.635 \times 10^5$ 
			& $7.561  \times 10^4$ & $7.387 \times 10^5$ \\
			$A_9$ & $2.057 \times 10^5$ & $1.432 \times 10^5$ & $2.268 \times 10^9$ 
			& $1.482 \times 10^5$ & $1.013 \times 10^5$ \\
			$A_{10}$ & $-$ & $9.507 \times 10^4$ & $-$ & $1.409 \times 10^5$ & 
			$2.414 \times 10^5$ \\
			$A_{11}$ & $-$ & $6.673 \times 10^5$ & $-$ & $8.554 \times 10^5$ & 
			$9.83 \times 10^5$ \\
			$A_{12}$ & $-$ & $8912             $ & $-$ & $4.41  \times 10^5$ & 
			$7.301 \times 10^5$ \\
			$A_{13}$ & $-$ & $152.3            $ & $-$ & $4.697 \times 10^4$ & 
			$8.411 \times 10^5$ \\
			$A_{14}$ & $-$ & $2.301 \times 10^5$ & $-$ & $1.566 \times 10^5$ & 
			$3.61 \times 10^5$ \\
			\noalign{\smallskip}\hline
		\end{tabular}
	\end{table}

		\begin{table}[!h]
			\centering
			\caption{Coefficients of the reduced order model for the 
			maximum horizontal stress at Point 5 (Equations 
			\eqref{eq:ROM_SigMax_P5_f0}--\eqref{eq:ROM_SigMax_P5_f18})}
			\label{App:Tab:ROM_Coeff_SigMax_P5}  
			\begin{tabular}{llllllll}
				\hline\noalign{\smallskip}
				Coeff & $f_0$ & $f_3$ & $f_5$ & $f_6$ & $f_7$ & $f_8$ & $f_{18}$\\ 
				\noalign{\smallskip}\hline\noalign{\smallskip}
				$A_0$ & $5.05 \times 10^7$ & $0.339$ & $3.209 \times 10^7$ & $2.718 
				\times 10^7$ & $2.28 \times 10^7$ & $0.6016$ & $1.434 
				\times 10^6$ \\
				$A_1$ & $-$ & $8.31 \times 10^6$ & $1.293 \times 10^7$ & $4.38 
				\times 10^6$ & $6 \times 10^7$ & $1.648 \times 10^6$ & $1.274 
				\times 10^6$ \\
				$A_2$ & $-$ & $-$ & $-$ & $-$ & $5.831 \times 10^7$ & $3.453 \times 
				10^6$ & $2.579 \times 10^5$ \\
				$A_3$ & $-$ & $-$ & $-$ & $-$ & $1.648 \times 10^7$ & $2.723 
				\times 10^6$ & $7.189 \times 10^5$ \\
				$A_4$ & $-$ & $-$ & $-$ & $-$ & $-$ & $1.789 \times 10^4$ & 
				$2.978 \times 10^6$ \\
				$A_5$ & $-$ & $-$ & $-$ & $-$ & $-$ & $1.512 \times 10^6$ & $1.117 
				\times 10^6$ \\
				$A_6$ & $-$ & $-$ & $-$ & $-$ & $-$ & $-$ & $3.386 \times 10^4$ \\
				$A_7$ & $-$ & $-$ & $-$ & $-$ & $-$ & $-$ & $6.242 \times 10^4$ \\
				$A_8$ & $-$ & $-$ & $-$ & $-$ & $-$ & $-$ & $7.895 \times 10^5$ \\
				$A_9$ & $-$ & $-$ & $-$ & $-$ & $-$ & $-$ & $4.039 \times 10^4$ \\
				$A_{10}$ & $-$ & $-$ & $-$ & $-$ & $-$ & $-$ & $2.015 
				\times 10^5$\\
				$A_{11}$ & $-$ & $-$ & $-$ & $-$ & $-$ & $-$ & $4.165 \times 10^5$\\
				$A_{12}$ & $-$ & $-$ & $-$ & $-$ & $-$ & $-$ & $1.227 \times 10^5$\\
				$A_{13}$ & $-$ & $-$ & $-$ & $-$ & $-$ & $-$ & $1.066 \times 10^6$\\
				$A_{14}$ & $-$ & $-$ & $-$ & $-$ & $-$ & $-$ & $3.404 \times 10^5$ 
				\\
				\noalign{\smallskip}\hline
			\end{tabular}
		\end{table}

	\begin{table}[!h]
	\centering
	\caption{Coefficients of the reduced order model for the pore 
	pressure at Point 6 (Equations \eqref{eq:ROM_Pp_P6_f0} - 
	\eqref{eq:ROM_Pp_P6_f48})}
	\label{App:Tab:ROM_Coeff_Pp_P6}  
	\begin{tabular}{llllll}
		\hline\noalign{\smallskip}
		Coeff & $f_0$ & $f_1$ & $f_8$ & $f_{38}$	& $f_{48}$ \\ 
		\noalign{\smallskip}\hline\noalign{\smallskip}
		$A_0$ & $2.65 \times 10^7$ & $0.824$ & $1.44 \times 10^8$ & $3.644 \times 
		10^6$ & $3.201 \times 10^6$ \\
		$A_1$ & $-$ & $1.789 \times 10^7$ & $0.4676$ & $2.137 \times 10^6$ & $8.985  
		\times 10^5$ \\
		$A_2$ & $-$ & $-$ & $2.499$ & $1.057 \times 10^6$ & $6.924 \times 10^5$ \\
		$A_3$ & $-$ & $-$ & $1.352 \times 10^8$ & $1.295 \times 10^5$ & $1.588 \times 
		10^6$ \\
		$A_4$ & $-$ & $-$ & $0.5053$ & $2.476 \times 10^6$ & $2.344 \times 10^6$ \\
		$A_5$ & $-$ & $-$ & $6.469$ & $7.658 \times 10^5$ & $3.294 \times 10^6$ \\
		$A_6$ & $-$ & $-$ & $-$ & $6.063 \times 10^4$ & $9.518 \times 10^5$ \\
		$A_7$ & $-$ & $-$ & $-$ & $1060$ & $1.285 \times 10^6$ \\
		$A_8$ & $-$ & $-$ & $-$ & $1.79 \times 10^6$ & $4.197 \times 10^6$ \\
		$A_9$ & $-$ & $-$ & $-$ & $3.814 \times 10^5$ & $3.788 \times 10^4$ \\
		$A_{10}$ & $-$ & $-$ & $-$ & $-$ & $8.773 \times 10^5$ \\
		$A_{11}$ & $-$ & $-$ & $-$ & $-$ & $2.672 \times 10^4$ \\
		$A_{12}$ & $-$ & $-$ & $-$ & $-$ & $1.711 \times 10^5$ \\
		$A_{13}$ & $-$ & $-$ & $-$ & $-$ & $1.057 \times 10^6$ \\
		$A_{14}$ & $-$ & $-$ & $-$ & $-$ & $9.727 \times 10^5$ \\
		$A_{15}$ & $-$ & $-$ & $-$ & $-$ & $3.976 \times 10^5$ \\
		$A_{16}$ & $-$ & $-$ & $-$ & $-$ & $2.538 \times 10^4$ \\
		$A_{17}$ & $-$ & $-$ & $-$ & $-$ & $1.006 \times 10^5$ \\
		$A_{18}$ & $-$ & $-$ & $-$ & $-$ & $7.33  \times 10^5$ \\
		$A_{19}$ & $-$ & $-$ & $-$ & $-$ & $1.885 \times 10^6$ \\
		$A_{20}$ & $-$ & $-$ & $-$ & $-$ & $2.617 \times 10^5$ \\
		\noalign{\smallskip}\hline
	\end{tabular}
	\end{table}

	\begin{table}[!h]
	\centering
	\caption{Coefficients of the reduced order model for the minimum 
	horizontal 
	stress at Point 6 (Equations 
	\eqref{eq:ROM_SigMin_P6_f0}--\eqref{eq:ROM_SigMin_P6_f12})}
	\label{App:Tab:ROM_Coeff_SigMin_P6}  
	\begin{tabular}{llll}
	\hline\noalign{\smallskip}
	Coeff & $f_0$ & $f_1$ & $f_{12}$ \\ 
	\noalign{\smallskip}\hline\noalign{\smallskip}
	$A_0$ & $5.53 \times 10^7$ & $5806$ & $3.886 \times 10^4$ \\
	$A_1$ & $-$ & $4.612 \times 10^4$ & $3.676 \times 10^4$ \\
	$A_2$ & $-$ & $1.266 \times 10^5$ & $5.765 \times 10^4$ \\
	$A_3$ & $-$ & $1.355 \times 10^5$ & $3.558 \times 10^4$ \\
	$A_4$ & $-$ & $4994$ & $1.259 \times 10^5$ \\
	$A_5$ & $-$ & $1.696 \times 10^5$ & $1.47 \times 10^4$ \\
	$A_6$ & $-$ & $2.416 \times 10^5$ & $2.868 \times 10^4$ \\
	$A_7$ & $-$ & $1.984 \times 10^5$ & $5.699 \times 10^4$ \\
	$A_8$ & $-$ & $3.771 \times 10^4$ & $1.4 \times 10^4$ \\
	$A_9$ & $-$ & $-$ & $3266$ \\
	\noalign{\smallskip}\hline
	\end{tabular}
	\end{table}

	\begin{table}[!h]
	\centering
	\caption{Coefficients of the reduced order model for the maximum 
	horizontal 
	stress at Point 6 (Equations 
	\eqref{eq:ROM_SigMax_P6_f0}--\eqref{eq:ROM_SigMax_P6_f18})}
	\label{App:Tab:ROM_Coeff_SigMax_P6}  
	\begin{tabular}{lllllllll}
	\hline\noalign{\smallskip}
	Coeff & $f_0$ & $f_1$ & $f_3$ & $f_5$	& $f_6$ & $f_7$ & $f_8$ & $f_{18}$ \\ 
	\noalign{\smallskip}\hline\noalign{\smallskip}
	$A_0$ & $5.05 \times 10^7$ & $1708$ & $0.335$ & $3.19 \times 10^7$ & $2.72 
	\times 10^7$ & $2.25 \times 10^7$& $2.565 \times 10^7$ & $2.528 \times 10^6$	\\
	$A_1$ & $-$ & $1.016 \times 10^5$ & $8.67 \times 10^6$ & $1.2 \times 10^7$ & 
	$5.1 \times 10^6$ & $2.565 \times 10^7$ & $1.271$ &	$3.111 \times 10^5$ \\
	$A_2$ & $-$ & $1.436 \times 10^6$ & $-$ & $-$ & $-$ & $1.46 \times 10^7$ & 
	$3.662$ & $3.678 
	\times 10^4$ \\
	$A_3$ & $-$ & $-$ & $-$ & $-$ & $-$ & $-$ & $1.506 \times 10^7$ & $4.567 
	\times 10^4$ \\
	$A_4$ & $-$ & $-$ & $-$ & $-$ & $-$ & $-$ & $1.411$ & $6.293 \times 10^5$ \\
	$A_5$ & $-$ & $-$ & $-$ & $-$ & $-$ & $-$ & $4.119$ & $3.488 \times 10^5$ \\
	$A_6$ & $-$ & $-$ & $-$ & $-$ & $-$ & $-$ & $1.153 \times 10^7$ & $6.563 
	\times 10^4$ \\
	$A_7$ & $-$ & $-$ & $-$ & $-$ & $-$ & $-$ & $1.016$ & $1.54 \times 10^5$ \\
	$A_8$ & $-$ & $-$ & $-$ & $-$ & $-$ & $-$ & $10.12$ & $1.734 \times 10^5$ \\
	$A_9$ & $-$ & $-$ & $-$ & $-$ & $-$ & $-$ & $-$ & $4.125 \times 10^4$ \\
	\noalign{\smallskip}\hline
	\end{tabular}
	\end{table}

\end{document}